\documentclass[11pt,a4paper]{amsart} 

\usepackage{amsmath}
\usepackage{amssymb}
\usepackage{enumerate}
\usepackage{appendix}

\usepackage[pdfborder={0 0 0},pagebackref,colorlinks,citecolor=darkgreen,linkcolor=darkred,urlcolor=darkblue]{hyperref}

\usepackage{multicol}
\usepackage{mathtools}
\usepackage[capitalise]{cleveref}
\usepackage{mathrsfs}
\usepackage{tikz}
\usepackage{graphicx}
\usepackage{tabularx}
\usepackage{ tipa }
\usepackage{a4wide}
\usetikzlibrary{calc}
\usepackage{tikz-cd}
\usepackage{adjustbox}

\hypersetup{
	colorlinks=true,
	linkcolor=blue,
	citecolor=magenta,      
	urlcolor=cyan,
	pdftitle={TITLE},
}
\usepackage{verbatim}
\usepackage[normalem]{ulem}
\usepackage{ dsfont }
\usepackage{ stmaryrd }

\newtheorem{mnthm}{Theorem}[section]

\newtheorem{thm}{Theorem}[section]
\newtheorem{prop}[thm]{Proposition}
\newtheorem{theorem}[thm]{Theorem}
\newtheorem{lem}[thm]{Lemma}

\newtheorem{corollary}[thm]{Corollary}
 
\newtheorem{setup}[thm]{Setup}

\theoremstyle{definition}
\newtheorem{definition}[thm]{Definition}

\newtheorem{lemma}[thm]{Lemma}
\newtheorem{proposition}[thm]{Proposition}

\theoremstyle{remark}
\newtheorem{construction}[thm]{Construction}
\newtheorem{remark}[thm]{Remark}

\newtheorem{clm}[thm]{Claim}

\numberwithin{equation}{section}

\newcommand{\Z}{\mathbb{Z}}

\newcommand{\CP}{\mathbb{CP}}
\newcommand{\RP}{\mathbb{RP}}

\newcommand{\ol}{\overline}

\newcommand{\wt}{\widetilde}

\newcommand{\mudash}{\smash{{\mkern0.75mu\raisebox{-1.52ex}{$\mathchar '26$}\mkern -5.8mu\mu}}}
\newcommand{\lambdadash}{{\mkern0.75mu\mathchar '26\mkern -9.75mu\lambda}}
\newcommand\smallO{
	\mathchoice
	{{\scriptstyle\mathcal{O}}}
	{{\scriptstyle\mathcal{O}}}
	{{\scriptscriptstyle\mathcal{O}}}
	{\scalebox{.7}{$\scriptscriptstyle\mathcal{O}$}}
}
\newcommand{\mudashj}{\mudash_{\raisemath{6pt}{J}}}
\newcommand{\mudashJ}{\mudash_{\raisemath{6pt}{J}}}

\makeatletter
\newcommand{\raisemath}[1]{\mathpalette{\raisem@th{#1}}}
\newcommand{\raisem@th}[3]{\raisebox{#1}{$#2#3$}}
\makeatother

\newcommand{\xtwoheadrightarrow}[2][]{%
	\xrightarrow[#1]{#2}\mathrel{\mkern-14mu}\rightarrow
}

\usetikzlibrary{decorations.pathmorphing,shapes}

\newcounter{sarrow}

\newcommand\setItemnumber[1]{\setcounter{enumi}{\numexpr#1-1\relax}}

\DeclareMathOperator\PD{PD}
\DeclareMathOperator\proj{proj}
\DeclareMathOperator\Id{Id}

\DeclareMathOperator\pt{pt}
	\DeclareMathOperator\Spin{Spin}

    \DeclareMathOperator\Top{Top}

	\DeclareMathOperator\tO{O}
	
	\DeclareMathOperator\BSpin{B\;\Spin}

	\DeclareMathOperator\BO{B\;\tO}
	
	\DeclareMathOperator\Pin{Pin}

	\DeclareMathOperator\Hom{Hom}

	\DeclareMathOperator\secondary{\mathfrak{sec}}
	\DeclareMathOperator\pri{\mathfrak{pri}}
	\DeclareMathOperator\ter{\mathfrak{ter}}
	\DeclareMathOperator\secG{\mathfrak{sec}^{geo}}
	
	\DeclareMathOperator\terG{\mathfrak{ter}^{geo}}
	\DeclareMathOperator\Sq{Sq}
	\DeclareMathOperator\im{im}
	\DeclareMathOperator\Tor{Tor}
	\DeclareMathOperator\Th{Th}
	\DeclareMathOperator\km{km}
	\DeclareMathOperator\colim{colim}
	\DeclareMathOperator\pr{pr}
	\DeclareMathOperator\red{red}
	\DeclareMathOperator\interior{int}
	\DeclareMathOperator\id{id}

	\DeclareMathOperator\Aut{Aut}

	\DeclareMathOperator\ev{ev}

	\title[$\xi$-fillings of 3-manifolds]{Geometric obstructions for $\xi$-fillings of 3-manifolds}
	
	\author[D. Galvin]{Daniel Galvin}
	\address{University of Texas at Austin, US}
	\email{daniel.galvin@austin.utexas.edu}
	\author[P. Teichner]{Peter Teichner}
	\address{Max Planck Institute for Mathematics, Bonn, Germany}
	\email{teichner@mpim-bonn.mpg.de}
	\author[S. Vesel\'{a}]{Simona Vesel\'{a}}
	\address{Max Planck Institute for Mathematics, Bonn, Germany}
	\email{vesela@mpim-bonn.mpg.de}

\begin{document}
		\begin{abstract}
			We consider the realisation problem for normal 1-types of 4-manifolds with a given boundary.  More precisely, given a normal 1-type $\xi$ and closed 3-dimensional $\xi$-manifold $Y$, does there exist a compact 4-dimensional $\xi$-manifold with boundary $Y$?  We describe a three stage obstruction theory for the existence of such a 4-manifold, with our main contribution being a `tertiary' obstruction that we describe geometrically via Wall's quadratic self-intersection form.
		\end{abstract}		
		\maketitle

		\section{Introduction}\label{sec:intro}

		An application of Kreck's modified surgery theory \cite{Kreck99} reduces the problem of classifying compact, smooth 4-manifolds up to $(S^2 \times  S^2)$-stable diffeomorphism to computing certain $\xi$-bordism groups. Here $\xi:B\to BO$ is the \emph{normal 1-type} of the manifolds $X$ in question,  a Postnikov-Moore approximation of the stable normal bundle $\nu_X:X\to BO$. So $\xi$ induces an isomorphism on homotopy groups $\pi_i$ for $i>2$ (and a monomorphism on $\pi_2$) and there is a lift $u_X: X \to B$ of $\nu_X$ that is an isomorphism on $\pi_i$ for $i<2$ (and an epimorphism on $\pi_2$).  
		
		We call $u_X$ a \emph{normal 1-smoothing} on $X$ and whenever $u_M\colon M \to B$ lifts the stable normal bundle of a manifold $M$ over $\xi$, we call $(M,u_M)$ a $\xi$-manifold and $u_M$ a $\xi$-structure. If $Y$ is a closed 3-manifold with $\xi$-structure $u_Y$, one can ask the following existence and uniqueness questions:
		\begin{enumerate}
			\item Does $(Y,u_Y)$ have a $\xi$-filling, i.e.\ a 4-dimensional $\xi$-manifold $(X,u_X)$ with boundary $(Y,u_Y)$?
			\item Can we choose a $\xi$-filling $(X,u_X)$ to be a normal 1-smoothing?
			\item How many such 1-smoothings $(X,u_X)$ exist up to stable diffeomorphism?
		\end{enumerate}
		
		By doing surgeries in the interior of $X$ $\xi$-fillings can be turned into  normal 1-smoothings (\Cref{lem:surgeryBelowMiddle}), answering (2) modulo (1). Kreck's theorem turns (3) into a question about 5-dimensional $\xi$-fillings, as follows. Two normal 1-smoothings $(X_i,u_i)$ with boundary $(Y,u_Y)$ are stably diffeomorphic (with boundary diffeomorphism $\varphi:Y\to Y$) if and only if $X_1 \cup_{\varphi} -X_2$ admits a 5-dimensional $\xi$-filling. 
		As a consequence,  (1) is a lower dimensional cousin of the stable classification question (3), addressing the \emph{existence} question for normal 1-smoothings $(X,u_X)$ with given boundary condition. It is the purpose of this paper. 
		
		We start with a $\xi$-manifold $(Y,u_Y)$ as above and ask whether it represents zero in the bordism group $\Omega_3^{\xi}$ of $\xi$-manifolds. One case is very easy, namely when $\pi_2(B)\neq 0$ (the \emph{totally non-spin} case) and hence $\pi_2(\xi)$ gives an isomorphism $\pi_2(B) \cong \pi_2(BO) \cong \Z/2$. Then $B$ is  homotopy equivalent to $BSO \times B\pi$ and $\Omega_3^{\xi}\cong H_3(\pi;\Z^{w_1^\pi})$. Here $\pi:=\pi_1(B)$ and $w_1^\pi\in H^1(\pi;\Z/2)$ is the first Stiefel-Whitney class of $\xi$ restricted to $B\pi$. In this case, $(Y,u_Y)$ has a $\xi$-filling if and only if $(p_2\circ u_Y)_*[Y]=0\in H_3(\pi;\Z^{w_1^\pi})$, giving a complete answer to (1).
		
		From now on we assume that $\pi_2(B)=0$, so the normal 1-type $\xi=\xi(\pi,w_1^{\pi},w_2^{\pi})\colon B\to BO$ is the pullback of the universal Stiefel-Whitney classes $(w_1, w_2):BO \to B\Z/2 \times B^2\Z/2$ and $(w_1^{\pi}, w_2^{\pi}):B\pi \to B\Z/2 \times B^2\Z/2$. In particular,
		there is a fibration $p:B \to B\pi$ such that the restriction of $\xi$ to a fibre $BSpin$ is the canonical projection.
		
		In this situation,  the second-named author \cite{teichnerthesis} (cf.\ \Cref{thm:jss}) constructed what he called the ``James spectral sequence'' for $\xi$. The skeleta of $B\pi$ induce a filtration 
		\[
		\{0\}=F_{0,3}\leq F_{1,2}\leq F_{2,1} \leq F_{3,0}=\Omega_3^{\xi},
		\]
		for which we will give new geometric interpretations.  
		
		Any $\xi$-structure on $M$ gives $p \circ u_M\colon M\to B\pi $ over $w_1^{\pi}$.  If $(Y,p\circ u_{Y})$ bounds a compact 4-manifold $(X,u_X^{\pi})$ with $u_X^{\pi}\colon X\to B\pi$ an extension of $p\circ u_Y$ then we call it a $(B\pi,w_1^{\pi})$-filling of $(Y,u_Y)$.
		There is a single obstruction $\smallO(X,Y,u_X^{\pi},u_Y)\in H^2(X,Y;\Z/2)$, to extending the  $\xi$-structure $u_Y$ over $X$, relative to $u_X^{\pi}$.

		\begin{mnthm}\label{mnthm:1}
			Let $\xi(\pi,w_1^{\pi},w_2^{\pi})\colon B\to BO$ be a normal 1-type with $\pi_2(B)=0$ and let $(Y,u_Y)$ be a 3-dimensional $\xi$-manifold.  For the above filtration of $\Omega_3^{\xi}$ the corresponding geometric conditions are as follows:
			
			\begin{center}
				\renewcommand{\arraystretch}{1.7}
				\setlength{\tabcolsep}{4pt}
				\begin{tabular}{rcp{9.4cm}}
					$(Y,u_Y)\in F_{2,1}$ & $\Longleftrightarrow$ & There exists a $(B\pi,w_1^{\pi})$-filling  of $(Y,u_Y)$. \\ 
					$(Y,u_Y)\in F_{1,2}$ & $\Longleftrightarrow$ & There exists a $(B\pi,w_1^{\pi})$-filling $(X,u_X^{\pi})$ of $(Y,u_Y)$, such that $\smallO(X,Y,u_X^{\pi},u_Y)$ has a spherical representative. \\  
					$(Y,u_Y)\in F_{0,3}$ & $\Longleftrightarrow$ & There exists a $\xi$-filling of $(Y,u_Y)$. Moreover,
				\end{tabular}
			\end{center}
			\begin{center}
				\renewcommand{\arraystretch}{1.7}
				\setlength{\tabcolsep}{4pt}
				\begin{tabular}{rp{8.6cm}cr}
					 $ (Y,u_Y)\in F_{2,1}\Longleftrightarrow$ & $(p\circ u_Y)_*[Y]=0\in H_3(\pi;\Z^{w_1^{\pi}})$. \\ 
					If $(Y,u_Y)\in F_{2,1}$, then $(Y,u_Y)\in F_{1,2}  \Longleftrightarrow $ & $(u_X^{\pi})_*\PD(\smallO(X,Y,u_X^{\pi},u_Y))=0\in H_2(\pi;\Z/2)/(\delta_2)$. \\  
					If $(Y,u_Y)\in F_{1,2}$, then $(Y,u_Y)\in F_{0,3} \Longleftrightarrow$ & For a spherical representative $c$ of $\smallO(X,Y,u_{X}^{\pi},u_Y)$ we have $[\mudashj(c)]=0\in H_1(\pi;\Z/2)/(\delta_2,\delta_3)$.    
				\end{tabular}
			\end{center}
			where $\mudashj$ denotes a certain modification of Wall's self-intersection number for spheres, which we introduce in \Cref{sec:intersections}, and $\delta_2$ and $\delta_3$ correspond to certain `geometric differentials' which we will define later.

		\end{mnthm}
		
		Wall's quadratic self intersection form $\mu\colon \pi_2(X)\to \Z[\pi]/(\sim)$ counts self-intersections (with fundamental group elements) of immersed spheres in a 4-manifold $X$, and is valued in a certain quotient of the group ring $\Z[\pi]$ (which we will make precise in \Cref{sec:intersections}).  We define an adjustment of this $\mudash\colon \pi_2(X)\to I/(\sim)$, where $I:=\ker(\Z[\pi]\to\Z)$ denotes the augmentation ideal.  We then define $\mudashj$ by postcomposing with a projection map $I/(\sim) \to H_1(\pi;\Z/2)$, which we will describe in \Cref{sec:intersections}.  To define the tertiary invariant in terms of the self-intersection of the class $c$ (as in the statement of \Cref{mnthm:1}) we need to make it not depend on various choices made therein. This quotient is in terms of certain `geometric differentials' $\delta_*$ that we now describe briefly.  They will mirror the differentials in the James spectral sequence, which we denote by $d_*$.  We define the first as 
		\begin{align*}
			\delta^{(4,0)}_2\colon H_4(\pi;\Z^{w_1^{\pi}})&\to H_2(\pi;\Z/2) \\
			(u_M^{\pi})_*[M] &\mapsto (u_M^{\pi})_*\PD(\smallO(M,\emptyset,u_M^{\pi},\emptyset))
		\end{align*}
		where $(M,u_M^{\pi})$ is a closed 4-manifold representing $a$ (implicitly this means that the image of the obstruction does not depend on this choice).  Again, let $(X,u_X^{\pi})$ be a filling of $(Y,u_Y)$ such that $u_X^{\pi}$ induces an isomorphism on fundamental groups (this can be achieved by surgery).  Then there is the following exact sequence corresponding to the fibration $u_X^{\pi}$ (see \Cref{sec:SerreExSeq}).
		\begin{equation}\label{eq:serre_spectral_sequence_induced_exact_sequence}
			H_3(X;Z/2)\xrightarrow{(u_{X}^{\pi})_*} H_3(\pi;\Z/2)\xrightarrow{\mathfrak{f}} \pi_2(X)\otimes \Z/2 \xrightarrow{h_2^{\Lambda}} H_2(X;\Z/2)\xrightarrow{(u_{X}^{\pi})_*} H_2(\pi;\Z/2)
		\end{equation}
		Hence, for $z\in H_3(\pi;\Z/2)$, choosing a lift for $\mathfrak{f}(z)$ produces an immersed sphere which is null-$\Z/2$-homologous (in fact, all such spheres are produced in this way).  We define the second geometric differential as
		\begin{align*}
			\delta^{(3,1)}_2\colon H_3(\pi;\Z/2)&\to H_1(\pi;\Z/2) \\
			z &\mapsto \mudashj(c) + \lambdadash(c,w)
		\end{align*}
		where $w$ is an immersed sphere representative for $\smallO(X',Y,u_{X'}^{\pi},u_Y)$, $c$ denotes a representative immersed sphere for $\mathfrak{f}(z)$, and $\lambdadash$ denotes a suitable normalisation of the equivariant intersection form on $X$ (in line with how $\mudashj$ is defined).  We show this is well-defined in \cref{subsec:delta_d2}.
		
		\begin{mnthm}\label{mnthm:2}
			Let $\xi(\pi,w_1^{\pi},w_2^{\pi})$ be a normal 1-type.  Then $\delta^{(4,0)}_2=d_2^{(4,0)}$, where $d_2^{(4,0)}$ is the corresponding differential in the James spectral sequence, and
			$\delta^{(3,1)}_2=d_2^{(3,1)}$, where $d_2^{(3,1)}$ is the corresponding differential in the James spectral sequence.
		\end{mnthm}
		
		The geometric differential $\delta_3^{(4,0)}$ is defined as
		\begin{align*}
			\delta_3^{(4,0)}\colon \ker(\delta_2^{(4,0)})\to& H_1(\pi;\Z/2)/(\delta_2^{(3,1)}) \\
			(u_M^{\pi})_*[M] \mapsto& \mudash_J(c_M)
		\end{align*}
		where $(M,u_M^{\pi})$ is a closed 4-manifold with $\delta_3^{(4,0)}(u_M^{\pi})_*[M])=0$ and $c_M\in\pi_2(M)$ is a spherical representative of $ \smallO(M,u_M^{\pi})$.  We show that this is well-defined in \Cref{sbs:delta_3_differential}.
		
		\begin{remark}
			It is a consequence of \Cref{mnthm:1} that the image of $\delta_3^{(4,0)}$ coincides with the image of $d_3^{(4,0)}$ but we do not know if they coincide as maps.  This fact is stated later as \Cref{prop:d3delta3Agree}.
		\end{remark}
		
		A well-known theorem due to Wu (see e.g.\ \cite[Theorem 11.14]{milnorstasheff}) applied to 4-manifolds gives the formula
		\[
		Q(\alpha, w) = Q(\alpha, \alpha) \mod{2},
		\]
		where $Q$ denotes the intersection form on homology of a closed 4-manifold $X$, $w$ denotes the Poincar\'{e} dual of the second Stiefel-Whitney class $w_2(X)$, and $\alpha\in H_2(X)$.  A corollary of \Cref{mnthm:2} is the following secondary version of this Wu formula, which may be of independent interest.
		
		\begin{mnthm}\label{mnthm:wu}$[$Secondary Wu formula$]$
			Let $X$ be a compact 4-manifold with boundary $Y$.  For any choice of normal spin structure $s$ on $Y$, let $w_2^{s}\in H^2(X,Y;\Z/2)$ denote the relative second Stiefel-Whitney class. Let $c$ be an immersed spherical representative for $w_2^{s}$, and let $b$ be any null-$\Z/2$-homologous immersed sphere in $X$. Fix any map $u_X^{\pi}\colon X\to B\pi$.  Then in $I$ we have
			\[
			\lambdadash_J(b,c) = \mudash_J(b) \mod{IJ}
			\]
			where $I$ and $J$ denote the $\Z$ and $\Z/2$ augmentation ideals of $\Z[\pi_1(X)]$, respectively and $\lambdadash$, $\mudash$ are defined in Section \ref{sec:intersections}.
		\end{mnthm}
		
		Interestingly, \Cref{mnthm:wu} applies outside of the relative setting that the rest of the paper is in.

        \begin{remark}\label{rmk:quadratic2types}
            All of the data that is required to state the Secondary Wu formula can be obtained from the \emph{quadratic 2-type} $(\pi_1(X),\pi_2(X),k(X),\lambda(X))$ of a spin 4-manifold $X$, where $k(X)$ denotes the $k$-invariant corresponding to the second Postnikov stage of $X$.  This is because, in the orientable case, $\mudash_J$ is determined by $\lambdadash_J$; by the Serre exact sequence (\ref{eq:serre_spectral_sequence_induced_exact_sequence}), being a $\Z/2$-null-homologous immersed sphere in $X$ is precisely determined by being in the image of $\mathfrak{f}$, the map from the exact sequence (\ref{eq:serre_spectral_sequence_induced_exact_sequence}), which is determined by $k(X)$ (see \Cref{sec:SerreExSeq}); and $c$ can be taken to be the unknot, since $X$ is spin.  This means that \Cref{mnthm:wu} could potentially be used to prove non-realisability of certain quadratic 2-types by spin 4-manifolds.  We have not attempted to explore this avenue here.
        \end{remark}

		\begin{remark}
			The reader who in interested in an example of a calculation of $\mudashj(c)$ (as in \Cref{mnthm:1}) can look at \cref{prop:realisation_for_terG}, where we calculate the prototypical non-trivial example of this obstruction.  This is for $\xi=\xi(\Z,w_1^\pi,0)$ and the $3$-manifold is either $Y=T^3$ or $S^1\times K$ for $K$ a Klein bottle, depending on $w_1^{\pi}$.
		\end{remark}
		
		\begin{remark}\label{rmk:smooth_vs_top}
			In this paper we work entirely in the smooth category, since working in the topological category provides no extra generality.  We elaborate on this now.  
			
			Firstly, every topological 3-manifold admits a smooth structure unique up to isotopy \cite{cairns_1940, moise_1952, bing_1954}, so we can freely assume that our input 3-manifolds are smooth.  Secondly, assume that we have a smooth $\xi$-3-manifold $(Y,u_Y)$ which admits a \emph{topological} $\xi$-filling $(X,u_X)$, where $\xi$ is a potential normal 1-type for a 4-manifold.  The relative Kirby-Siebenmann invariant $ks(X,Y)\in H_4(X,Y;\Z/2)\cong \Z/2$ can be modified to be zero by connected-summing with Freedman's $E_8$ manifold \cite{Freedman82}, if necessary, and since $E_8$ is a spin, simply-connected $4$-manifold the $\xi$ structure $u_X$ may be extended across the connected-sum (note that this is where we use that $\xi$ is a normal 1-type).  Then, by Freedman-Quinn's stable smoothing theorem \cite[Theorem 8.6]{Freedman-Quinn}, there exists a smooth structure on our filling after taking sufficiently many connected-sums with $S^2\times S^2$.  Again, since $S^2\times S^2$ is a spin, simply-connected 4-manifold, there is no issue in extending the $\xi$-structure and hence we have produced a \emph{smooth} $\xi$-filling of $(Y,u_Y)$.
		\end{remark}

		\subsection{Background}  Milnor \cite{Milnor-Spin} and Kirby \cite[VII, Theorem 3]{Kirby-4-manifold-book} proved that the spin bordism group $\Omega_3^{\Spin}$ is zero. In the proof, one actually produces a simply-connected filling.  In Kirby's proof, he starts by showing that there exists an oriented filling, and he then considers a representative for the second relative Stiefel-Whitney class.  He then uses an explicit geometric argument to complete the spin null-bordism.  One can interpret \Cref{mnthm:1} to be a generalisation of this result for all normal 1-types, though if one specialises our proof to the simply-connected case one obtains a slightly different proof of $\Omega_3^{\Spin}=0$ (the interested reader can follow through the proof of \Cref{mnthm:1}, taking the normal 1-type to be $\xi(\pi=\{1\},0,0)$, and see the differences with Kirby's proof).
		
		In \cite{KPT,KPT25}, the authors work on a similar problem to this paper but for $\Omega_4^{\xi}$, in the special case that $\xi=\xi(\pi,0,0)$.  Those works concern the uniqueness problem, i.e.\ whether two smooth 4-manifolds with normal 1-type $\xi$ are stably diffeomorphic, rather than the existence problem studied in this paper.  In that case, there are also primary, secondary and tertiary invariants, and the authors relate them to algebraic invariants in certain cases.  A particularly interesting case is discussed in \cite{KPT25}, where the authors relate the 4-dimensional tertiary invariant to the Kervaire-Milnor invariant, a secondary intersection invariant (see \Cref{sec:km}).  For us, our 3-dimensional tertiary invariant is a primary intersection invariant.  We do not have a theoretical justification for why this shift occurs, but it would be very interesting to know if there is some general phenomenon at play.
		
		In \cite{kasprowski_nicholson_vesela_24}, Kasprowski, Nicholson and the third-named author discuss geometric interpretations for each stage of the filtration of $\Omega_4^{\xi}$ given by the JSS, much like how we have stated \Cref{mnthm:1}.  Though again, much like in \cite{KPT} the emphasis is on the equivalence relation on the 4-manifolds involved, rather than on the 5-dimensional $\xi$-bordism itself.  One might try to also consider the filtration stages of $\Omega_3^{\xi}$ as giving equivalence relations on polarised 3-manifolds, but studying 3-manifolds up to these relations is potentially wrongheaded, since 3-manifolds are more naturally classified by geometric techniques developed separately to the high-dimensional theory.

		\subsection{Conventions and notation}\label{subsec:Conventions}
		
		\begin{enumerate}		
			\item By JSS we mean the {\it James spectral sequence}, introduced in \cite{teichnerthesis}, cf.\ \Cref{sec:jss}.
			
			\item We will be using a $B\Pin^{+}$-structure on the normal bundles of manifolds. Usually, manifolds are said to be $\Pin^{+}$ (resp.\ $\Pin^-$) if their {\it tangential} Stiefel-Whitney classes satisfy $w^t_2=0$ (resp.\ $w^{t}_2=(w_1^t)^2$). We deal with normal bundles and so say that a manifold has a $\Pin^{+}$ (resp.\ $\Pin^-$) structure if the Stiefel-Whitney classes of the normal bundle satisfy $w_2=0$ (resp.\ $w_2=w_1^{2}$). Note that $w_1=w_1^t$ and $w_2^{t}=w_1^2+w_2$ and that the Wu formula for normal bundle Stiefel-Whitney classes is
			\[
			x^2= w_2(M)\smile x\in H^2(M;\Z/2),
			\]
			for $M$ a possibly non-orientable 4-manifold.
				\item Let $A$ be a finitely presented abelian group. We write $B^n A$ for the $n$-th iterated bar construction.  This gives a model for the Eilenberg-MacLane space $K(A,n)$.
				\item Let $X$ be a compact space, let $\pi:\pi_1(X)$ be a group and set $\Lambda=\Z[\pi]$. We will use  Hurewicz homomorphisms $h_2,h_2^{\Z},h_2^{\Lambda}$ with domains $\pi_2(X), \pi_2(X)\otimes_\Z \Z/2, \pi_2(X)\otimes_{\Lambda} \Z/2$ respectively.

				\end{enumerate}
				
				\subsection{Organisation}
				In \Cref{sec:surgery} we recall the various notions of Kreck's modified surgery theory that we require.  In \Cref{sec:intersections} we define and discuss the modifications of the standard intersection invariants used in surgery theory, as well as recalling the Kervaire-Milnor invariant.  In \Cref{sec:obstructions} we use the modified intersection invariants to define our null-bordism obstructions.  In \Cref{sec:differentials} we show that these obstructions are well-defined by working with the geometric differentials.  In particular, we prove \Cref{mnthm:2} and its corollary \Cref{mnthm:wu}.  Finally, in \Cref{sec:vanishing} we show that our obstructions are complete obstructions to the null-bordism question, i.e.\ we prove \Cref{mnthm:1}. In \Cref{sec:SerreExSeq} we show a technical statement about a certain Serre exact sequence.   In \Cref{sec:jss} we extend the definition of the JSS to the non-orientable setting; this may be of independent interest.
				
				\subsection{Acknowledgements}
				We thank Daniel Kasprowski for helpful discussions.  The authors would like to thank the \emph{Max Plank Institute for Mathematics in Bonn} for their hospitality and support throughout much of the work.

\newpage

\begin{multicols}{2}
                \tableofcontents
                \end{multicols}
                

				
                        

				\section{Modified Surgery}\label{sec:surgery}

				\subsection{Normal \texorpdfstring{$1$}{1}--types}
				
				Consider the space $BO$, which classifies stable vector bundles. Kreck’s modified surgery theory \cite{Kreck99} provides a stable classification of manifolds by introducing certain approximation spaces $B$, which are equipped with a highly coconnected map to $BO$. The space $B$ serves as an approximation to a given manifold $M$ and is defined as a Moore–Postnikov approximation for the stable normal bundle $X\xrightarrow{\nu_X} BO$. We now give a brief overview of Kreck's theory.  Given a map $\xi\colon B\to BO$ a $\xi$-structure on a smooth manifold $X$ is a map $\ol{\nu}_X\colon X\to B$, which lifts the normal bundle $\nu_X$ over $\xi$. 
				\[
				\begin{tikzcd}
					& B \arrow[d,"\xi"] \\
					X \arrow[ur,"\ol{\nu}_X"] \arrow[r, ""{name=U, below}]{} \arrow[r,"\nu_X"'] & \BO
				\end{tikzcd}
				\]
				
				We call a pair $(X,\ol{\nu}_X)$ a $\xi$-manifold if $\ol{\nu}_X$ is a $\xi$-structure on $X$ although we sometimes drop $\ol{\nu}_X$ from the notation.  A $\xi$-diffeomorphism of two $\xi$-manifolds is a diffeomorphism of the underlying manifolds that commutes with the $\xi$-structures up to homotopy over $BO$.

				\begin{definition}
					Let $\xi\colon B\to \BO$ be a map from a CW-complex $B$ of finite type. We say that $\xi$ is a \emph{universal fibration} if it is a 2-coconnected fibration. In particular, this means that $\xi_*\colon \pi_k(B)\to \pi_k(\BO)$ is injective for $k=2$ and an isomorphism for $k\geq 3$.
				\end{definition}

				\begin{definition}
					Let $\xi\colon B\to \BO$ be a universal fibration and let $\nu_X\colon X\to \BO$ be the stable normal bundle for a compact, smooth $4$-manifold $X$.  We say that $\xi$ is a \emph{normal 1-type} for $X$ if there exists a $\xi$-structure $(X,\ol{\nu}_X)$ such that $\ol{\nu}_X\colon X\to B$ is 2-connected. We call such a  choice of $\ol{\nu}_X$ a \emph{normal 1-smoothing}. 
				\end{definition}

				The following definition together with \cref{thm:1types} classifies all possible normal $1$-types of $4$-manifolds.
				\begin{definition}\cite{teichnerthesis}
					Let $\pi$ be a finitely presented group, let $w_1^{\pi}\in H^1(\pi;\Z/2)$ and $w_2^{\pi}\in H^2(\pi;\Z/2)\cup \{\infty\}$ be given. Then define $\xi(\pi,w_1^{\pi},w_2^{\pi})$ to be given by the homotopy pullback of:
					\[
					\begin{tikzcd}
						B\arrow[d,"{\xi(\pi,w_1^{\pi},w_2^{\pi})}"'] \arrow[r] \arrow[dr, phantom, "\scalebox{1.5}{\color{black}$\lrcorner$}", very near start, color=black] & B\pi \arrow[d,"w_1^{\pi}\times w_2^\pi"] \\
						\BO \arrow[r," w_1\times w_2"'] & B\Z/2\times B^2\Z/2\end{tikzcd}
					\]
					if $w_2^{\pi}\neq \infty$ and
					\[
					\begin{tikzcd}
						B\arrow[d,"{\xi(\pi,w_1^{\pi},\infty)}"'] \arrow[r] \arrow[dr, phantom, "\scalebox{1.5}{\color{black}$\lrcorner$}" , very near start, color=black] & B\pi \arrow[d,"w_1^{\pi}"] \\
						\BO \arrow[r," w_1"'] & B\Z/2 \end{tikzcd}.
					\]
					if $w_2^{\pi}=\infty$. In both cases we first realise the map $w_1^{\pi}\times w_2^{\pi}$ resp.\ $w_1^{\pi}$ as fibrations, so the maps $\xi(\pi,w_1^{\pi},w_2^{\pi})$ resp.\ $\xi(\pi,w_1^{\pi},\infty)$ are fibrations.
				\end{definition}

                In this paper we sometimes refer to $\xi(\pi,w_1^{\pi},w_2^{\pi})$ as a normal $1$-type without referencing any particular manifold. By this we mean a normal structure  $\xi(\pi,w_1^{\pi},w_2^{\pi})\colon B\to BO$ defined by elements in the cohomology of $B\pi$ as above. We do this when we seek to consider $4$-manifolds for which $\xi(\pi,w_1^{\pi},w_2^{\pi})$ is the normal $1$-type.
				
				For a given manifold we use the following theorem to decide its normal structure.
				
				\begin{theorem}[Various normal 1-types, Teichner \cite{teichnerthesis}]\label{thm:1types}
					Let $X$ be a compact, smooth $4$-manifold with $\pi:=\pi_1(X)$  and let $\widetilde{X}$ denote the universal cover of $X$.
					
					\begin{enumerate}
						\item If $w_2(\wt{X})\neq 0$, then the normal 1-type of $X$ is $\xi(\pi,w_1^{\pi},\infty)$, where $w_1^{\pi}\in H^1(\pi;\Z/2)$ is the preimage of $w_1(X)$,
						\item If $w_2(\wt{X})=0$, then the normal $1$-type of $X$ is $\xi(\pi,w_1^{\pi},w_2^{\pi})$, where $w_1^{\pi}$ and $w_2^{\pi}$ are the preimages of $w_1(X)$ and $w_2(X)$, respectively.
					\end{enumerate}
				\end{theorem}
				
				\begin{remark}\label{rmk:SpinPinxi}
					Let $\pi$ be a finitely presented group and let $w_1^{\pi}\in H^1(\pi;\Z/2)$. Then a $\xi(\pi,w_1^{\pi},0)$-structure on a manifold $X$ is a $\Pin^{+}$-structure on $X$ along with any $\pi_1$-isomorphism  $u_X^{\pi}\colon X\to B\pi$ such that $u_X^{\pi*}(w_1^{\pi})=w_1(X)$. In the case that $w_1^{\pi}=0$ this is equivalent to a $\Spin$-structure on $X$ and any $\pi_1$-isomorphism  $u_X^{\pi}\colon X\to B\pi$.
				\end{remark}
				
				\begin{definition}\label{def:xifilling}
					Let $i\colon\xi\to \xi'$ be a map of normal structures, i.e.\ a map of total spaces over $BO$, and let $(Y^n,u_Y)$ be a $\xi$-manifold.  Then we say that a $\xi'$-manifold $(X^{n+1},u_X)$ with an identification $\partial N \xrightarrow[\cong]{\varphi} M$ is a $\xi'$\emph{-filling of} $M$ if $u_N\vert_{\partial N}$ is homotopic to $i\circ u_M\circ \varphi^{-1}$ rel $BO$.
				\end{definition}
				
				We will particularly use the above definition when $\xi=\xi(\pi,w_1^{\pi},w_2^{\pi})$ and $\xi'=\xi(\pi,w_1^{\pi},\infty)$, with the map between them given by the projection $B\Z/2\times B^2\Z/2 \to B\Z/2$ and using the universal property for homotopy pullbacks.  In this case, we refer to the $\xi'$-filling $(N,u_N^{\pi})$ as a $(B\pi,w_1^{\pi})$-\emph{filling} (in the case of a closed manifold we will simply refer to it as a $(B\pi,w_1^{\pi})$-\emph{manifold}).
				
				We have the following `surgery below the middle dimension' lemma in our relative setting.
				
				\begin{lemma}\label{lem:surgeryBelowMiddle}
					Let $\xi=\xi(\pi,w_1^{\pi},w_2^{\pi})$ be a normal $1$-type, $(Y^3,u_Y)$ a $\xi$-3-manifold and $X$ its $(B\pi,w_1^{\pi})$-filling. Then there exists $(X',u_X^{\pi}\colon X'\to  B\pi)$ a $(B\pi,w_1^{\pi})$-filling of $(Y,u_Y)$ such that $x_X^{\pi}$ is a $\pi_1$-isomorphism.
				\end{lemma}
				\begin{proof}
					First make the map $u_X^{\pi}$ surjective on $\pi_1$ as follows. Let $G:=\{g_1,g_2,\dots, g_k\}$ be a generating set for $\pi$.  Form the connected-sum of 
					\[
					\#_G S^1\times S^3\xrightarrow{\sqcup u_{g_i}} B\pi
					\] with $X$, where the $u_{g_i}$ are homotopy classes of maps which send the generators of the $\pi_1(S^1\times S^3)$-factors onto the $g_i$.
					
					Now surger out the kernel of $\pi_1(X)\to \pi$ following the usual procedure, see e.g\ \cite{Wall99,Kreck99}.
				\end{proof}
				
				\subsection{Obstruction theory for \texorpdfstring{$\xi$-structures}{xi-structures}}
				
				Let $\xi=\xi(\pi,w_1^{\pi},w_2^{\pi})$ with $w_2^{\pi}\neq \infty$ be a normal $1$-type. We wish to define an obstruction for $(B\pi,w_1^{\pi})$-manifolds to admit a $\xi$-structure extending the one given on a boundary.  
				
				\begin{proposition}\label{prop:relative_obstruction_class}
					Let $(Y,u_Y)$ be a, possibly empty,  $\xi$-3-manifold and let $(X,u_X^{\pi}\colon X\to B\pi)$ be its $(B\pi,w_1^{\pi})$-filling. Then there exists a unique obstruction $\smallO(X,Y,u_y,u_X^{\pi})\in H^2(X,Y;\Z/2)$ that vanishes if and only if $u_Y$ extends to a $\xi$-structure on $X$. 
					
					Furthermore, if $Y=\emptyset$ then $\smallO(X,\emptyset,\emptyset,u_X^{\pi})$ splits as a sum $w_2(M)+ u_M^{\pi*}(w_2^{\pi})$.
				\end{proposition}
				
				\begin{proof}
					
					Consider the following diagram.
					\[
					\begin{tikzcd}
						Y \arrow[d] \arrow[rr,"u_Y"] & & B  \arrow[r] \arrow[dr, phantom, "\scalebox{1.5}{\color{black}$\lrcorner$}" , very near start, color=black] & B\pi \arrow[d,"w_1^{\pi}\times w_2^\pi"] &\\
						X \arrow[rrru,"u_X^{\pi}" near start, end anchor= south west]\arrow[rr,"\nu_X"']&&\BO \arrow[from=u, crossing over, "\xi"' near start] \arrow[r," w_1\times w_2"'] & B\Z/2\times B^2\Z/2\ar[r]&B^2\Z/2.
					\end{tikzcd}
					\]
					The outermost square already commutes up to homotopy since $u_Y$ is a $\xi$-structure.  The triangle starting at $X$ does not necessarily commute, but it does after postcomposition with the projection $B\Z/2\times B^2\Z/2 \to B\Z/2$ since $X$ is a $(B\pi,w_1^{\pi})$-filling.  Hence we only have to consider the diagram after projecting to the $B^2\Z/2$-factor.
					
					The middle square comes with a homotopy $B\times I\to B^2\Z/2$ making the square commute. Define the homotopy $\mathscr{H}\colon Y\times I \to B^2\Z/2$ using precomposition. Using the $H$-space structure on $B^2\Z/2$, the difference of the maps $w_2(X)\colon X\to B^2\Z/2$ and $(u_X^{\pi})^*(w_2^{\pi})\colon X\to B^2\Z/2$ defines a map
					\[ 
					w_2(X) - (u_X^{\pi})^*(w_2^{\pi})\colon X\to B^2\Z/2.
					\]
					The homotopy $\mathscr{H}$ then gives a trivialisation of $w_2(X) - (u_X^{\pi})^*(w_2^{\pi})$ precomposed with the inclusion $Y\to X$, and hence defines a relative cohomology class $\smallO(X,Y,u_y,u_X^{\pi})\in H_2(X,Y;\Z/2)$.  It is clear from the construction that this cohomology class has the desired properties.
				\end{proof}
				
				\begin{remark}
					In particular, for $w_1^{\pi}=0=w_2^{\pi}$, a $u_Y$-structure on $Y$ includes a $\Spin$ structure on $Y$ and the obstruction $\smallO(X,Y,u_y,u_X^{\pi})$ coincides with the relative Stiefel-Whitney class $w_2^{R}$ defined by Kervaire \cite{Kervaire57}.
				\end{remark}
				
				\begin{definition}
					If $(X,u_X^{\pi})$ is a closed $(B\pi,w_1^{\pi})$-manifold then we set $\smallO(X,u_X^{\pi}):=\smallO(X,\emptyset,\emptyset,u_X^{\pi})$. If $(X_1,u_{X_1}^{\pi}), (X_2,u_{X_2}^{\pi})$ are two $(B\pi,w_1^{\pi})$-fillings of $Y$ when we use the relative Mayer-Vietoris sequence to show that the map 
					\[
					H^2(X_1\cup X_2,Y;\Z/2)\xrightarrow{\cong}H^2(X_1,Y;\Z/2)\oplus H^2(X_2,Y;\Z/2)
					\]
					is an isomorphism. We define $\smallO(X_1\cup X_2,Y,u_Y,u_{X_1\cup X_2}^{\pi})$ to be the preimage of the sum $\smallO(X_1,Y,u_Y,u_{X_1}^{\pi})+\smallO(X_2,Y,u_Y,u_{X_2}^{\pi})$.
				\end{definition}
				
				In \Cref{lem:doubleObstructionClass} we identify the obstruction $\smallO$ of manifolds obtained from gluing two manifolds along a common boundary.
				
				\begin{lemma}\label{lem:doubleObstructionClass}
					Let $\xi=\xi(\pi,w_1^{\pi},w_2^{\pi})$, $w_2^{\pi}\neq\infty$ be a normal $1$-type. Let $(Y,u_Y)$ be a $\xi$-3- manifold. Let $X_1,X_2$ be two $(B\pi,w_1^{\pi})$-fillings of $Y$. Then the inclusion maps $X_1,X_2\hookrightarrow X_1\cup \ol{X_2}$ induce an equality
					\[
					PD(\smallO(X_1,Y,u_Y,u_{X_1}^{\pi}))+PD(\smallO(\ol{X_2},Y,u_Y,u_{X_2}^{\pi}))=PD(\smallO(X_1\cup \ol{X_2},u_{X_1\cup X_2}^{\pi})).
					\]
				\end{lemma}

				\begin{proof}
					Write $M:=X_1\cup\ol{X_2}$. The following diagram, where the horizontal maps come from the relevant Mayer-Vietoris sequences and the unlabelled vertical map comes from the long exact sequence of the pair, commutes, (use naturality of cap product).
					\[
					\begin{tikzcd}
						0\ar[r]&H^2(M,Y;\Z/2)\ar[r,"\cong"]\ar[d]& H^2(X_1,Y;\Z/2)\oplus H^2(X_2,Y;\Z/2)\ar[dd,"PD"]\ar[r]&0\\
						&H^2(M;\Z/2)\ar[d,"PD"]&\\
						&H_2(M;\Z/2))&H_2(X_1;\Z/2)\oplus H_2(X_2;\Z/2)\ar[l,"+"]
					\end{tikzcd}
					\]
					
					Now the class $(\smallO(X_1,Y,u_Y,u_{X_1}^{\pi}),\smallO(X_2,Y,u_Y,u_{X_2}^{\pi}))$ maps, under the top horizontal isomorphism to the obstruction class $\smallO(M,Y,u_Y,u_M^{\pi})$ for $M$ to have a $\xi$-structure, which restricts to $u_Y$ on $Y$. This class in turn maps to $\smallO(M,u_M^{\pi})$, the obstruction to $M$ having a $\xi$-structure itself. Now the commutativity of the diagram concludes the proof since the bottom map, is just the summing map, since our coefficients are $\Z/2$. 
				\end{proof}
				
				\begin{definition}\label{def:sr_characteristic}
					A class $c\in\pi_2(X)$ is  {\it s-characteristic (spherical characteristic)}, if for the $\pmod{2}$ augmentation map $\varepsilon_2\colon \Lambda\to \Z\to \Z/2$ we have
					\[
					\varepsilon_2\lambda(c,\alpha)=\varepsilon_2\lambda(\alpha,\alpha),
					\]
					for all classes $\alpha\in\pi_2(X)$.  Similarly, $w\in \pi_2(X)$ is called \emph{r-characteristic} ($\RP^2$\emph{-characteristic}) if $Q(h_2(w),\alpha)=Q(\alpha,\alpha)$ for all classes $\alpha\in H_2(X;\Z/2)$ representable by maps $\alpha\colon \RP^2\to X$, where $Q$ denotes the $\Z/2$-valued intersection form on $X$.
				\end{definition}
				Note that for two spherical classes $\alpha,\beta\in\pi_2(X)$ we have $Q(h_2(\alpha),h_2(\beta))=\varepsilon_2(\lambda(\alpha,\beta))$, showing the analogy between the definitions of $s-$ and $r-$ characteristic elements above. 
				
				\begin{lemma}\label{lem:sphericalCharacteristic}
					Let $\xi=\xi(\pi,w_1^{\pi},w_2^{\pi})$ be a normal $1$-type. Let $(Y,u_Y)$ be a $\xi$-3-manifold and let $(X,u_X^{\pi})$ be a $(B\pi,w_1^{\pi})$-filling with $u_X^{\pi}$ a $\pi_1$-isomorphism (see \Cref{lem:surgeryBelowMiddle}). Let $w\in H_2(X;\Z/2)$ be the Poincar\'{e} dual of the class $\smallO(X,Y,u_Y,u_X^{\pi})$ defined above. Assume that $w$ is a spherical class. Then any lift $c\in\pi_2(X)$ of $w$ is {\it s-characteristic}.
				\end{lemma}
				
				\begin{proof}
					We will first show the lemma in the case that $Y=\emptyset$. Let $\alpha\in\pi_2(X)$. Let $Q$ be the $\Z/2$-intersection form on $H_2(M;\Z/2)$. Then $\varepsilon_2(\lambda(c,\alpha))=Q(w,h_2(\alpha))=\langle \smallO(X,u_X^{\pi}),h_2(\alpha)\rangle$. By \cref{prop:relative_obstruction_class},  $\smallO(X,u_X^{\pi})$ splits as a sum $(u_X^{\pi})^{*}w_2^{\pi}+w_2(X)$. The  term $\left< u_X^{\pi*}(w_2^{\pi }),h_2(\alpha)\right>=0$ vanishes  because $u_{M*}^{\pi}h_2=0$. Now
                    \[
					\left<u_X^{\pi*}(w_2^{\pi })+w_2(X),h_2(\alpha)\right>=\left<w_2(X),h_2(\alpha)\right>=Q(h_2(\alpha),h_2(\alpha))=\varepsilon_2\lambda(\alpha,\alpha)
					\]
					since $w_2(X)$ is characteristic, which proves the claim.
					
					For the general case $X$, with boundary $Y$, let $\alpha\in\pi_2(X)$ be any element and  consider the double $DX$ along with the map $u^{\pi}_{DX}:=u^{\pi}_X\cup u^{\pi}_X\colon DX\to B\pi$. Fix one of the inclusions $X\to DX$ and denote the switch map $DX\to DX $ by $D$. The dual of the obstruction class $\smallO(DX,u_{DX}^{\pi})$ is the sum $w_D=w+Dw$ according to \Cref{lem:doubleObstructionClass}. Since $w$ lifts to $c$, we can realise a lift of $w_D$ as the ambient connected sum $c\# Dc$.  We get 
					
					\begin{align*}
						&\varepsilon_2(\lambda^{X}(\alpha,c))=\varepsilon_2(\lambda^{DX}(\alpha,c))=\varepsilon_2(\lambda^{DX}(\alpha,c+Dc))=\varepsilon_2(\lambda^{DX}(\alpha,\alpha))=\varepsilon_2(\lambda^{X}(\alpha,\alpha)).
					\end{align*}
					Where the first and last equality follows by naturality of the map $X\to DX$, second equality holds since $\alpha$ and $Dc$ are disjoint and third by the closed case.
				\end{proof}

				\section{Intersections and self intersections}\label{sec:intersections}

				Throughout this section let $X$ be a compact 4-manifold, let $\pi:=\pi_1(X)$, let $\Lambda:= \Z[\pi]$, and let $w_1=w_1(X)\colon \pi\to \Z/2$ be the first Stiefel-Whitney class. 
				
				\subsection{Preliminaries}
				Recall there is a map  \[\mu\colon \pi_2(X)\to \Z[\pi]/(g- w_1(g)g^{-1},1)\] given by choosing an immersed representative for $\alpha$ and then calculating the self-intersections (counting group elements), ignoring those with the trivial group element (see \cite[Theorem 5.2]{Wall99}) to make this not depend on the choice of representative\footnote{If one is given an immersed representative $\wt{\alpha}$ for $\alpha\in\pi_2(X)$, we will write $\mu(\wt{\alpha})\in \Z[\pi]/(g-w_1(g)g^{-1})$, where we do not have to mod out by this indeterminacy.}. The indeterminacy $g\sim w_1(g)g^{-1}$ comes from there not being a preferred ordering of the sheets at a self-intersection point.  We define an altered version of this invariant.
				
				\begin{definition}
					Let $\wt{\alpha}\colon S^2\to X$ be an immersed sphere with $n$ self intersections. Let $\mathfrak{o}(\wt{\alpha})$ be a set of $n$  ordered pairs of points $a_i,b_i\in S^2$ for $i=1,\dots,n$ such that $\wt{\alpha}(a_i)=\wt{\alpha}(b_i)$. Assume $\mathfrak{o}$ spans all the double points. Then we call a pair $(\wt{\alpha},\mathfrak{o}(\wt{\alpha}))$ an \emph{immersed sphere with ordered double points}.
				\end{definition}
				
				Now for any immersed sphere with ordered double points $(\wt{\alpha},\mathfrak{o}(\wt{\alpha}))$ the self intersection $\mu(\wt{\alpha},\mathfrak{o}(\wt{\alpha}))$ is well-defined in $\Z[\pi]$.
				
				\begin{definition}\label{def:muI}  Let  \[
					\mudash\colon\pi_2(X)\to I/(g-1 - w_1(g)(g^{-1}-1)).
					\]
					be defined as follows. Let $\alpha\in\pi_2(X)$. Represent $\alpha$ as an immersed sphere with ordered double points $(\wt{\alpha},\mathfrak{o}(\wt{\alpha}))$.
					Now define $\mudash(\alpha) = [\mu(\wt{\alpha},\mathfrak{o}(\wt{\alpha})) -\varepsilon(\mu(\wt{\alpha},\mathfrak{o}(\wt{\alpha})))]$, where the brackets represent taking the quotient.
				\end{definition}
				We note that in the case that $X$ is orientable the above definition is equivalent to $\mudash(\alpha)=\mudash(\wt{\alpha})-\varepsilon(\mudash(\wt{\alpha}))$ since $\varepsilon$ factors through $\Z[\pi]/(g- w_1(g)g^{-1})$ if and only if $X$ is $w_1=0$.
				
				\begin{proposition}
					The function of \cref{def:muI}
					\[
					\mudash\colon\pi_2(X)\to I/(g-1 - w_1(g)(g^{-1}-1)).
					\]
					is well-defined.
				\end{proposition}
				\begin{proof}
					First we show that it is independent on the choice of $\mathfrak{o}(\wt{\alpha})$. Let  $\mathfrak{o}'(\wt{\alpha})$ be $\mathfrak{o}(\wt{\alpha})$ with one ordered pair $(a_i,b_i)$ swapped and assume that the group element at $(a_i,b_i)$ dictated by $\mathfrak{o}(\wt{\alpha})$ is $g\in\pi$. then the element at $(a_i,b_i)$ is dictated by $w_1(g)g^{-1}$ and we get 
					\begin{align*}
						&\mu(\wt{\alpha},\mathfrak{o}(\wt{\alpha})) -\varepsilon(\mu(\wt{\alpha},\mathfrak{o}(\wt{\alpha})))-\mu(\wt{\alpha},\mathfrak{o}'(\wt{\alpha})) +\varepsilon(\mu(\wt{\alpha},\mathfrak{o}'(\wt{\alpha})))=\\&=g-w_1(g)g^{-1}-\varepsilon(g)+\varepsilon(w_1(g)g^{-1})=\\
						&=g-w_1(g)g^{-1}-1+w_1(g)=g-1-w_1(g)(g^{-1}-1),
					\end{align*}
					which vanishes in $I/(g-1 - w_1(g)(g^{-1}-1))$.
					The proof that $\mudash$ does not depend on the choice of immersion is the same as in \cite{Wall99}.
				\end{proof}

				Analogously to $\mudash$, we consider the following $\Z$-intersection form 
				\[
				\lambdadash:\pi_2(X)\times \pi_2(X)\to I
				\]
				defined as $\lambdadash(\alpha,\beta)=\lambda(\alpha,\beta)-\varepsilon(\lambda(\alpha,\beta))$ where $\lambda:\pi_2(X)\times\pi_2(X)\to \Lambda$ is the usual equivariant intersection form. This form was considered in \cite{kosanovic2023newapproachlightbulb}.
				
				The main downside of this form is that we lose the $\Lambda$-hermitian structure, in particular the $\Lambda$-bilinearity of the usual form $\lambda$. 
				
				\begin{lemma}\label{lem:mudash_lambdadash_sum_formula}
					The form $\mudash$ is $\Z$-quadratic and the form $\lambdadash$ is the $\Z$-bilinear form associated to the (underlying) $\Z$-quadratic form $\mudash$, i.e.\ for all $\alpha,\beta\in \pi_2(X)$
					\[\mudash(\alpha+\beta)=\mudash(\alpha)+\mudash(\beta)+[\lambdadash(\alpha,\beta)]\in I/({g-1 -w_1(g)(g^{-1}-1))}.\]
				\end{lemma}
				
				\begin{proof}
					We leave to the reader to show that that $\lambdadash$ is a $\Z$-bilinear form and that $\mudash(n\alpha )=n^2\mudash(\alpha)$ for $n\in\Z$. Let $\alpha,\beta\in \pi_2(X)$, let $(\wt\alpha, \mathfrak{o}(\wt{\alpha}))$ and $(\wt{\beta},\mathfrak{o}(\wt{\beta}))$ be immersed sphere representatives of $\alpha$ and $\beta$ with ordered double points. We can naturally form an immersed sphere with ordered double points $(\wt{\alpha}+\wt{\beta},\mathfrak{o}(\wt{\alpha})\cup\mathfrak{o}(\wt{\beta})\cup\mathfrak{o}_{\wt{\alpha}+\wt{\beta}})$, where $\mathfrak{o}_{\wt{\alpha}+\wt{\beta}}$ is the ordering of the new intersection points between $\wt\alpha$ and $\wt\beta$ given by ordering the pair $(\wt\alpha,\wt\beta)$. We calculate
					\begin{align*}			
						\mudash(\alpha+\beta)&=[\mu(\wt{\alpha}+\wt{\beta},\mathfrak{o}(\wt{\alpha})\cup \mathfrak{o}(\wt{\beta})\cup\mathfrak{o}_{\wt{\alpha}+\wt{\beta}})-\varepsilon(\mu(\wt{\alpha}+\wt{\beta},\mathfrak{o}(\wt{\alpha})\cup\mathfrak{o}(\wt{\beta})\cup\mathfrak{o}_{\wt{\alpha}+\wt{\beta}}))]\\
						&=[\mu(\wt\alpha, \mathfrak{o}(\wt{\alpha}))+\mu(\wt{\beta},\mathfrak{o}(\wt{\beta}))+\lambda(\alpha,\beta)-\varepsilon(\mu(\wt\alpha, \mathfrak{o}(\wt{\alpha})))-\varepsilon(\mu(\wt{\beta},\mathfrak{o}(\wt{\beta})))-\varepsilon(\lambda(\alpha,\beta))]\\
						&=\mudash(\alpha)+\mudash(\beta)+[\lambda(\alpha,\beta)-\varepsilon(\lambda(\alpha,\beta))]\\
						&=\mudash(\alpha)+\mudash(\beta)+\lambdadash(\alpha,\beta),
					\end{align*}
					where the square brackets represent the corresponding quotient.
				\end{proof}
				
				We will use the involution on the group ring $\ol{\cdot}\colon\Z[\pi]\to \Z[\pi]$ defined by linearly extending $\ol{g}=w_1(g)g^{-1}$ for $g\in \pi$ to the whole group ring.  Recall that, given $\alpha\in\pi_2(X)$, the sum $\mu(\alpha)+\ol{\mu(\alpha)}$ lies in $\Z[\pi]/(1)$, since the indeterminacy due to the sheet choices cancels out.  As an easy consequence, the sum $\mudash(\alpha)+\ol{\mudash(\alpha)}$ lies in $I$.  Using this fact, the forms $\mudash$ and $\lambdadash$ satisfy the following relation.
				
				\begin{lemma}\label{lem:PropertyOfLambdadash}
					For any $\alpha\in \pi_2(X)$ we have \[\lambdadash(\alpha,\alpha)=\mudash(\alpha)+\ol{\mudash(\alpha)}\in I.\]
				\end{lemma}
				
				\begin{proof}
					Pick an immersed representative $\wt{\alpha}$ for $\alpha$.  By Wall \cite[Theorem 5.2]{Wall99}, we have the corresponding version for the original $\mu$ and $\lambda$, i.e.\
					\begin{equation}\label{eq:wall_mu_lambda}
						\lambda(\alpha,\alpha)=\mu(\wt{\alpha}) +\ol{\mu(\wt{\alpha})} +e(\wt{\alpha}),
					\end{equation}
					where $e(\wt{\alpha})$ is the Euler number of the immersion $\wt\alpha$.  This gives
					\begin{align*}
						\lambdadash(\alpha,\alpha) &= \lambda(\alpha,\alpha) - \varepsilon(\lambda(\alpha,\alpha))\\
						&= \mu(\wt{\alpha}) + \ol{\mu(\wt{\alpha})} +e(\wt{\alpha}) - \varepsilon (\lambda(\alpha,\alpha))\\
						&= \mudash(\alpha) + \varepsilon(\mu(\wt{\alpha})) + \ol{\mudash(\alpha)} +\varepsilon(\ol{\mu(\wt{\alpha})}) +e(\wt{\alpha}) -\varepsilon(\lambda(\alpha,\alpha)) \\
						&=\mudash(\alpha) +\ol{\mudash(\alpha)},
					\end{align*}
					where the first equality is by the definition of $\lambdadash$, the second by \Cref{eq:wall_mu_lambda}, the third by the definition of $\mudash$, and the fourth by applying $\varepsilon$ to \Cref{eq:wall_mu_lambda}.
				\end{proof}

				\subsection{The form \texorpdfstring{$\mudashj$}{mudashj}}
				
				We are interested in making the self-intersection form $\mudash$ valued in $H_1(\pi;\Z/2)$.  To do this, we need a map $I\to H_1(\pi;\Z/2)$ which vanishes on the ideal $(g-1 - w_1(g)(g^{-1}-1))\subset I$.  In that direction, we prove the following.

				\begin{lemma}\label{lem:tensorIsos}
					Let $\pi$ be a group, $I$ resp.\ $J$ be the $\Z$ resp.\ $\Z/2$ augmentation ideals. Then we have the following: 
					
					\begin{enumerate}
						\item\label{it:H1group1} $IJ=JI$.
						\item\label{it:H1group2} For any right $\Lambda$-module $M$ there are isomorphisms
						\begin{itemize}
							\item $M\otimes_{\Lambda}\Z/2\cong M/MJ$ given by $m\otimes 1\mapsto m$,
							\item and $M\otimes_{\Lambda}\Z\cong M/MI$ given by $m\otimes 1\mapsto m$.
						\end{itemize}
						\item\label{it:H1group2b} For all $i\in I$ we have $i^2, 2i\in IJ$.
						\item\label{it:H1group3} There is an isomorphism of $\Z/2$-vector spaces $H_1(\pi;\Z/2)\cong I/IJ$ induced by the map sending $g\in \pi$ to $g-1\in I$. The inverse of this map is given as follows. The ideal $I$, as an abelian group is generated by elements $\{g-1\mid g\in \pi\}$. We define the map by sending $\sum_i(g_i-1)+\sum_j(1-g_j)\in I$ to $\prod_ig_i\prod_jg_j^{-1}$.
					\end{enumerate}
					
				\end{lemma}
				\begin{proof}
					
					For (\ref{it:H1group1}), let $i\in I$ and $j\in J$.  Then 
					\[ij=i(j-\varepsilon(j))+i\varepsilon(j)=i(j-\varepsilon(j))+\varepsilon(j)i\]
					
					Since $j-\varepsilon(j)\in I\subset J$ we find that $ij\in JI$. We apply an analogous calculation to show that $ji\in IJ$.
					
					For the first isomorphism of (\ref{it:H1group2}), the $\Lambda$-bilinear map $M\times \Z/2\to M/MJ$ given by $(m,k)\mapsto mk$ is well-defined since $2\in J$ and so $m2\in MJ$. This map descends to a map $M\otimes_{\Lambda}\Z/2\to M/MJ$. It is clearly surjective. If $m\otimes 1$ is mapped to zero then $m\in MJ$ and hence $m=m'j$. Now $m'j\otimes 1=m'\otimes \varepsilon_2(j)=0$. 
					
					The proof of the second isomorphism is analogous.
					
					For (\ref{it:H1group2b}) we have $I\subset J$ and so $i^2\in IJ$. Additionally we have $2\in J$ which proves the claim.
					
					For (\ref{it:H1group3}) consider the short exact sequence:
					\[0\to I\to \Lambda\to \Z. \]
					
					Tensoring the sequence with $\Z/2$ gives the long exact sequence
					\[\Tor_1^{\Lambda}(\Lambda,\Z/2)\to \Tor_1^{\Lambda}(\Z,\Z/2)\to I\otimes_\Lambda \Z/2\to \Lambda\otimes_{\Lambda}\Z/2\to \Z\otimes_\Lambda \Z/2\to 0\]
					
					where the second term is isomorphic to $H_1(\pi;\Z/2)$.  On the right side of the sequence we have the isomorphism $\Z/2\cong \Lambda\otimes_{\Lambda}\Z/2\to \Z\otimes \Z/2\cong \Z/2$, so we obtain an exact sequence
					\[\Tor_1^{\Lambda}(\Lambda,\Z/2) \to H_1(\pi;\Z/2)\to I\otimes_{\Lambda}\Z/2\to 0.\]
					
					Since the term $\Tor_1^{\Lambda}(\Lambda,\Z/2)$ vanishes, the map $H_1(\pi;\Z/2)\to I\otimes_\Lambda \Z/2$ is an isomorphism.  Postcomposing with the isomorphism from (\ref{it:H1group2}), we obtain an isomorphism $\psi\colon H_1(\pi;\Z/2)\to I/IJ$, which we now describe. Let $B$ be a CW-complex with a single $0$ cell and, and one $1$-cell for each generator in some fixed generating set $G$ of $\pi$. Let $E\pi$ be the universal cover, which is a contractive space and we can take its chain complex $C_3(E\pi)\to C_2(E\pi)\to C_1(E\pi)\to C_0(E\pi)$ as the free resolution of $\Z$. Note that $C_0(E\pi)\cong \Lambda$ and $C_1(E\pi)\cong \prod_{g\in G}\Lambda$.
					
					We have the following commutative diagram
					\[
					\begin{tikzcd}
						C_2(E\pi)\ar[d,equals]\ar[r]&C_1(E\pi)\ar[d,"\cong"]\ar[r]&C_0(E\pi)\ar[d,"\cong"]\ar[r]&\Z\ar[d,equals]\\
						C_2(E\pi)\ar[d]\ar[r]&\prod_{g\in G}\Lambda\ar[d,"e_1^g\to g-1" ]\ar[r]&\Lambda\ar[r]\ar[d,equals]&\Z\ar[d,equals]\\
						0\ar[r]&I\ar[r]&\Lambda\ar[r]&\Z.
					\end{tikzcd}
					\]
					Upon tensoring the rows with $\Z/2$ and taking homology we find that the map
					\[
					H_1(\pi;\Z/2)\cong\Tor^{\Lambda}_1(\Z,\Z/2)\to \ker(I\otimes_{\Lambda}\Z/2\to \Lambda\otimes_{\Lambda} \Z/2)\cong I\otimes_{\Lambda}\Z/2 
					\]
					is given by $g\to (g-1)\otimes 1$, so $\psi\colon H_1(\pi;\Z/2)\to I/IJ$ is given by $g\mapsto g-1$. The map $\psi$ is a map of $\Z/2$-vector spaces because 
					\begin{align*}
						g_1+g_2=g_1g_2\mapsto g_1g_2-1&=(g_1-1)(g_2+1)+g_1-g_2\\&\equiv g_1+g_2
						\equiv  (g_1-1)+(g_2-1)\pmod{IJ}.\qedhere
					\end{align*}
				\end{proof}
				
				This produces a map $I\to H_1(\pi;\Z/2)$ given by quotienting by $IJ$.  We prove the following fact about the kernel of this quotient map, which will be useful in \Cref{sec:vanishing}.  The reader can skip this for now and move onto the definition of $\mudashj$.
				
				\begin{lemma}\label{lem:KernelOfTheMapItoHomology}
					Let $\pi$ be a finitely presented group. Then the kernel of the map $$I/({g-1 - w_1(g)(g^{-1}-1)})\to H_1(\pi;\Z/2)$$ defined in \cref{lem:tensorIsos}(\ref{it:H1group3}) is generated as a $\Lambda$-module by \[\{g_ig_j-g_i-g_j+1,2g_i-2\mid g_i,g_j\in \pi\},\] where the module structure is given by conjugation. 
				\end{lemma}
				\begin{proof}
					From \Cref{lem:tensorIsos}(\ref{it:H1group3}) we know that the kernel of the map $I\to H_1(\pi;\Z/2)$ is $IJ$. Now we have a factorisation
					\[
					I\xtwoheadrightarrow{q} I/({g-1 - w_1(g)(g^{-1}-1)})\xrightarrow{f} H_1(\pi;\Z/2)
					\]
					Using the surjectivity of $q$ we show that $q(IJ)=\ker(f)$. The module $I$ is generated (as an abelian group) by the set $\{g-1\mid g\in\pi\}$ and the module $J$ is generated (as an abelian group) by $\{g-1\mid g\in\pi\}\cup\{2\}$. It follows that the modules $IJ$ and $p(IJ)$ are generated by the product of pairs of the above generators $\{g_ig_j-g_i-g_j+2\mid g_i,g_j\in\pi\}\cup\{2g-2\mid g\in\pi\}$.
				\end{proof}
				
				We now define $\mudashj$ using \Cref{lem:tensorIsos}.  Note that since $g-g^{-1}=(g+1)(g^{-1}-1)\in IJ$ and $g+g^{-1}-2=(g+1)(1-g^{-1})\in I^2\subset IJ$ we have that $(g-1-w_1(g)(g^{-1}-1))\subset IJ$.  Hence we can define $\mudashj\colon \pi_2(X)\to I/IJ$, simply by postcomposing $\mudash$ by the quotient map.
				
				\begin{definition}\label{def:muJlambdaJ}
					Let \[
					\mudashj\colon \pi_2(X)\to I/IJ
					\]
					be defined as $\mudash$ postcomposed with the quotient map $I/((g-1-w_1(g)(g^{-1}-1))\to I/IJ$.  Similarly define
					\[
					\lambdadash_J: \pi_2(X)\times \pi_2(X) \to I/IJ
					\]
					as $\lambdadash$ postcomposed with the quotient map $I\to I/IJ$.
				\end{definition}

				We have the following refinement of these forms.
				
				\begin{proposition}\label{prop:mudashj_lamdadashj_well-defined}
					The forms $\mudashj$ and $\lambdadash_J$ of \Cref{def:muJlambdaJ} induce a well-defined $\Z/2$-quadratic form and (an underlying) $\Z/2$-bilinear form on $\pi_2(X)\otimes_{\Z}\Z/2$ and $\left(\pi_2(X)\otimes_{\Z}\Z/2\right) \times \left(\pi_2(X)\otimes_{\Z}\Z/2\right)$, respectively.
				\end{proposition}
				
				\begin{proof}
					Let $\alpha\in\pi_2(X)\otimes_{\Z}\Z/2$. Using that $\pi_2(X)\otimes_{\Z}\Z\cong \pi_2(X)$ and the long exact sequence obtained by tensoring the sequence $\Z\xrightarrow{2}\Z\rightarrow \Z/2$ using $\pi_2(X)\otimes_{\Z}-$ shows that different lifts $\ol{\alpha}\in\pi_2(X)$ of $\alpha$ differ by $2\pi_2(X)$. Let $2\beta\in 2\pi_2(X)$. By the following calculation
					\[
					[\mudash(\ol{\alpha}+2\beta)]=[\mudash(\ol{\alpha})]+[\mudash(2\beta)]+[\lambdadash(\ol{\alpha},2\beta)]=[\mudash(\ol{\alpha})]
					\]
					
					We see that $\mudashj$ does not depend on the choice of lift.  Here $[-]$ denotes the projection $I/g\sim g^{-1}\to I/IJ$, and the last equality holds since $I/IJ$ is a $\Z/2$-vector space by \Cref{lem:tensorIsos}(\ref{it:H1group3}).
					
					Similarly for $\alpha, \alpha'\in \pi_2(X)\otimes_{\Z}\Z/2$, their lifts $\ol{\alpha},\ol{\alpha'}\in \pi_2(X)$ their lifts, pick $2\beta,2\beta'\in 2\pi_2(X)$ and so we have
					\[
					[\lambdadash(\ol{\alpha}+2\beta,\ol{\alpha'}+2\beta')]=\lambdadash(\ol{\alpha},\ol{\alpha'})
					\]
					again using $I/IJ$ is a $\Z/2$- vector space.  The rest of the claim follows from \Cref{lem:mudash_lambdadash_sum_formula}.
				\end{proof}
				
				The induced forms from \Cref{prop:mudashj_lamdadashj_well-defined} will still be denoted by $\mudashj$ and $\lambdadash_J$.  We have the following algebraic fact which will be useful when establishing the properties of $\mudashj$ and $\lambdadash_J$.
				
				\begin{lemma}
					Let $u\in J^2$ and $u\in I$. Then $u\in IJ$.
				\end{lemma}
				
				\begin{proof}
					We can write $u=j_1j_2$ for $j_1,j_2\in J$. Since $u\in I$ we have that  $\varepsilon(j_1j_2)=0$ and so either $\varepsilon(j_1)=0$ or $\varepsilon(j_2)=0$ and so either $j_1\in I$ or $j_2\in I$. Now use that $IJ=JI$ (\Cref{lem:tensorIsos}(\ref{it:H1group1})).
				\end{proof}
				
				We now state and prove the various properties that $\mudashj$ and $\lambdadash_J$ satisfy.
				
				\begin{lemma}\label{lem:PropertiesOfMuLambdaIJ}
					Let $g\in \pi$, $a,b \in \pi_2(X)\otimes_{\Z}\Z/2$, $j\in J$ and $i\in I$.  Then the following equations hold in $I/IJ$.
					\begin{enumerate}
						\item\label{it:propertiesMuLambda1} $\mudashj(ga)=\mudashj(a)$.
						\item\label{it:propertiesMuLambda2} $\lambdadash_J(a,a)=0$.
						\item\label{it:propertiesMuLambda3} $\lambdadash_J(ga,a)=(g-1)\varepsilon_2(\lambda(a,a))$.
						\item\label{it:3AxiomLambda} $\lambdadash_J(ja,b)=(j-\varepsilon(j))\varepsilon_2(\lambda(a,b))$.
						\item\label{it:3AxiomMu} $\mudashj(ja+b)=(j-\varepsilon(j))\varepsilon_2(\lambda(a,a)+\lambda(a,b))+\mudashj(b)$.
					\end{enumerate}
				\end{lemma}

				\begin{proof}
					For (\ref{it:propertiesMuLambda1}), note that $g^{-1}-g=(g^{-1}-1)(g+1)\in IJ$  and so $g\equiv g^{-1}\pmod{IJ}$.  We have that
					\[
					\mudashj(ga)=g\mudashj(a)g^{-1}=(g-1)\mudashj(a)g+\mudashj(a)g \equiv \mudashj(a)(g-1)+\mudashj(a)\equiv \mudashj(a)\pmod{IJ},
					\]
					by using \cref{prop:mudashj_lamdadashj_well-defined} and that $(g-1)\in J, \mudashj(a)\in I$.  This completes the proof of (\ref{it:propertiesMuLambda1}).
					
					For (\ref{it:propertiesMuLambda2}) we calculate $\mudashj(2a)=2\mudashj(a)+\lambdadash_J(a,a)$. Since $I/IJ$ is $2$-torsion immediately obtain that~$\lambdadash_J(a,a)=0$.  This completes the proof of (\ref{it:propertiesMuLambda2}).
					
					For (\ref{it:propertiesMuLambda3}) we calculate that
					\begin{align*}
						\lambdadash(ga,a)&=g\lambda(a,a)-\varepsilon g\lambda(a,a)\\
						&=g(\lambdadash(a,a)+\varepsilon(\lambda(a,a)))-\varepsilon\lambda(a,a)\\
						&=g\lambdadash(a,a)+(g-1)(\varepsilon(\lambda(a,a))-\varepsilon_2(\lambda(a,a)))+(g-1)\varepsilon_2(\lambda(a,a))\\
						&=(g-1)\varepsilon_2(\lambda(a,a)),
					\end{align*}
					where the first and second equalities follow from the definition of $\lambdadash$, the third by using $(g-1)(\varepsilon(\lambda(a,a))-\varepsilon_2(\lambda(a,a)))\in IJ$, and the last by (\ref{it:propertiesMuLambda2}).  This completes the proof of (\ref{it:propertiesMuLambda3}).
					
					For (\ref{it:3AxiomLambda}) we  let $j\in J$. We can express $j$ as $i-y$ for $y=\varepsilon(j)$ an even integer and $i\in I$. Then $\lambdadash(ja,b)=\lambdadash(ia,b)+2\lambdadash(a,b)=\lambdadash(ia,b)$ in $IJ$, so it is sufficient to assume that $j=i\in I$.  Write $i$ as a sum of $n$ expressions of the form $(g-1)$ and $(1-g)$ for $g\in \pi$. We will proceed by the induction on $n$. If $n=0$ then $i=0$ and there is nothing to prove.
					
					If $n=1$ then $i=1-g$ or $i=g-1$. The proof will be the same in all cases so we assume the former. We  obtain
					\begin{align*}
						\lambdadash((g-1)a,b)&=\lambdadash(ga,b)-\lambdadash(a,b)\\
						&=\lambda(ga,b)-\varepsilon\lambda(ga,b)-(\lambda(a,b)-\varepsilon(\lambda(a,b))\\
						&=g\lambda(a,b)-\varepsilon\lambda(a,b)-\lambda(a,b)+\varepsilon\lambda(a,b)\\
						&=(g-1)\lambda(a,b)\\
						&=(g-1)(\lambda(a,b)-\varepsilon_2(\lambda(a,b)))+(g-1)\varepsilon_2\lambda(a,b)\\
						&\equiv (g-1)\varepsilon_2\lambda(a,b) \pmod{IJ},
					\end{align*}
					where in the sixth line we used that $\lambda(a,b)-\varepsilon(\lambda(a,b))\in J$.
					
					For the induction step we take $i=(g-1)+i'$ or $i=(1-g)+i''$ where the number of summands of $i',i''\in I$ is smaller than the number of summands in $i$. In the former case the induction step holds since
					\[
					\lambdadash(((g-1)+i')a,b)=\lambdadash((g-1)a,b)+\lambdadash(i'a,b)\equiv(g-1+i')\varepsilon_2\lambda(a,b),
					\] by using $\Z$-bilinearity of $\lambdadash$.  The other case is similar.  This completes the proof of (\ref{it:3AxiomLambda}).
					
					For (\ref{it:3AxiomMu}) we can proceed as above, reducing to the case that $j\in I$ since 
					\[
					\mudashj(-\varepsilon(j)a+ja+b)=\varepsilon(j)^2\mudash(a)+\mudashj(ja+b)+\varepsilon(j)\lambdadash(a,ja+b)=\mudashj(ja+b),
					\] using that $\varepsilon(j)$ is even.

					We will do the same induction on the number of summands in $j=i$ as in (\ref{it:3AxiomLambda}). Again, the claim is trivial for $i=0$.
					
					Assume $i=g-1$ (if $i=1-g$ the proof is analogous). We calculate
					\begin{align*}
						\mudashj((g-1)a+b)&=\mudashj((g-1)a)+\mudashj(b)+\lambdadash((g-1)a,b)\\
						&=\mudashj(ga)+\mudashj(-a)+\lambdadash(ga,-a)+\mudashj(b)+(g-1)\varepsilon_2(\lambda(a,b)) \\
						&=\mudashj(a)+\mudashj(a)+\lambdadash(ga,a)+\mudashj(b)+(g-1)\varepsilon_2(\lambda(a,b)) \\
						&= (g-1)\varepsilon_2((\lambda(a,a)+\lambda(a,b))+\mudashj(b),
					\end{align*}
					
					where the first and second equalities follow from \Cref{prop:mudashj_lamdadashj_well-defined}, the third from (\ref{it:propertiesMuLambda1}) and that $I/IJ$ is a $\Z/2$-vector space, and the fourth by (\ref{it:propertiesMuLambda3}) and that $I/IJ$ is a $\Z/2$-vector space.  This proves the base case.
					
					For the induction step assume $i=(g-1)+i'$ for a non-zero element $i'\in I$ which has fewer summands then $i$.  We obtain that
					\begin{align*}
						\mudashj(ia+b)&=\mudashj(i'a+b)+\mudashj((g-1)a)+\lambdadash_J(j'a+b,(g-1)a)\\
						&=\mudashj(i'a+b)+\mudashj(ga)+\mudashj(a)+\lambdadash(ga,a)+\lambdadash(j'a,(g-1)a)+\lambdadash(b,(g-1)a)\\
						&=\mudashj(i'a+b)+2\mudashj(a)+\lambdadash(ga,a)+\lambdadash(i'a,(g-1)a)+\lambdadash(b,(g-1)a)\\
						&=\mudashj(i'a+b)+(g-1)\varepsilon_2(\lambda(a,a))+\lambdadash(i'a,(g-1)a)+(g-1)\varepsilon_2(\lambda(a,b)),
					\end{align*}
					where the first and second equality follow from \Cref{prop:mudashj_lamdadashj_well-defined}, the third follows from (\ref{it:propertiesMuLambda1}) and (\ref{it:propertiesMuLambda2}), and the last follows from (\ref{it:3AxiomLambda}). Showing $\lambdadash(i'a,(g-1)a)= 0$ will complete the proof of the claim, which we do in the following sequence of equations.
					\[
					\lambdadash(i'a,(g-1)a)=i'\varepsilon_2(\lambda(a,(g-1)a))=i'\varepsilon_2(g\lambda(a,a))-\lambda(a,a))=0.
					\]
					This completes the proof of (\ref{it:3AxiomMu}).
				\end{proof}

				\subsection{The forms \texorpdfstring{$\mudashj$}{mudashj} and \texorpdfstring{$\lambdadash_J$}{lambdadashj} in relation to normal 1-types}
				
				Ideally we would like $\mudashj$ and $\lambdadash_J$ to be forms on the group $\pi_2(X)\otimes_{\Lambda} \Z/2$.  This will not be possible, but if we make some assumptions on our inputs we can push these forms to be well-defined on $\pi_2(X)\otimes_{\Lambda} \Z/2$. 
				\begin{prop}\label{prop:formsDescendTo}
					Fix a normal $1$-type $\xi=\xi(\pi,w_1^{\pi},w_2^{\pi})$, let $(Y,u_Y)$ be a (potentially empty) $\xi$-3-manifold and let $(X,u_X^{\pi})$ be a $(B\pi,w_1^{\pi})$-filling.  Then the following hold.
					
					\begin{enumerate}
						\item\label{it:indepentOn1}
						Let $c\in \pi_2(X)$ be s-characteristic.  Then $\mudashj(c)$ only depends on the image $[c]\in \pi_2(X)\otimes_{\Lambda}\Z/2$.
						\item\label{it:indepentOn2} Let $\omega\in \pi_2\otimes_{\Lambda}\Z/2$ such that $h^{\Lambda}_2(\omega)=0$. Then $\lambdadash_J(-,\omega)\colon \pi_2(X)\otimes_{\Lambda}\Z/2\to I/IJ$ is well-defined.
					\end{enumerate}
				\end{prop}

				\begin{proof}
					
					For (\ref{it:indepentOn1}) note that different lifts $c\in \pi_2(X)$ of elements in $\pi_2(X)\otimes_\Lambda\Z/2$ differ by elements in $J\pi_2(X)$. Let $\beta\in \pi_2(X)$ and $j\in J$. Then
					\[
					\mudash(c+j\beta)=\mudash(c)+j-\varepsilon(j)\varepsilon_2((\lambda(\beta,\beta))+\lambda(\beta,c))=\mudash(c)
					\]
					where the first equality is by (\ref{it:3AxiomMu}) of \Cref{lem:PropertiesOfMuLambdaIJ} and the second holds since $c$ is s-characteristic, implying that $\varepsilon_2((\lambda(\beta,\beta))+\lambda(\beta,c))=0$.
					
					For (\ref{it:indepentOn2}) choose a lift $\ol{\omega}\in\pi_2(X)$ of $\omega$ and an element $\alpha\in \pi_2(X)\otimes_{\Lambda}\Z$ along with its lift $\ol{\alpha}\in \pi_2(X)$. We need to show that the value $\lambdadash_J(\ol{\alpha},\ol{\omega})$ does not depend on the choices above. Different choices of $\ol{\omega},\ol{\alpha}$ differ by $j\gamma$ and $j'\beta$ for $\beta\in\ker(h_2^{\Lambda})$, $\gamma\in\pi_2(X)$ and $j,j'\in J$.  We obtain that
					\[
					\lambdadash(\alpha+j\beta,\omega+j'\gamma)=\lambdadash(\alpha,\omega)+\lambdadash(i\beta,\omega)+\lambdadash(\alpha,i'\gamma)+\lambdadash(i\beta,i'\gamma),
					\]
					and using (\ref{it:3AxiomLambda}) of \Cref{lem:PropertiesOfMuLambdaIJ} we get $\lambdadash(i\beta,\omega)=i\varepsilon_2(\lambda(\beta,\omega))$, $\lambdadash(\alpha,i'\gamma)=i'\varepsilon_2(\lambda(\alpha,\gamma))$ and $\lambdadash(i\beta,i'\gamma)=i'\varepsilon_2(\lambda(\alpha,i'\gamma))$, all of which vanish in $I/IJ$.  Since $\omega$ and $\gamma$ vanish in homology, we get that $\varepsilon_2(\lambda(\beta,\omega))=0$.  An analogous method proves that the other two $\lambdadash$ terms vanish.
				\end{proof}
				
				An immediate consequence of (\ref{it:indepentOn2}) of \Cref{prop:formsDescendTo} is that we have a well-defined form
				\begin{equation}\label{eq:lambdadashzeroform}
					\lambdadash_J^{X}\colon H_3(\pi;\Z/2)\times H_3(\pi;\Z/2)\to I/IJ,
				\end{equation} 
				using the exact sequence from \cref{sec:SerreExSeq}.  We now aim to show that this form is actually the zero form.
				
				\begin{prop}\label{prop:lambdacheckpizero}
					The map $\lambdadash_J^{X}\colon H_3(\pi;\Z/2)\times H_3(\pi;\Z/2)\to I/IJ$ from (\ref{eq:lambdadashzeroform}) is zero.
				\end{prop}
				
				\begin{remark}\label{rmk:lambdadash_on_homology}
					An immediate consequence of \Cref{prop:lambdacheckpizero} is that 
					\[
					\lambdadash_J^X: H_3(\pi;\Z/2)\times \im(h_2^{\Z/2})\to I/IJ
					\]
					is well-defined, i.e.\ for $z\in H_3(\pi;\Z/2)$ and $c\in \pi_2(X)$, $\lambdadash_J(\mathfrak{f}(z),c)$ only depends on the homology class $[c]\in H_2(X;\Z/2)$.
				\end{remark}
				
				We postpone the proof of this Proposition until we have established the following lemma which might be of independent interest.

				\begin{lemma}\label{lem:collapsetensoredDown}
					Let $V_{k,\ell}:=\ker(h_2^{\Lambda})\subset \pi_2(X\#k\CP^2\#\ell\ol{\CP^2})\otimes_{\Lambda}\Z/2$ for non-negative integers $k,\ell$. Then the map induced by collapsing the extra connected-summands
					\[
					f_{\text{coll}}\colon V_{k,\ell}\to V_{0,0}
					\] is an isomorphism.
				\end{lemma}
				\begin{proof} 
					
					We show the result for $k=1,\ell=0$; the general case is a straightforward induction. 
					
					Note first that the map $f_{\text{coll}}$ is surjective. Since $\pi_2(X\#\CP^2)\cong \pi_2(X)\oplus \Lambda$ we get that each $a\in V_{1,0}$ can be written as a sum $a=b+c$, $b\in V_{0,0}$, $c=\gamma[\CP^1]$ for $\gamma\in \Lambda$, such that $\varepsilon_2(\gamma)=0$.  The element $a\otimes 1\in \pi_2(X\#\CP^{2})\otimes_{\Lambda}\Z/2$ satisfies
					\[
					a\otimes 1=b\otimes 1+ \gamma[\CP^1]\otimes 1=b\otimes 1+[\CP^1]\otimes \varepsilon_2(\gamma)=b\otimes 1,
					\]
					which proves the injectivity of $f_{\text{coll}}$. 
				\end{proof}
				
				\begin{lemma}\label{lem:disjoint2complexes}
					There are integers $k,\ell$ such that there are two $2$-complexes $K_1,K_2$ {\it disjointly} mapped into $X\#_k\CP^2\#_l\ol{\CP}^2$, such that under the collapse map $X\#_k\CP^2\#_{\ell}\ol{\CP}^2\to X$ the $2$-complexes $K_1,K_2$ map to the $2$-skeleton of $X$.  
				\end{lemma}
				\begin{proof}
					Let $K_1\hookrightarrow X$ be any $2$-skeleton of $X$. Let $K_2^{(1)}$ be a wedge of circles embedded in $X$ disjoint to $K_1$. To make $K_2^{(1)}$ into a $2$-skeleton, we need to attach some number of $2$-cells to $K_2^{(1)}$. We can attach each such 2-cell such that it intersects $K_1$ in a finite number of points in the $2$-cells of $K_1$. For each such intersection we attach $\CP^2\#\CP^2\#\ol{\CP^2}\cong S^2\times S^2\#\CP^2$ and tube the intersection into the $S^2\times S^2$ to get rid of it (i.e.\ perform the Norman trick \cite{norman_69}). Iterating this process over all $2$-cells finishes the proof.    
				\end{proof}
				
				\begin{proof}[Proof of \cref{prop:lambdacheckpizero}]
					Let $k,\ell$ be integers and $K_1,K_2$ be two complexes as in \Cref{lem:disjoint2complexes}. Let $X_{k,\ell}:=X\# k\CP^2\# \ell\ol{\CP^2}$ and let $a,b\in H_3(\pi;\Z/2)$. From the map of fibrations
					\[
					\begin{tikzcd}
						\wt{K_i}\ar[r]\ar[d]&K_i\ar[r]\ar[d]& B\pi\ar[d]\\
						\wt{X_{k,\ell}}\ar[r]&X_{k,\ell}\ar[r]& B\pi
					\end{tikzcd}
					\]
					we get a diagram involving the maps introduced in \cref{sec:SerreExSeq}
					\[
					\begin{tikzcd}
						H_3(\pi;\Z/2)\ar[r,"\mathfrak{f}^{K_1}"]\ar[d,equal]&\pi_2(K_1)\otimes_{\Lambda}\Z/2\ar[d]\\
						H_3(\pi;\Z/2)\ar[r,"\mathfrak{f}^{X_{k,\ell}}"]\ar[d,equals]&\pi_2(X_{k,\ell})\otimes_{\Lambda}\Z/2\\
						H_3(\pi;\Z/2)\ar[r,"\mathfrak{f}^{K_2}"]&\pi_2(K_2)\otimes_{\Lambda}\Z/2\ar[u]
					\end{tikzcd}
					\]
					The above diagram together with \Cref{lem:disjoint2complexes} shows that in $\pi_2(X_{k,\ell})$ we can make $\mathfrak{f}^{X_{k,\ell}}(a)$ and $\mathfrak{f}^{X_{k,\ell}}(b)$ disjoint, and so 
					\[
					\lambdadash_J^{X\# k\CP^2\# \ell\ol{\CP^2}}(a,b)=0.
					\]
					
					Now the collapse map $X\# k\CP^2\# \ell\ol{\CP^2}\to X$ induces an isomorphism on the kernels of the Hurewicz homomorphisms $h_2^{\Lambda}$ (\cref{lem:collapsetensoredDown}) and thus on the images of $\mathfrak{f}^{X\# k\CP^2\# \ell\ol{\CP^2}}$ and $\mathfrak{f}^X$, which shows that $\mathfrak{f}^{X\# k\CP^2\# \ell\ol{\CP^2}}(a)$ and $\mathfrak{f}^{X\# k\CP^2\# \ell\ol{\CP^2}}(b)$ can both be represented by spheres in $X$ and so $\lambdadash_J^{X}(a,b)=\lambdadash_J^{X\# k\CP^2\# \ell\ol{\CP^2}}(a,b)=0$.
				\end{proof}
				
				\subsection{Secondary intersections}\label{sec:km}
				
				Since we will need it in \Cref{lem:mu_vanish_group_ring}, we briefly review the Kervaire-Milnor invariant.  The Kervaire-Milnor invariant was originally defined in \cite[Definition 10.8A]{Freedman-Quinn}, generalising the embedding obstruction first investigated in \cite{kervaire_milnor} (see also \cite{tristram}).  For the history of this secondary intersection invariant we direct the reader to \cite[Section 3]{kasprowski_powell_ray_teichner} where it is given a full treatment.  
				
				\begin{definition}\cite[Definition 1.4]{kasprowski_powell_ray_teichner}\label{def:kervaire_milnor}
					Let $b\subset X$ be an immersed sphere in a smooth 4-manifold $X$ and assume that $\mudash(b)=0\in I$, the augmentation ideal of $\Lambda$.  Then the \emph{Kervaire-Milnor invariant} $\km(b)\in\Z/2$ is zero if and only if there exists a collection of Whitney discs $D_i$ pairing the double points of $b$ such that the interiors of the $D_i$ are all disjoint from $b$, after potentially changing $b$ by a regular homotopy.
				\end{definition}
				
				For the statement of the following theorem the reader should recall \Cref{def:sr_characteristic} of an s-characteristic class and r-characteristic class.
				
				\begin{theorem}{\cite{stong_1994}\cite[Theorem 1.6]{kasprowski_powell_ray_teichner}}\label{thm:kervaire_milnor}
					Let $b$ be an $s$-characteristic immersed sphere in a smooth 4-manifold $X$ and assume that $\mudash(b)=0$, and hence denote the collection of Whitney discs pairing up the double points of $b$ by $D_i$.  Further assume that $b$ has a (potentially unframed) algebraic dual sphere.  If $b$ is $r$-characteristic then
					\[
					\km(b) \equiv \sum_i \vert \interior D_i\cap b\vert \pmod{2},
					\]
					otherwise $\km(b)=0$.
				\end{theorem}
				
				We now demonstrate that the Kervaire-Milnor invariant is the stable obstruction to changing an immersed sphere $b$ by regular homotopy to an embedding.
				
				\begin{lem}\label{lem:km_stable_embedding}
					Let $b$ be an immersed sphere in a 4-manifold $X$ with $\mudash(b)=0$ and $\km(b)=0$.  Then $b$ is regularly homotopic to an embedding in $X\#k(S^2\times S^2)$ for some $k\geq 0$, i.e.\ $b$ is stably regularly homotopic to an embedding.
				\end{lem}
				
				\begin{proof}
					By \Cref{def:kervaire_milnor}, we can choose a collection of framed Whitney discs $D_i$ for $b$ whose interiors are disjoint from $b$.  We wish to perform Whitney moves over all of the $D_i$ such that $b$ becomes embedded, but to do this we need to remove all of the intersections between the various Whitney discs.  Let $p\in D_i\cap D_j$.  Form the connected-sum $X\# (S^2\times S^2)$ and tube $D_i$ into $S^2\times \{\pt\}\subset S^2\times S^2$ (the standard Norman trick \cite{norman_69}).  Note that this operation preserves the discs being framed.  Now $D_i$ has a geometric dual sphere $\{\pt\}\times S^2$, and so we can remove the intersection point $p$ by tubing $D_j$ into the geometric dual (again this preserves the framing).  Repeating this process for all intersection points makes all of the Whitney discs disjointly embedded and framed; hence performing the corresponding Whitney moves takes $b$ to an embedding.
				\end{proof}
				
				We give an example, which will be used in \Cref{lem:mu_vanish_group_ring}, of an immersed sphere with non-trivial Kervaire-Milnor invariant, which is due to \cite{kervaire_milnor}.
				
				\begin{lem}\label{lem:km_example}
					Let $w\in H_2(\CP^2)$ be a generator.  Then the class $3w$ can be represented by an immersed sphere $b$ with $\mudash(b)=0$ and $\km(b)=1$.
				\end{lem}
				
				\begin{proof}
					By adding trivial cusps, it is clear that one can arrange for an immersion $b$ with $\mudash(b)=0$ representing $3w$.  One can verify that the argument given in \cite[Example 2]{kervaire_milnor} works even if we allow stabilisations (stabilisations do not change the signature), and hence $3w$ cannot be represented by an embedding, even stably.  Hence, $\km(b)=1$.
				\end{proof}

				\section{Obstruction theory for filling 3-manifolds}\label{sec:obstructions}

				\subsection{The primary obstruction}\label{subsec:primary}
				
				Let $\xi=\xi(\pi,w_1^{\pi},w_2^{\pi})$ be a normal $1$-type and let $(Y,u_Y)$ be a $\xi$-3-manifold. Recall (\cref{def:xifilling}) that we have the natural map $\xi\to \xi':=\xi(\pi,w_1^{\pi},\infty)$ and by a $(B\pi,w_1^{\pi})$-filling $(X,u_X^{\pi})$ of $Y$ we mean a null-bordism of $Y$ in $\Omega_3^{\xi}$. In other words $X$ is a $\xi'$-manifold with boundary $Y$, whose $\xi'$-structure agrees with the $\xi$ structure on the boundary.

				Recall that for a $\xi$-3-manifold $(Y,u_Y)$ there is a well-defined primary obstruction $\pri(Y,u_Y)\in H_3(\pi;\Z^{w_1^{\pi}})$. The following is a well-known result.

				\begin{lem}\label{lem:primary}
					Let $\xi=\xi(\pi,w_1^{\pi},w_2^{\pi})$ be a normal $1$-type where we allow the totally non-spin case, i.e.\ $w_2^{\pi}\in H^2(\pi;\Z/2)\cup\{\infty\}$. Let $(Y,u_Y)$ be a $\xi$-$3$-manifold. Then the primary obstruction $\pri(Y,u_Y)$ equals $(u_{Y})_*[Y] \in \ker d_2 \subset H_3(\pi;\Z^{w_1^{\pi}})$, where $[Y]\in H_3(Y;\Z^{w_1(Y)})$ is the twisted fundamental class and $d_2$ denotes the appropriate second differential in the James spectral sequence.
				\end{lem}
				
				For the proof see, for example, \cite[Proposition 3.2.1]{teichnerthesis} which can be readily adapted for the non-orientable case.
				
				\begin{prop}\label{prop:PrimaryObstruction}
					Let $\xi=\xi(\pi,w_1^{\pi},w_2^{\pi})$ be a normal $1$-type. Let $(Y,u_Y)$ be a $\xi$-3-manifold. Then $\pri(Y,u_Y)=0$ if and only if there if a $(B\pi,w_1^{\pi})$-filling $(X,u_X^{\pi}\colon X\to B\pi)$ of $(Y,u_Y)$ such that $u_X^{\pi}$ induces an isomorphism on the fundamental groups.
				\end{prop}  
				
				\begin{proof}
					
					The proof is by comparison $\xi=\xi(\pi,w_1^{\pi},w_2^{\pi})\to \xi'=\xi(\pi,w_1^{\pi},\infty)$. Naturality of the James spectral sequence gives a diagram of edge homomorphisms:
					
					\[\begin{tikzcd}
						\Omega_3^{\xi}\ar[r]\ar[dr]&\Omega_3^{\xi'}\ar[d,"\cong"]
						\\
						&H_3(\pi;\Z^{w_1^{\pi}})
					\end{tikzcd}\]
					where the vertical map is an isomorphism since the group $\Omega_*^{\xi'}$ can be calculated using the James spectral sequence 
                    \[
                    H_p(\pi;(\Omega_3^{SO})^{w_1^{\pi}})\Rightarrow \Omega_{p+q}^{\xi'},
                    \]
					and the bordism group $\Omega_*^{SO}$ vanishes in dimensions $1,2$ and $3$.
					
					This shows that $\pri(Y,u_Y)=0$ if and only if $Y$ vanishes in $\Omega_3^{\xi'}$ and we finish the proof by applying \cref{lem:surgeryBelowMiddle}.  
				\end{proof}

				\subsection{The secondary obstruction}
				
				Let $(Y,u_Y)$ be a $3$-dimensional $\xi$-manifold. Assume that $\pri(Y)$ vanishes and choose a $(B\pi,w_1^{\pi})$-filling $(X,u_X^{\pi})$. Recall that $\smallO(X,Y,u_Y,u_X^{\pi})\in H^2(X,Y;\Z/2)$ is the obstruction to extending the $\xi$-structure over $X$ rel $Y$ (see \Cref{prop:relative_obstruction_class}).  We denote by $w\in H_2(X;\Z/2)$ the Poincar\'{e} dual of $\smallO(X,Y,u_Y,u_X^{\pi})$. Using the James spectral sequence for $\xi(\pi,w_1^{\pi},w_2^{\pi})$, recall that there is now a well-defined secondary invariant $\secondary(Y,u_Y)\in H_2(\pi;\Z/2)/\im(d_2)$ (there is no outgoing differential, see \Cref{prop:JSSdifferential}).  We now define the secondary invariant geometrically using $(X,u_X^{\pi})$.  Later in \Cref{sec:differentials} we will quotient out by the choice of this filling to get a well-defined invariant that will be equal to $\secondary(Y,u_Y)$.
				
				\begin{definition}
					Let $\xi=\xi(\pi,w_1^{\pi},w_2^{\pi})$ be a normal $1$-type, let $(X,u_X^{\pi})$ be a $(B\pi,w_1^{\pi})$-filling of a $\xi$-3-manifold $(Y,u_Y)$, and let $w\in H_2(X;\Z/2)$ be the obstruction class as above.  Then we define the \emph{geometric secondary obstruction} to be
					\[
					\secG(X,Y,u_Y,u_X^{\pi}):=(u_{X}^{\pi})_*(w)\in H_2(\pi;\Z/2).
					\]
				\end{definition}
				
				In the above definition if $Y$ is empty we denote $\secG(X,\emptyset,\emptyset,u_X^{\pi})=\secG(X,u_X^{\pi})$.
				
				\begin{lem}\label{lem:sec_vanishing_implies_obstruction_spherical}
					Let $\xi=\xi(\pi,w_1^{\pi},w_2^{\pi})$ be a normal 1-type, let $(Y,u_Y)$ be a $\xi$-3-manifold with $(X,u_X^{\pi})$ a $(B\pi,w_1^{\pi})$-filling and $u_X^{\pi}$ inducing an isomorphism on fundamental groups.  Then the Poincar\'{e} dual of $\smallO(X,Y,u_X^{\pi},u_Y)$ is spherical if and only if $\secG(X,Y,u_X^{\pi},u_Y)=0$.
				\end{lem}
				
				\begin{proof}
					The result follows from considering the Serre exact sequence and the surjectivity of the map $\pi_2(X)\to \pi_2(X)\otimes_{\Lambda}\Z/2$.
				\end{proof}

				\subsection{The tertiary obstruction}

				We now use the self-intersection invariant $\mudashj\colon \pi_2(X)\otimes_{\Z} \Z/2\to I/IJ\cong H_1(\pi;\Z/2)$ from \Cref{sec:intersections} to define a geometric tertiary obstruction $\terG$.  
				
				\begin{definition}
					Let $\xi=\xi(\pi,w_1^{\pi},w_2^{\pi})$ be a normal $1$-type, let $(X,u_X^{\pi})$ be a $(B\pi,w_1^{\pi})$-filling of a $\xi$-3-manifold $(Y,u_Y)$ such that $\smallO(X,Y,u_Y,u_X^{\pi})\in H_2(X;\Z/2)$ is a spherical homology class, i.e.\ $\secG(X,Y,u_Y,u_X^{\pi})=0$. Let $c\in \pi_2(X)$ satisfy $h_2(c)=\smallO(X,Y,u_Y,u_X^{\pi})$, where $h_2$ denotes the $\Z/2$ Hurewicz map.  Then we define the \emph{geometric tertiary obstruction} to be
					\[
					\terG(X,Y,u_Y,u_X^{\pi},c):=\mudashj(c)\in H_1(\pi;\Z/2).
					\]
				\end{definition}
				This value depends on the bounding manifold $X$ as well as the choice of $c$ representing $\smallO(X,Y,u_Y,u_X^{\pi})$. To get rid of these dependencies we have to take an appropriate quotient.  We will show how to do this in the next section. 
				
				\section{Differentials and well-definedness of the obstructions}\label{sec:differentials}
				
				\subsection{The geometric differential \texorpdfstring{$\delta_2^{(4,0)}$}{delta2,4,0}}\label{subsec:delta_d2_sec}
				
				We now define the geometric differential
				
				\[\delta_2^{(4,0)}\colon H_4(\pi;\Z^{w_1^{\pi}})\to H_2(\pi;\Z/2)\] for the secondary invariant which was briefly described in \Cref{sec:intro}.  We then show that it is equal to the corresponding $d_2$ differential in the James spectral sequence.  We start with the definition.
				
				\begin{definition}\label{def:sec_delta_2}
					Let $x\in H_4(\pi;\Z^{w_1^{\pi}})$ be any element. Represent $x$ as a closed $(B\pi,w_1^{\pi})$-manifold $(M,u_M^{\pi})$ with $u_M^{\pi}$ a $\pi_1$-isomorphism (its existence follows from the proof of \Cref{prop:PrimaryObstruction}). Let $\xi=\xi(\pi,w_1^{\pi},w_2^{\pi})$ be a normal 1-type and let $w\in H_2(M;\Z/2)$ be Poincar\'{e} dual to $\smallO(M,u_M^{\pi})$.  We define
					\[
					\delta_2^{(4,0)}(x):= \secG(M,u_M^{\pi}).
					\]
				\end{definition}
				
				It is currently not clear that the above definition is well-defined, i.e.\ does not depend on the choice of representative manifold.  We deal with this now by showing that $\delta_2^{(4,0)}$ is \emph{equal} to $d_2^{(4,0)}$ from the JSS.  Since the latter is already well-defined, it follows that $\delta_2^{(4,0)}$ is also.
				
				\begin{lemma}\label{lem:d2sec}
					Let $\xi=\xi(\pi,w_1^{\pi},w_2^{\pi})$ be a normal 1-type, let $(M,u_M^{\pi})$ be a closed $(B\pi,w_1^{\pi})$-manifold and let $x = (u_M^{\pi})_*[M]$ be the image of the twisted fundamental class $[M]\in H_4(\pi;\Z^{w_1(M)})$.  Then $\delta_2^{(4,0)}(x)=d_2^{(4,0)}(x)$.
				\end{lemma}
				
				\begin{proof}
					We will prove the following string of equalities.
					\[
					d_2^{(4,0)}(x)=(u^{\pi}_{M})_*(PD(w_2(M))+PD((u_M^{\pi})^*(w_2^{\pi}))=(u_M^{\pi})_*(PD(\smallO(M,u_M^{\pi})))=\delta_2^{(4,0)}(x).
					\]
					The second equality follows from \cref{prop:relative_obstruction_class}, and the third is by \Cref{def:sec_delta_2}.  We now concern ourselves with the first equality.
					
					The Universal coefficient theorem gives an isomorphism $H_2(\pi;\Z/2)\cong\Hom(H^2(\pi;\Z/2),\Z/2)$ and hence the first equality of the above holds if and only if   
					\begin{align*}
						\langle w_2^{\pi}\smile \alpha+w_1^{\pi}\smile \Sq^1(\alpha) +&\Sq^2(\alpha),(u_M^{\pi})_*[M]\rangle =\\ 
						&=\langle\alpha, (u_M^{\pi})_*(w_2(M)\frown [M])\rangle+ \langle\alpha,(u_M^{\pi})_*((u_M^{\pi})^*w_2^{\pi}\frown [M])\rangle
					\end{align*}
					for all $\alpha\in H^2(\pi;\Z/2)$ (above we have used the description of the second differential from \Cref{prop:JSSdifferential}).
					
					We calculate
					\begin{align*}
						\langle &w_2^{\pi}\smile \alpha+w_1^{\pi}\smile \Sq^1(\alpha) +\Sq^2(\alpha),(u_M^{\pi})_*[M] \rangle = \\ 
						&=\langle w_2^{\pi}\smile \alpha,(u_M^{\pi})_*[M]\rangle+\langle w_1^{\pi}\smile\Sq^1(\alpha)\rangle+\langle\alpha^2,(u_M^{\pi})_*[M]\rangle\\
						&=\langle (u_M^{\pi})^*\alpha,(u_M^{\pi})^*w_2^{\pi} \frown [M]\rangle+\langle (u_M^{\pi})^*(w_1^{\pi})\smile \Sq^1((u_M^{\pi})^*(\alpha)),[M]\rangle + \langle(u_M^{\pi})^*(\alpha)^2,[M]\rangle\\
						&=\langle\alpha,(u_M^{\pi})_*((u_M^{\pi})^*w_2^{\pi}\frown [M])\rangle+\langle w_1(M)\smile \Sq^1((u_M^{\pi})^*(\alpha))\rangle+\langle(u_M^{\pi})^*(\alpha)\smile w_2(M),[M]\rangle\\
						&=\langle\alpha,(u_M^{\pi})_*((u_M^{\pi})^*w_2^{\pi}\frown[M])\rangle+\langle\Sq^1(\Sq^1((u_M^{\pi})^*(\alpha)))[M]\rangle+\langle(u_M^{\pi})^*(\alpha), w_2(M)\frown [M]\rangle\\
						&=\langle\alpha,(u_M^{\pi})_*((u_M^{\pi})^*w_2^{\pi}\frown [M])\rangle+\langle\alpha, (u_M^{\pi})_*(w_2(M)\frown [M])\rangle.
					\end{align*}
					
					Here the first equality follows by linearity of the Kronecker pairing, the second by the cup-cap formula, and the third equality used that the normal second Stiefel-Whitney class $w_2$ equals the second tangential Wu class $v^t_2$ and so acts a characteristic element.\footnote{In general we have for tangential Wu classes of $n$-dimensional manifolds $W$ that
						\(v_i^{t}x=\Sq^i(x)\) for \(x\in H^{n-i}(W;\Z/2)\).} The fourth equality follows from $w_1=w_1^{t}=v_1^{t}$, and the last equality uses the Adem relation $\Sq^1\Sq^1=0$.
				\end{proof}
				
				\begin{prop}\label{prop:SecXYuY}
					Let $\xi=\xi(\pi,w_1^{\pi},w_2^{\pi})$ be a normal 1-type. Let $(X,u_X^{\pi})$ be a $(B\pi,w_1^{\pi})$-filling of $(Y,u_Y)$, with $\pi_1(u_X^{\pi})$ an isomorphism. Let $w\in H_2(X;\Z/2)$ be Poincar\'{e}-Lefschetz dual to the obstruction class $\smallO(X,Y,u_Y,u_X^{\pi})$.  Then 
					\[
					H_2(\pi;\Z/2)/\im d_2^{(4,0)}\ni \secondary(Y,u_Y)=[\secG(X,Y,u_Y,u_X^{\pi})]\in H_2(\pi;\Z/2)/\im\delta_2^{(4,0)}
					\]
					and hence $[\secG(X,Y,u_Y,u_X^{\pi})]\in H_2(\pi;\Z/2)/\im\delta_2^{(4,0)}$ does not depend on the choice of $(B\pi,w_1^{\pi})$-filling.
				\end{prop}
				
				\begin{proof}
					This is due to the second-named author \cite[Proposition 3.2.3]{teichnerthesis}.
				\end{proof}
				
				We now finish this subsection by showing the following improvement of \Cref{lem:sec_vanishing_implies_obstruction_spherical}.
				
				\begin{lem}\label{lem:sec_vanishing_better}
					Let $\xi=\xi(\pi,w_1^{\pi},w_2^{\pi})$ be a normal 1-type.  Let $(Y,u_Y)$ be a $\xi$-3-manifold with $\pri(Y,u_Y)=0$.  Then there exists $(X,u_{X}^{\pi})$ a $(B\pi,w_1^{\pi})$-filling of $Y$ with the Poincar\'{e}-Lefschetz dual to $\smallO(X,Y,u_Y,u_X^{\pi})$ spherical if and only if $\secondary(Y,u_Y)=0$.
				\end{lem}
				
				\begin{proof}
					The `only if' direction follows from \Cref{lem:sec_vanishing_implies_obstruction_spherical}.

					Conversely, assume $\secondary(Y,u_Y)=0$.  By \Cref{lem:primary} and \Cref{prop:SecXYuY} there exists a~$(B\pi,w_1^{\pi})$-filling $(X',u_{X'}^{\pi})$ of $(Y,u_Y)$ such that $\secG(X',Y,u_Y,u_{X'}^{\pi})$ is in the image of $\delta_2^{(4,0)}$.  This implies that there exists a closed~$(B\pi,w_1^{\pi})$-manifold $(M,u_M^{\pi})$ such that 
					\[
					\delta_2^{(4,0)}(u_{M*}^{\pi}([M]))=\secG(X',Y,u_Y,u_{X'}^{\pi}).
					\] 
                    Form the connected sum $X'':=M\#X'$ and  $u^{\pi}_{X''}:=(u_M^{\pi}\# u_X^{\pi})\colon M\#X'\to B\pi$. By \Cref{lem:doubleObstructionClass} the Poincar\'{e} dual of $\smallO(X'',Y,u_Y,u_{X''}^{\pi})$ can be written as a sum $w_{X''}=w_{M}+w_{X'}$ of the Poincar\'{e} duals of the obstructions of the summands.  We have 
					\[
					u_{M*}^{\pi}(w_M)+(u_{X'}^{\pi})_*(w_{X'})=0\in H_2(\pi;\Z/2),
					\]
                    Realise $w_{X''}$ by some surface $\Sigma$ in $X''$. Using \Cref{lem:surgeryBelowMiddle} we can surger $X''$ to turn the classifying map $u^{\pi}_{X''}$ to a $\pi_1$-isomorphism, while preserving $\xi$ structure in the complement of $\Sigma$. Call the resulting $(B\pi,w_1^{\pi})$ manifold $X$. Now $\Sigma$ is a surface in $X$ representing $\secG(X,Y,u_Y,u_X^{\pi})=0$. Since $u_X^{\pi}(\Sigma)=0$, applying \Cref{lem:sec_vanishing_implies_obstruction_spherical}  implies that the class $[\Sigma]$ is spherical. 
				\end{proof}

				\subsection{The geometric differential \texorpdfstring{$\delta_2^{(3,1)}$}{delta2,3,1}}\label{subsec:delta_d2}
				
				We now set out to define the geometric differential $\delta_2^{(3,1)}$ for the tertiary invariant that was briefly described in \Cref{sec:intro}.  We then show that it is equal to the corresponding $d_2$ differential in the JSS.  Due to an abundance of superscripts, we will drop the $(3,1)$ superscript in this section, and so the unaltered differential will simply read as $\delta_2$.
				
				\begin{proposition}\label{prop:delata_2properties}
					Let $\xi=\xi(\pi,w_1^{\pi},w_2^{\pi})$ be a normal $1$-type, let $(Y,u_Y)$ be a (possibly empty) $\xi$-$3$-manifold, and let $X$ be its $(B\pi,w_1^{\pi})$-filling. Let $c\in \pi_2(X)$ be s-characteristic and let $w:=[c]\in H_2(X;\Z/2)$. Then we have the following 
					\begin{enumerate}
						\item\label{it:delta2_1} 
						The map
						\[
						\delta_2^{X,c}:=\mudashj(-)+\lambdadash_J(-,c):\ker(h_2\colon\pi_2(X)\to H_2(X;\Z/2))\to I/IJ
						\] 
						only depends on $[c]\in H_2(X;\Z/2)$, and we denote the map above as $\delta_2^{X,w}$. 
						\item\label{it:delta2_2}
						The map $\delta_2^{X,w}$ factors as
						\[
						\delta_2^{X,w}=\mudashj(-)+\lambdadash_J(-,w)\colon \ker(h_2^{\Lambda}\colon\pi_2(X)\otimes_{\Lambda}\Z/2 \to H_2(X;\Z/2))\to I/IJ,
						\] 
						and so the following map is well-defined 
						\[
						\delta_2^{X,w}=\mudashj(-)+\lambdadash_J(-,w)\colon H_3(\pi;\Z/2)\to I/IJ.
						\]
						\item\label{it:delta2_3} The map $\delta_2^{X,w}\colon H_3(\pi;\Z/2)\to I/IJ$ only depends on  the $\CP^2$-stable equivalence class of the manifold $X$, i.e.\
						\[
						\delta_2^{X\#\CP^2,w+[\CP^1]}=\delta_2^{X,w}=\delta_2^{X\#\ol{\CP^2},w+[\ol{\CP^1}]}.                 
						\]
					\end{enumerate}
				\end{proposition}

				\begin{proof}
					(\ref{it:delta2_1}) For a fixed homology class $w\in H_2(X;\Z/2)$, different lifts $c\in \pi_2(X)$ differ by elements $\beta\in\pi_2(X)$ such that $h_2(\beta)=0$. Let $\alpha\in \ker(h_2)\subset \pi_2(X)$.
					\[
					\mudashj(\alpha)+\lambdadash_J(\alpha,c+\beta)=\mudashj(\alpha)+\lambdadash_J(\alpha,c)+\lambdadash_J(\alpha,\beta)
					\]
					Now from \Cref{prop:lambdacheckpizero} we get $\lambdadash_J(\alpha,\beta)=0$. 
					
					For (\ref{it:delta2_2}) let $\alpha\in \ker(h_2^{\Lambda})\subset \pi_2(X)\otimes_{\Lambda}\Z/2$ and let $\ol{\alpha}\in\pi_2(X)$ be any lift. Note that a different lift differs by $j\beta\in J\ker(h_2^{\Lambda})$ by \Cref{lem:tensorIsos}(\ref{it:H1group2}).  Observe that
					\begin{align*}          
						\mudashJ(\ol{\alpha}+j\beta)+\lambdadash_J(\ol{\alpha}+&j\beta,c)=\\
						&=\mudashJ(\ol{\alpha})+(j-\varepsilon(j))(\varepsilon_2(\lambda(\beta,\beta))+\varepsilon_2(\lambda(\beta,\ol{\alpha})))+\lambdadash(\ol{\alpha},c)+\lambdadash(j\beta,c)\\
						&=\mudashJ(\ol{\alpha})+\lambdadash(\ol{\alpha},c)+\lambdadash(j\beta,c)\\
						&=\mudashJ(\ol{\alpha})+\lambdadash(\ol{\alpha},c)+(j-\varepsilon(j))\varepsilon_2\lambda(\beta,c)\\
						&=\mu_J(\ol{\alpha})+\lambdadash_J(\ol{\alpha},c),
					\end{align*}
					where we used \Cref{lem:PropertiesOfMuLambdaIJ}(\ref{it:3AxiomLambda}) and (\ref{it:3AxiomMu}) and the fact that $\varepsilon_2(\lambda(-,-))$ is zero if at least one of the arguments is null-homologous.
					
					For (\ref{it:delta2_3}), we do the proof for $\CP^2$. The proof for $\ol{\CP}^2$ is analogous. Let $z\in H_3(\pi;\Z/2)$. We use the commutative diagram
					\[
					\begin{tikzcd}
						H_3(\pi;\Z/2)\ar[r,"\mathfrak{f}^{X\#\CP^2}"]\ar[d,equal]&\pi_2(X\#\CP^2)\otimes_{\Lambda}\Z/2\ar[d]\ar[r,"h_2^{\Lambda}"] &H_2(X\#\CP^2;\Z/2)\ar[d]\\
						H_3(\pi;\Z/2)\ar[r,"\mathfrak{f}^{X}"]&\pi_2(X)\otimes_{\Lambda}\Z/2\ar[r,"h_2^{\Lambda}"] &H_2(X;\Z/2),
					\end{tikzcd}
					\]
					induced by the collapse map $X\#\CP^2\to X$. From \Cref{lem:collapsetensoredDown} we find that the middle vertical map induces an isomorphism on the images of $\mathfrak{f}^{X\#\CP^{2}}$ and $\mathfrak{f}^{X}$. It follows that $\mathfrak{f}^{X\#\CP^{2}}(z)$ can be represented by a sphere in $X$, which proves the claim. 
				\end{proof}
				
				\begin{proposition}\label{prop:delta2_closed}
					Let $\xi=\xi(\pi,w_1^{\pi},w_2^{\pi})$ be a normal $1$-type, let $(Y,u_Y)$ be a (possibly empty) $\xi$-$3$-manifold, let $X$ be its $(B\pi,w_1^{\pi})$-filling and assume the obstruction class $\smallO(X,Y,u_Y,u_X^{\pi})$ is spherical. Then there exists a closed $(B\pi,w_1^{\pi})$-manifold $(M,u_M^{\pi})$, with $\pi_1(u_M^{\pi})$ an isomorphism, $u_{M*}^{\pi}[M]=0$ for the twisted fundamental class $[M]\in H_4(M;\Z^{w_1^{\pi}})$ and $w_M:=\PD(\smallO(M,u_{M}^{\pi}))$ is spherical such that 
					\[
					\delta_2^{X,w_X}=\delta_2^{M,w_M}\colon H_3(\pi;\Z/2)\to H_1(\pi;\Z/2),
					\]
					where $w_X$ is the Poincar\'{e} dual of $\smallO(X,Y,u_Y,u_X^{\pi})$.
				\end{proposition}
				\begin{proof}
					If $X$ is already closed, work instead with punctured $X$. Now assume $X$ is a compact 4-manifold with non-empty boundary. Consider the inclusion into the double $i\colon X\hookrightarrow DX$ and the following commutative diagram.
					\[
					\begin{tikzcd}
						H_3(\pi;\Z/2)\ar[r,"\mathfrak{f}^X"]\ar[d,equal]&\pi_2(X)\otimes_{\Lambda}\Z/2\ar[d]\\
						H_3(\pi;\Z/2)\ar[r,"\mathfrak{f}^{DX}"]&\pi_2(DX)\otimes_{\Lambda}\Z/2
					\end{tikzcd}
					\]
					Let $z\in H_3(\pi;\Z/2)$. \Cref{lem:doubleObstructionClass} implies that the obstruction to a $\xi$-structure on $DX$ is given by $i_*w_X+Dw_X\in H_2(DX;\Z/2)$, where $Dw_X$ is the image of $w_X$ in the other half of $DX$. The commutative diagram above shows that the image of $a$ in $\pi_2(DX)\otimes_{\Lambda}\Z/2$ can be realised by a sphere entirely in $X$.  It immediately follows that
                    \begin{align*}
                        \delta_2^{DX,i_*w+Dw}(z)&=\mudashj^{DX}(\mathfrak{f}^X(z))+\lambdadash^{DX}_J(\mathfrak{f}^X(z),i_*w+Dw)\\
                        &=\mudashj^{X}(\mathfrak{f}^X(z))+\lambdadash_J^X(\mathfrak{f}^X(z),i_*w)=\delta_2^{X,w}(z).\qedhere
                    \end{align*}
				\end{proof}
				
				The proof of the following proposition will comprise the next ten pages.
				
				\begin{proposition}\label{prop:delta2isd2}
					Let $\xi=\xi(\pi,w_1^{\pi},w_2^{\pi})$ be a normal $1$-type, let $X$ be a $(B\pi,w_1^{\pi})$-$4$-manifold with possibly non-empty boundary, a $\xi$-3-manifold $(Y,u_Y)$. Let $w\in H_2(X;\Z/2)$ be the Poincar\'{e} dual of the obstruction class $\smallO(X,Y,u_Y,u_X^{\pi})\in H^2(X,Y;\Z/2)$. Then we have that~$\delta_2^{X,w}\colon H_3(\pi;\Z/2)\to I/IJ$ defined in \Cref{prop:delata_2properties}(\ref{it:delta2_2}) satisfies
					\[\delta^{X,w}_2 = d_2\colon H_3(\pi;\Z/2)\to H_1(\pi;\Z/2).\]
					In particular, it does not depend on the choice of manifold $X$ (only on the normal 1-type $\xi$).
				\end{proposition}
				
				We will first solve the $\Spin$ and $\Pin^+$ case of \Cref{prop:delta2isd2}, which we will then use to solve the general case. Recall that for $\xi=\xi(\pi,w_1^{\pi},0)$ the differential $d_2\colon H_3(\pi;\Z/2)\to H_1(\pi;\Z/2)$ is the Hom dual of $\Sq^1(-)\smile w_1$, since the second Steenrod square vanishes in this case (see \Cref{sec:jss}).
				
				\begin{proposition}\label{prop:spinCasedelta2d2}
					Let $\xi=\xi(\pi,w_1^{\pi},0)$ be a $\Spin$ or $\Pin^+$ normal $1$-type (depending on whether $w_1^{\pi}$ is zero or not). Let $(M,u_M^{\pi})$ be a $(B\pi,w_1^{\pi})$-manifold with $\secG(M,u_M^{\pi})=0$. Then \[\delta_2^{M,w}=(w_1^{\pi}\smile \Sq^1(-))^*=d_2\colon H_3(\pi;\Z/2)\to H_1(\pi;\Z/2).\]
				\end{proposition}
				
				The proof of \Cref{prop:spinCasedelta2d2}  relies on first proving it in a few special cases that we will then compare to.
				
				\begin{lemma}\label{lem:spinCasedelta2d2}
					Let
					
					\begin{enumerate}
						\item[(i)] $\pi=\Z/2$ and $\xi=\xi(\pi,0,0)$, or
						\item[(ii)] $\pi=\Z/2$ and $\xi=\xi(\pi,x_1,0)$ for $x_1\in H^1(\Z/2;\Z/2)$ the generator, or
						\item[(iii)]$\pi=\Z/2\times\Z/2$  and $\xi=\xi(\pi,x_1+x_2,0)$ for $x_1,x_2$ the standard generators of $ H^1(\Z/2\times\Z/2;\Z/2)$.
					\end{enumerate}

					Then for a closed $(B\pi,w_1^{\pi})$-manifold $u_M^{\pi}\colon M\to B\pi$, with $\secG(M,u_M^{\pi})=0$ the conclusion of \cref{prop:spinCasedelta2d2} holds.
				\end{lemma}
				
				We will use the following construction multiple times during the proof, so we record it here.

				\begin{construction}\label{constr:SurgeryOnManifold}
					For a manifold $M$ with $\pi_1(M)\cong \Z/2*\Z/2$ we denote by $(M)_*$ the result of a surgery on the commutator of the generators of $\pi_1(M)$. If $M$ has a normal $\Pin^+$ structure away from some characteristic surface $w$ this surgery can be preformed disjointly from $w$ and in a $\Pin^+$-preserving way away from $w$.
				\end{construction}

				\begin{proof}[Proof of \Cref{lem:spinCasedelta2d2}]

					By \Cref{prop:delata_2properties}(\ref{it:delta2_3}) we only need to consider each $\CP^2$-stable class $(B\pi,w_1^{\pi})$. By \cite{Kreck99} there is $\vert H_4(\pi;\Z^{w_1^{\pi}})/\Aut(\pi,w_1^{\pi})\vert$ many $\CP^2$-stable classes\footnote{Allowing connected summing with both $\CP^2$ and $\ol{\CP^2}$.} of manifolds with fundamental group $\pi$ and orientation class $w_1^{\pi}$. The secondary invariant does not depend on the stable $\CP^2$-class.
					
					Case (i), $\xi=\xi(\Z/2,0,0)$: There is exactly $\vert H_4(B\Z/2)/\Aut(\Z/2)\vert=1$ such stable class, represented by a manifold commonly denoted by $F$; this is a closed spin manifold with $\pi_1(F)\cong \Z/2$ and universal cover $S^2\times S^2$.  We will use some topological properties of $F$ listed for example in \cite[End of Section 4]{Ham-kreck-finite}. The deck transformation on the universal cover of $F$ is 
					\begin{align*}
						S^2\times S^2&\rightarrow S^2\times S^2\\
						(x,y)&\mapsto (-x,r(y)),
					\end{align*}
					The $k$-invariant $k_F\in H^3(\Z/2;\pi_2(F))$ is  $(1,0)$ under the obvious isomorphism $\pi_2(S^2\times S^2)\cong \pi_2(S^2)\times \pi_2(S^2)\cong \Z\left<a,b\right>$.  By \Cref{prop:f_all_Agree}, the evaluation map of the reduced $k$-invariant $k^r\in H^3(\Z/2;\pi_2(F)\otimes_{\Lambda}\Z)$ give the first map of the Serre exact sequence for $F$ (see \Cref{sec:SerreExSeq}).
					\[
					\begin{tikzcd}
						H_3(\pi;\Z/2)\ar[d,"\cong"]\ar[r,"ev(k_F)"]&\pi_2(F)\otimes_\Lambda\Z/2\ar[r,"h_2^{\Lambda}"]\ar[d,"\cong"]&H_2(F;\Z/2)\ar[d,"\cong"]\\
						\Z/2\ar[r,"{(1,0)}"]&\Z/2\oplus \Z/2\ar[r,"{\proj_2}"]&\Z/2
					\end{tikzcd}
					\]
					Since $F$ is Spin we need to check that $\delta_2^{F,0}(\alpha)=\mudashj(\alpha)+\lambdadash_j(\alpha,0)$ is zero for any $\alpha\in \ker(h_2^{\Lambda})$.  The kernel $\ker(h_2^{\Lambda})$ is generated by the class $\alpha=[*\times S^2]$.  From the deck transformation we see that $\mudashj(\alpha)=0$ and so $\delta_2^{F,0}=0$. The expression $(w_1^{\pi}\smile \Sq^1(-))^*$ vanishes as well.
						
						Case (ii), $\xi=\xi(\Z/2,w_1^{\Z/2}=x_1,0)$: Since $H_4(\Z/2;\Z^{w_1^{\Z/2}})\cong \Z/2$ we need to consider two stable $\CP^2$-classes. These are $\RP^4$ and $\RP^2\times S^2$ for the non-trivial and the trivial-class respectively, as can be seen from their classifying maps.
						
						The obstruction to these manifolds being {\it normally} $\Pin^+$ equals the obstruction to them being tangentially $\Pin^{-}$. Since $w^t_1(\RP^4)^2+w^t_2(\RP^4)=x_1^2\in H^2(\RP^2;\Z/2)$, its Poincar\'{e} dual represents the non-trivial class in $H_2(\RP^2;\Z/2)$. Since the classifying map $u_{\RP^{4}}^{\pi}\colon\RP^4\to B\Z/2\simeq \RP^{\infty}$ is the inclusion the secondary invariant $\secG(\RP^4,u_{\RP^{4}}^{\pi})$ is non-zero.

						On the other hand, the manifold $\RP^2\times S^2$ is normally $\Pin^{+}$ and so the obstruction class $\smallO(\RP^2\times S^2,u_{\RP^2\times S^2}^{\pi})$ can be represented by a null homologous element and so $\secG(\RP^2\times S^2,u_{\RP^2\times S^2}^{\pi})$ vanishes. To show that  $\delta_2^{\RP^2\times S^2,0}=(x_1\smile \Sq^1(-))^*$ we first calculate the right-hand side. The map $H^1(\Z/2;\Z/2)\to H^3(\Z/2;\Z/2)$ given by $x_1\smile \Sq^1(-)$ is an isomorphism and so is its dual.
						To calculate the left-hand side we need to examine the relevant portion of the Serre exact sequence, which is
						\[
						H_3(\Z/2;\Z/2)\xrightarrow{\mathfrak{f}^{\RP^2\times S^2}} \pi_2(\RP^2\times S^2)\otimes \Z/2\xrightarrow{h_2^{\Lambda}} H_2(\RP^2\times S^2).
						\]
						The kernel of the second map is $\Z/2$ generated by the universal cover $a\colon S^2\to \RP^2$ postcomposed with the inclusion $\iota_{\RP^2}\colon \RP^2\to \RP^2\times S^2$.  We show that $\mudashj(\iota_{\RP^2}\circ a)\in H_1(\Z/2;\Z/2)$ is the nontrivial element as follows. Let $f_{\text{fold}}\colon S^2\to D^2$ denote the fold map which sends the equator of $S^2$ into the boundary of $D^2$ and the north and south-pole to $0\in D^2$. The following map is homotopic to $\iota_{\RP^2}\circ a$.
						\begin{equation}\label{eq:foldProof}
							a'\colon S^2\xrightarrow{\id_{S^2}\times f_{\text{fold}}} S^2\times D^2\xrightarrow{a\times \id_{D^2}} \RP^2\times D^2\hookrightarrow \RP^2\times S^2
						\end{equation}
						
						The image of this map has a single double point, namely it sends the south-pole and the north-pole to the same element. Since the line which joins the north-pole to the south-pole is the generator of $\pi_1(\RP^2)$ we conclude that the double-point loop associated to this intersection is the generator of $\pi_1(\RP^2\times S^2)$ and hence $ \mudashj(a')=x_1\neq 0$.
						
						Case (iii), $\xi=\xi(\Z/2\times \Z/2,w_1^{\Z/2\times \Z/2}=x_1+x_2,0)$: 
						
						We use the model $B(\Z/2\times\Z/2)\simeq B\Z/2\times B\Z/2\simeq \RP^{\infty}\times \RP^{\infty}$ to calculate $H_4(\Z/2\times \Z/2;\Z^{x_1+x_2})$. The universal cover of $\RP^\infty$ is the contractible space $S^\infty$. The singular chain complex for $S^{\infty}$ is $C_i(S^{\infty})\cong \Z[\Z/2]$ with the differential  $d_i:C_i(S^{\infty})\to C_{i-1}(S^{\infty})$ given by $t+(-1)^{i}$ which can be used to calculate that \[H_4(\Z/2\times \Z/2;\Z^{x_1+x_2})\cong (\Z/2)^3\cong \Z/2\left<e_0\times e_4,e_2\times e_2, e_4\times e_0\right>\] and  \[H_4(\Z/2\times \Z/2;\Z^{x_1+x_2})/\Aut_{x_1+x_2}(\Z/2\times \Z/2)\cong (\Z/2)^2\cong \Z/2\left<e_0\times e_4,e_2\times e_2\right>,\]
						
						Where we take the cellular decomposition of $S^\infty$ to have two cells in each dimension and the universal covering projection $S^\infty\to \RP^2$ to be cellular.
						We now aim to find manifolds to represent the classes  $\left<e_0\times e_4,e_2\times e_2\right>$.  The class $e_2\times e_2$ is given by $\RP^2\times \RP^2$. The class $e_0\times e_4$ is given by $(S^2\times \RP^2\#\RP^4)_*$, since this class is given by $S^2\times \RP^2\#\RP^4\to B\Z/2\times B\Z/2$ and cobordism preserves this image.  Finally, the zero class is given by $(S^2\times\RP^2\#S^2\times\RP^2)_*$.  For all of the above manifolds we use the obvious map to $\RP^{\infty}\times \RP^{\infty}$.  Since these computations will be lengthy, we split into sub cases.
						
						Case (iii.a), $X=((S^2\times\RP^2)\#(S^2\times\RP^2))_*$: The Kirby diagram of $X$ is given in \Cref{fig:rp2_x_s2_rp2_x_s2_kirby diagram}. Here the secondary invariant of $X$ vanishes, since $X$ is normally $\Pin^{+}$ and so the obstruction class $\smallO(X,u_X^{\pi})$ vanishes.  
						
						\begin{figure}[ht]
							\centering
							\includegraphics[width=0.5\linewidth]{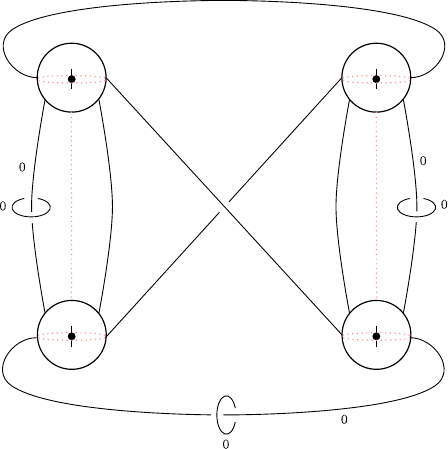}
							\caption{A Kirby diagram for the manifold $X=((\RP^2\times S^2)\#(\RP^2\times S^2))_*$.  The notation for the non-orientable 1-handles is due to Akbulut \cite[Sections 1.1, 1.5]{akbulut}}\label{fig:rp2_x_s2_rp2_x_s2_kirby diagram}
						\end{figure}
						
						We will define some elements $\alpha_1,\alpha_2\in\pi_2(X)\otimes_{\Lambda}\Z/2$ that we will need later. Let $C\cong D^2$ be the core of the $S^2\times D^2$ which was attached during the surgery to form $X$. We define annuli $A_t,A_s\colon S^1\times D^2\hookrightarrow X$ such that each boundary of the annuli maps to $\partial C$, and define $\alpha_1=C\cup A_t\cup \ol{C}$ and $\alpha_2=C\cup A_s\cup \ol{C}$.  We describe the annulus $A_t$ using an inertial regular homotopy of a push-off $\gamma$ of $C$ (see \Cref{fig:annulus_picture}).  The annulus $A_s$ is defined analogously.
						
						Now we make the following claim about the properties of $X$ that we will need.    
						
						\begin{clm}\label{clm:computation}
							Let $X=((\RP^2\times S^2)\#(\RP^2\times S^2))_*$ and let $\pi:=\pi_1(X)$.  Then the Serre exact sequence (Appendix \ref{sec:SerreExSeq}) sequence associated to $X\to B(\Z/2\times \Z/2)$ is
							\begin{equation}\label{eq:claim_sequence}
								\begin{tikzcd}[cramped]
									H_3(X;\Z/2) \arrow[r,"u_{X*}^{\pi}"] & H_3(\pi;\Z/2) \arrow[r,"\mathfrak{f}"] & \pi_2(X)\otimes \Z/2 \arrow[r,"h_2^{\Lambda}"] & H_2(X;\Z/2) \arrow[r,"u_{X*}^{\pi}"] & H_2(\pi;\Z/2)    
								\end{tikzcd}
							\end{equation}
							Then the following holds. The groups involved are given by 
							
							\begin{enumerate}
								\item \label{it:1ExProp}
								\[
								H_k(X;\Z/2)\cong \begin{cases}
									\Z/2\ \text{if}\ k=0,4 \\ (\Z/2)^2\ \text{if}\ k=1,3 \\ (\Z/2)^6\ \text{if}\ k=2 \\ 0\ \text{else} 
								\end{cases}
								H_k(\pi;\Z/2)\cong \begin{cases}
									\Z/2\ \text{if}\ k=0 \\ (\Z/2)^2\ \text{if}\ k=1\\ (\Z/2)^3\ \text{if}\ k=2 \\ (\Z/2)^4\ \text{if}\ k=3
								\end{cases},
								\]
								\item\label{it:2ExProp} and $\pi_2(X)\otimes \Z/2 \cong (\Z/2)^7$.
							\end{enumerate}
							We collect information about the generators:
							\begin{enumerate}
								\setItemnumber{3}
								\item\label{it:3ExProp} The group $H_3(\pi;\Z/2)$: generated by $e_3\times e_0, e_2\times e_1, e_1\times e_2$, and $e_0\times e_3$, where $e_i$ denotes the $i$-th cell in $\RP^{\infty}$.
								\item\label{it:4ExProp} The group $H_2(\pi;\Z/2)$: generated by $e_2\times e_0, e_1\times e_1$ and $e_0\times e_2$.
								\item\label{it:5ExProp} The group $\pi_2(X)\otimes_{\Lambda} \Z/2$: generated by the two covering maps $S^2\to \RP^2$, the two $\{\pt\}\times S^2$ and the belt sphere to the surgery, as well as two  classes $\alpha_1$ and $\alpha_2$ defined above.
								\item\label{it:6ExProp} The group $H_2(X;\Z/2)$: generated by the two copies of $\RP^2$, the two copies of $S^2$, the belt sphere to the surgery, as well as another class $R$.
							\end{enumerate}
							and the maps:
							\begin{enumerate}
								\setItemnumber{7}
								\item\label{it:7ExProp} The map $u_{X*}^{\pi}$ on $H_3(X;\Z/2)$: this is the zero map.
								\item\label{it:8ExProp} The map $\mathfrak{f}$: this sends $e_3\times e_0$ and $e_0\times e_3$ to the two covering maps $S^2\to \RP^2$, respectively, and sends $e_2\times e_1$ and $e_1\times e_2$ to $\alpha_1$ and $\alpha_2$, respectively.
								\item\label{it:9ExProp} The Hurewicz map $h_2^{\Lambda}$: this sends the two copies of $S^2$ to the two copies of $S^2$, and the belt sphere of the surgery to the belt sphere of the surgery.
								\item\label{it:10ExProp} The map $u_{X*}^{\pi}$ on $H_2(X;\Z/2)$: this sends the subgroup $(\Z/2)^3$ generated by the two copies of $\RP^2$ and some class $R$ injectively to $H_2(\pi;\Z/2)$.

							\end{enumerate}
							This finishes the statement of the claim.  We delay its proof until after we finish the current proof.
						\end{clm}

\begin{figure}[h!]
    \centering
    \includegraphics[width=0.87\linewidth]{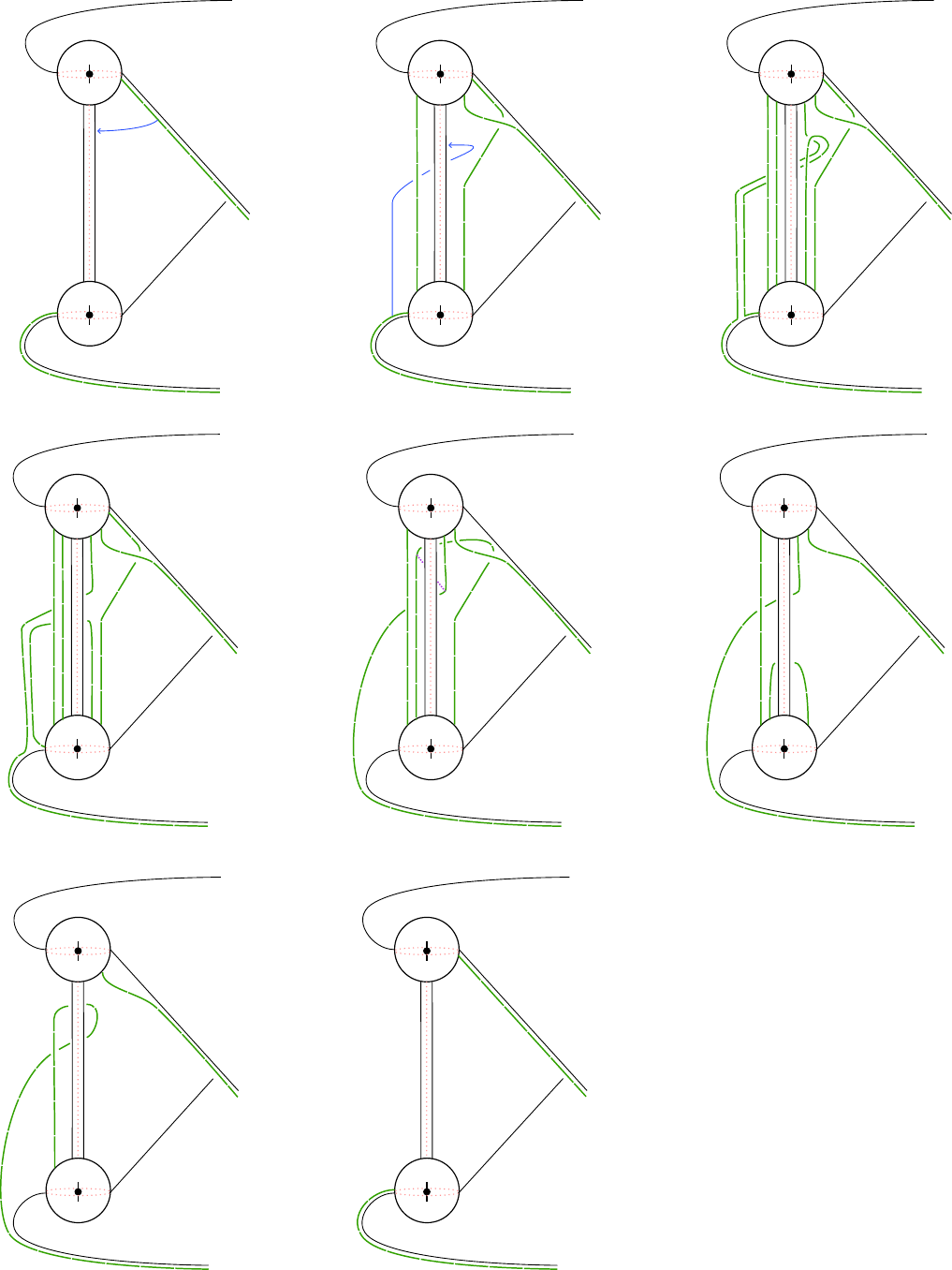}
    \caption{The annulus $A_t$ described using a regular homotopy of $\gamma$, shown in sporadically dashed green, whereas the relevant parts of the Kirby diagram are shown in solid black.  The guide for the single intersection point is shown in dashed purple.  The handle slides are guided by the blue arrows.  Only relevant portions of the Kirby diagram from \Cref{fig:rp2_x_s2_rp2_x_s2_kirby diagram} are shown.}\label{fig:annulus_picture}
    \end{figure}
    
    Since $w=\PD(X,u_X^{\pi})=0$, we have that $\delta_2^{X,0}=\mudashj\circ \mathfrak{f}^X$. By \Cref{clm:computation}(\ref{it:8ExProp}) we need to calculate the self intersection of elements in the image of $\mathfrak{f}^X$. One can compute from \Cref{fig:annulus_picture} that $\alpha_1$ and $\alpha_2$ have a single double-point and furthermore that $\mudashj(\alpha_1)=t$ and $\mudashj(\alpha_2)=s$, where $s,t\in H_1(B\Z/2\times B\Z/2;\Z/2)$ are generators dual to $x_1$ and $x_2$ respectively.  The two covering maps $S^2\to \RP^2$ have $\mudashj =t$ and $s$, respectively (the proof is identical to the earlier calculation for $\RP^2\times S^2$ in this proof, see map \ref{eq:foldProof}).
    
    The map $d_2$ is dual to $(w_1^{\pi}\smile \Sq^1(-))\colon H^1(B\Z/2\times B\Z/2;\Z/2)\to H^3(B\Z/2\times \Z/2;\Z/2)$, which is given on the generators by \begin{align*}
        x_1\mapsto x_1^3+x_1^2x_2\\
        x_2\mapsto x_1x_2^2+x_2^{3}
    \end{align*}
    since $w_1^{\pi}=x_1+x_2$. Now $d_2$ is 
    \begin{align*}
        x_1^3\mapsto t=(x_1)^*\\
        x_1^2x_2\mapsto t=(x_1)^*\\
        x_2^3\mapsto s=(x_2)^*\\
        x_1x_2^2\mapsto s=(x_2)^*
    \end{align*}
    which agrees with $\delta_2^{X,0}$.
						
    Case (iii, b), the manifold $X=\RP^2\times \RP^2$ is \emph{not} normally $\Pin^+$ and the normal $w_2$ is given by $x_1x_2$. This element does not vanish in $H_2(B\Z/2\times B\Z/2)$ and so 
    \[
        \secG(\RP^2\times \RP^2;\Z/2)=x_1x_2\neq 0.
        \]
        Case (iii, c), $X=((S^2\times\RP^2)\#\RP^4)_*$:  A calculation shows that $S^2\times\RP^2\#\RP^4$ has the normal second Stiefel-Whitney class $w_2= x_2^2$. Let $w\in H_2(S^2\times\RP^2\#\RP^4;\Z/2)$ be the Poincar\'{e} dual of $x_2^2$.  By preforming surgery like in \cref{constr:SurgeryOnManifold} we obtain   $(S^2\times\RP^2\#\RP^4)_*$, which still has $w$ as a characteristic element.  Furthermore, an easy Mayer-Vietoris argument shows that the class $w$ is still-nontrivial and we obtain the following commutative diagram.
        \begin{equation*}
            \begin{tikzcd}
                H_2(S^2\times\RP^2\#\RP^4;\Z/2)\ar[r]\ar[d,tail]&H_2(B(\Z/2*\Z/2);\Z/2)\ar[d,tail]\\
                H_2((S^2\times\RP^2\#\RP^4)_*;\Z/2)\ar[r]&H_2(B(\Z/2\times\Z/2);\Z/2)
            \end{tikzcd}
        \end{equation*}
        The left vertical map is actually given by a zigzag of Mayer-Vietoris sequences.  The commutativity of the above diagram then gives that $\sec$ is non-zero, since $w$ maps non-trivially into $H_2(\Z/2\times \Z/2;\Z/2)$.
    \end{proof}

    \begin{proof}[Proof of \Cref{clm:computation}]
        We first show (\ref{it:1ExProp}) Calculating the (generators for the) homology of $B(\Z/2\oplus \Z/2)\simeq \RP^{\infty}\times \RP^{\infty}$ is a basic exercise in algebraic topology, proving (\ref{it:3ExProp}) and (\ref{it:4ExProp}).
        
        Let $X_0$ obtained by removing the attaching region of the surgery 
        \[
        (S^2\times\RP^2)\#(S^2\times\RP^2)\rightsquigarrow (S^2\times\RP^2\#S^2\times\RP^2)_*=:X,
        \] 
        which was performed away from some fixed representative of $w$.  Then $X_0$ admits inclusions $S^2\times\RP^2\#S^2\times\RP^2\hookleftarrow X_0\hookrightarrow X$. 
        
        Consider the following Mayer-Vietoris sequence for the first inclusion (with $\Z/2$-coefficients omitted).
        \begin{equation*}
            \begin{tikzcd}[cramped]
                0\ar[r]&H_2(S^1\times S^2)\ar[r,"\text{belt}"]&H_2(X_0)\ar[r]& H_2(S^2\times\RP^2\#S^2\times\RP^2)\ar[r]&0
            \end{tikzcd}
        \end{equation*}
        The last map is zero because the next map is an injection $H_1(S^1\times S^2)\to H_1(S^1\times D^2)\oplus H_1(X_0)$. The image of the map labelled ``belt'' gives a class in $X_0$ and in turn in $X$ is given by the belt sphere of the handle.  We now state the corresponding sequence for the second inclusion.
        \begin{equation*}
            \begin{tikzcd}[cramped]
                0\ar[r]& H_2(S^1\times S^2)\ar[r]&H_2(X_0)\oplus H_2(D^2\times S^2)\ar[r]& H_2(X)\ar[r]&H_1(S^1\times S^2)\ar[r]&0
            \end{tikzcd}
        \end{equation*}
        Here the last map is zero because the generator of $H_1(S^1\times S^2)$ maps to a commutator in $X_0$ and hence vanishes in homology.
        
        Via the two sequences above we conclude that we have a (non-canonical) isomorphism $H_2(X;\Z/2)\cong (\Z/2)^6\cong  H_2(S^2\times\RP^2\#S^2\times\RP^2;\Z/2)\oplus (\Z/2)^2$ where the $(\Z/2)^{2}$ class is generated by the belt sphere as defined above and some class $R$ and $H_2(S^2\times\RP^2\#S^2\times\RP^2;\Z/2)\cong (\Z/2)^4$ is generated by two copies of $S^2$ and two copies of $\RP^2$. This proves (\ref{it:1ExProp}) and (\ref{it:6ExProp}).
        
        Consider the Serre exact sequence (\ref{eq:claim_sequence}).  The last map $(u_X^{\pi})_*\colon H_2(X;\Z/2)\to H_2(\pi;\Z/2)$ is surjective; its kernel is $(\Z/2)^3$, generated by the two embedded $2$-spheres and the belt sphere of the attached handle, since spherical classes necessarily vanish in $H_2(B\pi)$. This proves (\ref{it:10ExProp}).  We show that the first map $(u_X^{\pi})_*$ is zero. The map $S^2\times\RP^2\#S^2\times\RP^2\to B\Z/2\times B\Z/2$ factors through the 2-skeleton $(B\Z/2\times B\Z/2)_{(2)}$ and so does the map $X\to B\Z/2\times B\Z/2$, since, up to homotopy equivalence, only a 2-handle was attached, which proves (\ref{it:7ExProp}).
        
        Since  $H_3(B\Z/2\times B\Z/2;\Z/2)\cong (\Z/2)^4$ we have that $\pi_2(X)\otimes \Z/2\cong (\Z/2)^7$, proving (\ref{it:2ExProp}). Since the image of $h_2^{\Lambda}$ is $(\Z/2)^{3}$ its kernel must be $(\Z/2)^4$. To identify the generators, we consider the map $\mathfrak{f}$. By \Cref{prop:f_all_Agree}, an alternate description of $\mathfrak{f}$ is the following: a model for $B\pi$ can be constructed by attaching $n$-cells for $n\geq 3$ to $X$ to kill all higher homotopy groups.  In this model any $z\in H_3(\pi;\Z/2)$ is given by a linear combination of 3-cells, which attach along 2-spheres in $X$.  Then $\mathfrak{f}(z)$ is this linear combination of attaching maps.  We use this description now.

        Factor the classifying map of $X$ as 
        \[
        X=(S^2\times\RP^2\#S^2\times\RP^2)_*\to (S^2\times \RP^2)\vee (S^2\times\RP^2)\cup f_2\to \RP^2\vee \RP^2\cup f_2\xrightarrow{\text{incl.}} \RP^{\infty}\times\RP^{\infty},
        \]
        where $f_2$ denotes a 2-cell that attaches along the commutator.  Since $\RP^2\vee \RP^2\cup f_2$ is precisely the 2-skeleton of $\RP^{\infty}\times\RP^{\infty}$, the map $\mathfrak{f}$ is easy to see on that space; it is given by the attaching maps of the 3-cells in $\RP^{\infty}\times \RP^{\infty}$.  We obtain the diagram, which commutes by naturality of Serre spectral sequences (see \Cref{sec:SerreExSeq}).
            \[
        \begin{tikzcd}
            \pi_2(X)\otimes_{\Lambda} \Z/2 \arrow[r] & \pi_2(\RP^2\vee\RP^2\cup f_2)\otimes_{\Lambda}\Z/2 \\
            H_3(\Z/2\times\Z/2;\Z/2)\arrow[r,equals]\arrow[u,"\mathfrak{f}_1"]& H_3(\Z/2\times\Z/2;\Z/2)  \arrow[u,"\mathfrak{f}_2"] 
        \end{tikzcd}
        \]
        By the characterization above, the map $\mathfrak{f}_2$ sends $e_3\times e_0$ to $2 e_2\times e_0$, $e_0\times e_3$ to $e_0\times 2e_2$, $e_2\times e_1$ to $2e_1\times e_1 \cup e_2\times 2e_0$, and $e_1\times e_2$ to $2e_0\times e_2 \cup e_1\times 2e_1$.  Note that the two covering maps  $(S^2\to \RP^2)\in \pi_2(X)\otimes\Z/2$ map to $2e_2\times e_0$ and $e_0\times 2e_2$ in $\pi_2(\RP^2\vee\RP^2\cup f_2)\otimes_{\Lambda}\Z/2$  The maps $\mathfrak{f}_2$ is an isomorphism and so $\mathfrak{f}_1(e_3\times e_0)$ and $\mathfrak{f}_2(e_0\times e_3)$ map to the covering maps, since the covering maps vanish in $\Z/2$-homology and so lift to $H_3(\Z/2\times \Z/2;\Z/2)$.
        
    We now show that $\alpha_1=\mathfrak{f}_1(e_2\times e_1)$ (showing $\alpha_2=\mathfrak{f}_1(e_1\times e_2)$ is the same).  Like above, it suffices to show that $\alpha_1$ maps to $\partial (e_2\times e_1)\in \pi_2(\RP^2\vee\RP^2\cup f_2)\otimes \Z/2)$.  The two discs $C$ and $\ol{C}$ both map to $f_2$, so we only have to conclude that the annulus $A_1$ maps to $2e_2\times e_0$.  Recall that the annulus was formed by two handle slides, both over the 2-handle whose core maps to $e_2\times e_0$.  This implies that the trace of the annulus maps to $2e_2\times e_0$, and hence $\alpha_1$ maps to $\partial(e_2\times e_1)$. This proves (\ref{it:8ExProp}). Using the above and (\ref{it:10ExProp}) we show (\ref{it:5ExProp}). Combining  (\ref{it:5ExProp}) and (\ref{it:6ExProp}) proves (\ref{it:9ExProp}).
\end{proof}

\begin{proof}[Proof of \Cref{prop:spinCasedelta2d2}]
    Let $\xi=\xi(\pi,w_1^{\pi},0)$ be given and let $(M,u_M^{\pi})$ be a $(B\pi,w_1^{\pi})$-manifold with $\secG(M,u_M^{\pi})=0$.  We show that we can always set up a comparison with one of the 1-types considered in \Cref{lem:spinCasedelta2d2} that will suffice to describe the differential.
    
    Let $z\in H_3(\pi;\Z/2)$, let $\ol{\mathfrak{f}^M(z)}\colon S^2\looparrowright X$ be an immersed 2-sphere such that  $[\ol{\mathfrak{f}^M(z)}]\in \pi_2(M)$ is a lift of $\mathfrak{f}^M(z)$, and let $[c\colon S^2\looparrowright M]\in \pi_2(M)$ be a lift of $\PD (w_2(M))$, which exists since the secondary invariant of $M$ vanishes.  Assume towards a contradiction that $x:=\delta_2^{M,w}(z)-d_2(z)\in H_1(\pi;\Z/2)$ is non-zero.  Then there exists a map $\varphi_x\colon \pi\to \Z/2$ which maps $x$ to the generator of $H_1(\Z/2;\Z/2)$.  We now split into cases.  
    \begin{enumerate}
        \item If $w_1^\pi=0$ then surger out the kernel of $\varphi_x$ to form a $(B\Z/2,w_1^{\pi}=0)$-manifold $M'$.
        \item If $w_1^{\pi}(x)\neq 0$ then surger out the kernel of $w_1^{\pi}$ to form a $(B\Z/2,w_1^{\pi}\neq 0)$-manifold $M'$.
        \item If $w_1^{\pi}\neq 0$ but $w_1^{\pi}(x)=0$ then surger out a kernel of $w_1^{\pi}\oplus \varphi_x\colon\pi\to \Z/2\oplus\Z/2 $ to form a $(B(\Z/2\times \Z/2),w_1^{\pi}=(1,0))$-manifold $M'$.
    \end{enumerate}
    In all cases perform these surgeries away from $\ol{\mathfrak{f}^M(z)}$ and $c$, which can be done by general position.  Denote this $5$-dimensional surgery trace by $W$.  The class $c$ still represents the second Stiefel-Whitney class of $M'$ since we may assume that our surgeries preserve $\Spin/\Pin^{+}$-structures in the complement of $c$.
    
    Let $i_M\colon M\to W$, $i_{M'}\colon M'\to W$ denote the two inclusions.  Consider the following diagram where the rows are formed from the Serre exact sequence for $M$, $W$ and $M'$, respectively, which commutes by naturality.
    \[
    \begin{tikzcd}
        H_3(\pi;\Z/2)\ar[r,"\mathfrak{f}^M"]\ar[d,"(i_M)_*"]&\pi_2(M)\otimes_{\Lambda}\Z/2\ar[r]\ar[d,"(i_M)_*"]&H_2(M;\Z/2)\ar[d,"(i_M)_*"]\\
        H_3(\pi_1(W);\Z/2)\ar[r,"\mathfrak{f}^W"]&\pi_2(W)\otimes_{\Lambda}\Z/2\ar[r]&H_2(W;\Z/2)\\
        H_3(\pi_1(M');\Z/2)\ar[r,"\mathfrak{f}^{M'}"]\ar[u,"\cong", "(i_{M'})_*"']&\pi_2(M')\otimes_{\Lambda}\Z/2\ar[r]\ar[u,"(i_{M'})_*"']&H_2(M';\Z/2)\ar[u,"(i_{M'})_*"']
    \end{tikzcd}\]
    
						Note that by construction $\pi_1(\varphi_x)=\pi_1(i_{M'}^{-1}\circ (i_M))$ and so $\varphi_x=i_{M'}^{-1}\circ i_M\colon H_3(\pi;\Z/2)\to H_3(\pi_1(M');\Z/2)$.
						Since we performed the surgeries away from $\ol{\mathfrak{f}^M(z)}$, this immersed sphere can be viewed as an immersed sphere in $M'$, but for clarity we denote $y:=[\ol{\mathfrak{f}(z)}]\in \pi_2(M')\otimes_{\Lambda} \Z/2$, which then satisfies $(i_M')_*(y) = (i_M)_*(\mathfrak{f}(z))$. It is not necessarily true that $y$ and $\mathfrak{f}^{M'}(\varphi_x(z))$ agree. The relative handle chain complex for $(W,M')$ is
						\[ 
						C_*(W,M')=\begin{cases}
							\Z^n\quad *=3\\
							0\quad \text{else}
						\end{cases}
						\]
						
						for some non-negative integer $n$ and hence \[\ker \left((i_{M'})_* \colon H_2(M';\Z/2)\to H_2(W;\Z/2)\right) \cong H_3(W,M';\Z/2)\cong (\Z/2)^n,\] generated by the belt $3$-discs $B_{\ell}$ of the $1$-handles, and their images in $H_2(M')$ are the boundaries $\partial B_{\ell}$, for ${\ell}=1,\dots, n$.  By an application of the snake lemma \[\ker ((i_{M'})_* \colon \pi_2(M')\otimes_{\Lambda} \Z/2\to \pi_2(W)\otimes_{\Lambda} \Z/2) \cong (\Z/2)^n,\] again generated the $\partial B_{\ell}$ and so there are integers $k_{\ell}$ such that 
                        \[
                        y+ \sum_{\ell=1,\dots, n}k_{\ell}\partial B_{\ell} = \mathfrak{f}^{M'}((\varphi_x)_*(z)).
                        \]  
                        We obtain that
						\begin{align*}(\varphi_x)_*d_2(z) &= d_2((\varphi_x)_*(z)) \\ 
							&= \delta_2^{M',w_2(M')}((\varphi_x)_*(z)) =\mudashj(\mathfrak{f}^{M'}((\varphi_x)_*(z))+\lambdadash_J(\mathfrak{f}^{M'}((\varphi_x)_*(z)),c) \\
							&= \mudashj(y+ \sum_{\ell=1,\dots, n}k_{\ell}\partial B_{\ell})+\lambdadash_J(y+ \sum_{\ell=1,\dots, n}k_{\ell}\partial B_{\ell},c) = \mudashj(y)+\lambdadash_J(y,c) \\
							&= (\varphi_x)_*[\mudash_J(\mathfrak{f}^M(z))+\lambdadash_J(\mathfrak{f}^M(z),c)]=(\varphi_x)_*\delta_2^{M,w}(z),
						\end{align*}
						where the first equality comes from naturality of the JSS; the second by \Cref{lem:spinCasedelta2d2}; the third by the definition of $\delta_2^{M',w_2(M')}$; the fourth from the above argument; the fifth since the $\partial B_{\ell}$ are all framed, embedded and disjoint from $c$; the sixth is from the definition of the element $y$; and the seventh is by the definition of $\delta_2^{M,w}$.
						
						This implies that $(\varphi_x)_*(\delta_2^{M,w}(z)-d_2(z))=0$, which contradicts the definition of $\varphi_x$.  Hence $x=\delta_2^{M,w}(z)-d_2(z)=0$.
					\end{proof}

					\begin{lemma}\label{lem:CapLambdaCheck}
						For any closed $(B\pi,w_1^{\pi})$-manifold $(M,u_M)$ with a $\pi_1$-isomorphism $u_M^{\pi}\colon M\to B\pi$, with $u^{\pi}_{M*}([M])=0$ such that $\PD(w_2(M))$ can be represented by an immersed sphere,  we have an equality in $I/IJ$
						\[
						a\frown z=\lambdadash_J(\mathfrak{f}^M(z),\PD(u_M^{\pi*}a))
						\] 
						for any $z\in H_3(\pi;\Z/2)$ and any $a\in H^2(\pi;\Z/2)$. 
					\end{lemma}
					
					\begin{proof}
						First, note that the right-hand side of the above equation is well-defined, see \Cref{rmk:lambdadash_on_homology}.  
						
						Form a model for $B\pi$ by attaching $k$-cells to $M$ where $k\geq 3$. This way the map $\mathfrak{f}^M$ has a geometric interpretation as a certain boundary map (see \Cref{def:mathfrakf2}). Let  $z\in H_3(\pi;\Z/2)$ and $a\in H^2(\pi;\Z/2)$ be any elements. Represent $z$ as a map from a closed 3-dimensional manifold $z\colon W^3\to B\pi$, which is always possible by Thom \cite{Thom1954}. Assume that $W$ has a handlebody decomposition with a single 3-handle and denote by $W^0$ to be the manifold $W$ with this handle removed. Homotope the map~$z$ such that $W^0$ lands in $M$. 
						
						Choose a model for $D:=\PD((u^{\pi}_M)^*(a))\in H_2(M;\Z/2)$ given by an immersed sphere in $M$. This is possible since by assumption the lower horizontal map of the following commutative diagram is zero
						\[
						\begin{tikzcd}
							H_2(M;\Z/2)\ar[d,"u_{M*}^{\pi}"]&\ar[l,"\PD"] H^2(M;\Z/2)\\
							H_2(\pi;\Z/2)&\ar[l,"{\frown u_M^{\pi*}[M]}"] H^2(\pi;\Z/2)\ar[u,"u_M^{\pi*}"]
						\end{tikzcd}
						\]
						and therefore $\PD((u^{\pi}_M)^*(a))$ is spherical.  The transverse intersection of $W^0$ and $D$ is given by a $1$-dimensional manifold, which we denote by $A$, that is some collection of $\Z/2$-nullhomologous circles and proper arcs $\{a_i\}_{i=1,\dots,k}$ in $W^0$. Close up all the arcs in $A$ by arcs $\{\wt{a_i}\}_{i=1,\dots,k}$ in $\partial W^0$ (there is a unique way to do this since different choices are homologous through the top cell of $W$) to obtain a manifold $\ol{A}$ that, by construction, is homologous to~$a\frown z$.
						
						By construction the element $\partial W^0\rightarrow M$ is $\mathfrak{f}^M(a)\in \pi_2(M)\otimes\Z/2$. Hence the intersection $P:=\partial W^0\pitchfork D $, when counted with groups elements, gives the value $\lambdadash_J(\mathfrak{f}^M(z),a)$. Note that $\partial A=P$ and using this we have a natural decomposition of $P$ into sets of size two, $P=\{p_1,q_1\}\cup\{p_2,q_2\}\cup\dots\{p_k,q_k\}$, where each pair $\{p_i,q_i\}$ is connected by two arcs $a_i\subseteq A,\wt{a_i}\subseteq \ol{A}\setminus\ol{{\ol{A}\setminus A}}$.  What is left to show is that the $1$-dimensional manifold $\ol{A}$ is homologous to the double point loops of $\partial W^0$ and $D$. 
						
						Fix a basepoint $x_0\in M$. For every pair $\{p_i,q_i\}\subset P$ we fix paths $\gamma_i$, $\wt{\gamma_i}$ from $x_0$ to  $p_i$ in $D$ and $\partial W^0$, respectively.  For $p\in P$ let $\theta(p)$ be the $\pi_1$ element associated to the point $p$.  Then
						\[
						\lambdadash_J(\mathfrak{f}(z),a)=\sum_{p\in P}\theta(p)= \sum_{i=1,\dots,k} \theta(p_i) +\theta(q_i) =\sum_{i=1,\dots,k}\wt{\gamma_i}\gamma_i+\wt{\gamma_{i}}\wt{a_i}a_{i}\gamma_i=\sum_{i=1,\dots,k}\wt{a_i}a_{i}=\ol{A},
						\]
						where the first equality is by the definition of $\lambdadash_J$, the second is by reindexing, the third is by the observation above, the fourth holds since we are working in homology, and the last holds by the decomposition of $\ol{A}$ into arcs $\widetilde{a}_i a_i$ and null-homologous circles.
					\end{proof}
					
					We can now finally prove \Cref{prop:delta2isd2}.
					
					\begin{proof}[Proof of \Cref{prop:delta2isd2}]
						We need to show that 
						\[
						\delta_2^{X,w}(z)=d_2(z)
						\] 
						for all $z\in H_3(\pi;\Z/2)$.  By \Cref{prop:JSSdifferential}, $d_2$ is dual to 
						\[
						d^2(-)= -\frown w_2^{\pi} + Sq^1(-)\frown w_1^{\pi}.
						\]
						By \Cref{prop:delta2_closed} there exists a closed $(B\pi,w_1^{\pi})$-manifold $M$ such that $\delta_2^{X,w}=\delta_2^{M,w_M}$, $(u_M^{\pi})_*([M])=0$ and the Poincar\'{e} dual of the obstruction class $w_M=\PD(\smallO(M,u_M^{\pi}))$ spherical.  Since the diagram
						\[\begin{tikzcd}
							H_2(M;\Z/2)\ar[d,"(u_{M}^{\pi})_*"]&&\ar[ll,"\PD"] H^2(M;\Z/2)\\
							H_2(\pi;\Z/2)&&\ar[ll,"{\frown (u_M^{\pi}[M]})^*"] H^2(\pi;\Z/2)\ar[u,"(u_M^{\pi})^*"]
						\end{tikzcd}\]
						commutes and $(u_M^{\pi})_*([M])=0$, we have that $\PD((u^{\pi}_M)^*(w_2^{\pi}))$ is spherical.  By decomposing $w_M$ as $w_M=\PD(w_2(M)+(u^{\pi}_M)^*(w_2^{\pi}))$ we see that $\PD(w_2(M))$ is also spherical.  Furthermore, by \Cref{lem:sphericalCharacteristic} we have that both $\PD(w_2(M))$ and $\PD((u^{\pi}_M)^*(w_2^{\pi}))$ are $s$-characteristic.  These facts give us the following description of $\delta_2^{M,w_M}$.
						\begin{align*}
							\delta_2^{M,w_M}(z) &= \mudashj(\mathfrak{f}(z)) + \lambdadash(\mathfrak{f}(z), \PD(w_2(M)+(u_M^{\pi})^*(w_2^{\pi})))\\
							&=\delta_2^{M,\PD(w_2(M))}(z)+\lambdadash(\mathfrak{f}(z),(u_M^{\pi})^*(w_2^{\pi}))\\
							&= (w_1^{\pi}\smile \Sq^1(-))^*(z) + \lambdadash(\mathfrak{f}(z),(u_M^{\pi})^*(w_2^{\pi})) \\
							&=(w_1^{\pi}\smile \Sq^1(-))^*(z)+ (u_M^{\pi})^*(w_2^{\pi})\frown z\\
							&=(w_1^{\pi}\smile \Sq^1(-))^*(z)+ ((u_M^{\pi})^*(w_2^{\pi})\smile-)^{*} (z),
						\end{align*}
						where the first equality is by the definition of $\delta_2^{M,w_M}$, the second is by the above and the definition of $\delta_2^{M,\PD(w_2(M))}$, the third is by \Cref{prop:spinCasedelta2d2}, the fourth is by \Cref{lem:CapLambdaCheck}, and the last is by the cup-cap formula.
					\end{proof}
					
					Everything we have done so far implies \Cref{mnthm:2}.  For the reader's convenience, we spell this out now.
					
					\begin{proof}[Proof of \Cref{mnthm:2}]
						Let $x\in H_4(\pi;\Z^{w_1^{\pi}})$. By an analogous argument to that in \Cref{subsec:primary} there exists a $(B\pi,w_1^{\pi})$-manifold $M$ such that $u_M^{\pi}[M]=x$. Then immediately \Cref{lem:d2sec} implies that $\delta_2^{4,0}(x) =d_2^{4,0}(x)$.  
						
						Let $z\in H_3(\pi;\Z/2)$.  Pick $(Y,u_Y)$ a $\xi$-3-manifold such that $\secondary(Y,u_Y)=0$.  By \Cref{lem:sec_vanishing_better} there exists a $(B\pi,w_1^{\pi})$-filling $(X,u_X^{\pi})$ with obstruction class $\smallO(X,Y,u_Y,u_X^{\pi})$ Poincar\'{e}-Lefschetz dual to some $w\in H_2(X;\Z/2)$ with a spherical representative.  By \Cref{prop:delta2isd2} we have that $\delta_2^{X,w}(z) = d_2^{(3,1)}(z)$.  This implies that we can define $\delta_2^{(3,1)}:= \delta_2^{X,w}$ since it is independent of $X$ and $w$; tautologically this implies that $\delta_2^{(3,1)}=d_2^{(3,1)}$. 
					\end{proof}
					
					Similarly, this implies \Cref{mnthm:wu}.
					
					\begin{proof}[Proof of \Cref{mnthm:wu}]
						Let $\pi:=\pi_1(X)$ and let $\xi=\xi(\pi,0,0)$. The spin structure $s$ gives us a $\xi$-structure on $Y$. From \Cref{prop:relative_obstruction_class} one can compute that $\smallO(X,Y,u_Y,u_X^{\pi}) = w_2^{s}$, the relative Stiefel-Whitney class. Let $b$ be a $\Z/2$-null-homologous immersed sphere in $X$.  By the Serre exact sequence we have that there exists $z\in H_3(\pi;\Z/2)$ such that $b$ is a lift of $\mathfrak{f}(z)$.  For the normal 1-type $\xi$ we have that
						\[
						\lambdadash(b,c) - \mudash(b) = \delta_2^{X,w}(z) = \delta_2^{(3,1)}(z) = d_2^{(3,1)}(z) = 0
						\]
						where all equalities are modulo $IJ$.  The first follows from \Cref{prop:delata_2properties}; the second and the third follow from \Cref{mnthm:2}; and the fourth follows from \Cref{prop:JSSdifferential}, since $w_2^{\pi}$ and $w_1^{\pi}$ are both zero.  The result is now immediate.
					\end{proof}

					\subsection{The geometric differential \texorpdfstring{$\delta_3^{(4,0)}$}{delta3,4,0}}\label{sbs:delta_3_differential}
					
					In this subsection fix $\xi=\xi(B\pi,w_1^{\pi}, w_2^{\pi})$ a normal $1$-type.  We will now define the third geometric differential 
					\[
					\delta_3^{(4,0)}\colon\ker(\delta^{(4,0)}_2\colon H_4(\pi;\Z^{w_1^{\pi}})\to H_2(\pi;\Z/2))\to H_1(\pi;\Z/2)/\im(\delta_2).
					\]
					
					\begin{definition}\label{def:geometricDelta3}
						Represent $x\in \ker (\delta_2^{(4,0)}\colon H_4(\pi;\Z^{w_1^{\pi}})\to H_2(\pi;\Z/2))$ by a closed $(B\pi,w_1^{\pi})$-manifold $(M,u_M^{\pi})$, with a $c_M$ a spherical representative for $\smallO(M,u_M^{\pi})$ with respect to $\xi$.  Then we define 
						\[
						\delta_3^{(4,0)}(x) := \terG(M,u_M^{\pi},c_M) = [\mudashj(c_M)] \in H_1(\pi;\Z/2)/\im (\delta_2^{(3,1)}).
						\]
					\end{definition}
					
					Note that the above definition is  analogous to the definition of $\delta_2^{(4,0)}$ (see \Cref{def:sec_delta_2}) in that both are defined as the secondary/tertiary invariant of some closed $4$-manifold. 
					
					Similarly to $\delta_2^{(4,0)}$ we have to prove that this is well-defined, i.e.\ does not depend on the choice of a manifold representing $x$.  For $\delta_2^{(4,0)}$ we did this by showing that it was equal to the JSS $d_2$.  We will not do that here.
					
					\begin{proposition}\label{prop:geometricDelta3}
						The above definition of $\delta_3^{(4,0)}$ is well-defined, i.e.\ it defines a map
						\[
						\delta_3^{(4,0)}:\ker(\delta^{(4,0)}_2\colon H_4(\pi;\Z^{w_1^{\pi}})\to H_2(\pi;\Z/2))\to H_1(\pi;\Z/2)/\im(\delta_2^{(3,1)}).
						\]
					\end{proposition} 	
					\begin{proof}
						First we show that $\delta_3$ does not depend on the choice of $c_M$. Let $c_M,c'_M\in\pi_2(X)$ be elements with $h_2(c_M)=h_2(c'_M)=\PD(\smallO(M,u_M))$ and set $a=c_M-c_M'\in \ker h_2$. Now we obtain
						\begin{align*}
							\mudashj(c_M)-\mudashj(c'_M)&=\mudashj(c'_M+a)-\mudashj(c'_M)\\
							&=\mudashj(c'_M)+\mudashj(a) +\lambdadash_J(a,c_M)-\mudashj(c'_M)\\
							&=\mudashj(a) + \lambdadash(c'_M,a)=\delta^{(3,1)}_2(a),
						\end{align*}
						where the second equality follows from $\lambdadash_J$ being the $\Z$-bilinear form corresponding to $\mudashJ$ (\Cref{prop:mudashj_lamdadashj_well-defined}) and the rest are immediate.
						
						Secondly we need to show that $\mudashj(c_M)\in H_1(\pi;\Z/2)/\im(\delta_2)$ only depends on $u_{M*}^{\pi}[M]$. Two closed $4$-dimensional manifolds $M_1,M_2$ with identified fundamental groups $\pi$ are $\CP^2$-stable diffeomorphic\footnote{i.e.\ diffeomorphic after connected-summing with copies of $\CP^2$ and $\ol{\CP^2}$ on both sides.} if and only if their fundamental classes agree in $H_4(\pi;\Z^{w_1^{\pi}})$. Thus, it suffices to show that $\mudashj(c_M)$ only depends on the $\CP^{2}$-stable class of $M$, i.e.\ that if $N=\CP^{2}$ or $N=\ol{\CP^{2}}$ and we take any class $c_{M\# N}\in \pi_2(M\# N)$ such that $h_2(c_{M\# N})=\PD(\smallO(M\# N,u_{M\# N}^{\pi}))$ we get $\mudashj(c_{M\# N})=\mudashj(c)$ in $H_1(\pi;\Z/2)/\im(\delta_2)$.  Recall that by \cref{prop:delata_2properties} the choice of $c_{M\# N}$ does not matter, as long as it represents the characteristic element.  We only carry out the proof for $N=\CP^2$; the other case is analogous.  Choose a representative for $c_{M\# \CP^{2}}$ to be a connected sum of a class $c\in \pi_2(M)$ and the class $[\CP^{1}]\in \pi_2(M\# \CP^{2})$. Now
						\[
						\mudashj(c+[\CP^1])=\mudashj(c)+\mudash([\CP^1])+\lambdadash_J(c,[\CP^{1}])=\mudashj(c)\in H_1(\pi;\Z/2)/\im(d_2),
						\]
						where the first equality follows from \Cref{prop:mudashj_lamdadashj_well-defined}, and the second since $\CP^1$ is embedded and disjoint from $c$.
					\end{proof}
					
					\subsection{Well-definedness of the geometric tertiary invariant}
					
					We now take the definition of the geometric tertiary obstruction given in \Cref{sec:obstructions} and show that we can make it only depend on the bounding $\xi$-3-manifold once we quotient out by the differentials $\delta_2^{(3,1)}$ and $\delta_3^{(4,0)}$.  This is encapsulated by the following theorem.
					
					\begin{theorem}\label{thm:terInvariants}
						Let $(Y,u_Y)$ be a $\xi$-3-manifold with two $(B\pi,w_1^{\pi})$-fillings $(X,u^{\pi}_X)$ and $(X',u^{\pi}_{X'})$ such that the obstruction $\smallO(X,Y,u_Y,u^{\pi}_X)$ has two spherical representatives $c_1$ and $c_2$, and such that $\smallO(X',Y,u_Y,u^{\pi}_{X'})$ has some spherical representative $c'$.  Then we have the following.
						\begin{enumerate}
							\item\label{it:ter1}
							$\terG(X,Y,u_Y,u_X^{\pi},c_1)-\terG(X,Y,u_Y,u_X^{\pi},c_2)\in \im(\delta_2),$
							\item\label{it:ter2}$\terG(X,Y,u_Y,u_X^{\pi},c_1)-\terG(X',Y,u_Y,u_{X'}^{\pi},c')\in \im(\delta_2,\delta_3).$
						\end{enumerate}
					\end{theorem}
					
					A direct consequence of \Cref{thm:terInvariants} is that the following are well-defined.
					
					\begin{enumerate}[(a)]
						\item$\terG(X,Y,u_Y,u_X^{\pi}):=[\mudashj(c)]\in H_1(\pi;\Z/2)/\im(\delta_2)$ for some choice of $c$ representing $\smallO(X,Y,u_Y,u^{\pi}_X)$.
						\item$\terG(Y,u_Y):=[\mudashj(c)]\in H_1(\pi;\Z/2)/(\im(\delta_2),\im(\delta_3))$ for some choice of $(B\pi,w_1^{\pi})$-filling $(X,u_X^{\pi})$ and $c$ representing $\smallO(X,Y,u_Y,u^{\pi}_X)$.
					\end{enumerate}
					
					\begin{proof}[Proof of \Cref{thm:terInvariants}]
						Note that (\ref{it:ter1}) was already proved in the proof of \Cref{prop:geometricDelta3}.
						
						We now prove (\ref{it:ter2}). Let $w\in H_2(X;\Z/2)$ and $w'\in H_2(X';\Z/2)$ be the Poincar\'{e}-Leftschetz duals of $\smallO(X,Y,u_Y,u_X^{\pi})$ and $\smallO(X',Y,u_Y,u_{X'}^{\pi})$, respectively.  Let 
						\[
						(M,u_M^{\pi}):= (X\cup_Y \ol{X'}, u_X^{\pi}\cup_{u_Y}u_{X'}^{\pi})
						\] 
						be the closed $(B\pi,w_1^{\pi})$-manifold obtained by gluing together $(X,u_X^{\pi})$ and $(X',u_{X'}^{\pi})$ together along $(Y,u_Y)$.  By \Cref{lem:doubleObstructionClass} the obstruction to $X\cup X_2$ having a $\xi$-structure is the sum $w+w'\in H_2(X_1\cup \ol{X_2};\Z/2)$, represented by $c:=c_1\# c'$. Surger the map $u_M^{\pi}$ away from $c$ to be an isomorphism on fundamental groups, such that $c$ still represents the obstruction to $M$ admitting a $\xi$-structure. Using \Cref{prop:geometricDelta3} we see that $\mudashj(w+w')\in \im(\delta_3)$.
					\end{proof}

					\section{Vanishing of the geometric tertiary obstruction}\label{sec:vanishing}
					
					This section is first concerned with proving \Cref{mnthm:1}.  The main result in that direction is the following.
					
					\begin{theorem}\label{thm:ter=0}
						Let $\xi=\xi(\pi,w_1^{\pi},w_2^{\pi})$ be a normal 1-type and let $(Y,u_Y)$ be a $\xi$-3-manifold with $\pri(Y,u_Y)=0=\secondary(Y,u_Y)$.  Then $\terG(Y,u_Y)=0$ if and only if there exists a $\xi$-filling of $(Y,u_Y)$.
					\end{theorem}
					
					Recall that we have already seen that the vanishing of primary obstruction is equivalent to the existence of a $(B\pi,w_1^{\pi})$-filling (\Cref{prop:PrimaryObstruction}).  We have also seen that the vanishing of the secondary obstruction is equivalent to the existence of a $(B\pi,w_1^{\pi})$-filling $X$ such that the obstruction class to extending the $\xi$-structure is spherical (\Cref{lem:sec_vanishing_better}).  We will take the consequences of these two lemmas to be our starting point for proving \Cref{thm:ter=0}, i.e.\ we will refer to the following setup.
					
					\begin{setup}\label{stp:vanishing_of_ter}
						We have a $\xi$-3-manifold $(Y,u_Y)$ and a $(B\pi,w_1^{\pi})$-filling $(X,u^{\pi}_X)$ such that the Poincar\'{e} dual of the obstruction $\smallO(X,Y,u^{\pi}_X,u_Y)\in H^2(X,Y;\Z/2)$ has a spherical representative $c\in \pi_2(X)$. 
					\end{setup}
					
					We now begin a sequence of lemmas which concern the various consequences if $\mudashj(c)$ `vanishes' to different extents. 
					
					\begin{lemma}\label{lem:mu_vanish_group_ring}
						Assume we are in the situation described by \Cref{stp:vanishing_of_ter}.  If $\mudash(c)=0\in I\subset \Z[\pi]$ then either 
						\begin{enumerate}[(a)]
							\item The class $c$ can be represented by an embedding in $X\#k (S^2\times S^2)$ for some $k$.
							\item There exists an immersion $c'\subset \CP^2\#\CP^2$ such that $[c']\in H_2(\CP^2\# \CP^2;\Z/2)$ is characteristic and such that $c\# c'$ is represented by an embedding in $X\# 2\CP^2\#k (S^2\times S^2)$ for some $k$.
						\end{enumerate}
						In either case, we can find a new null-bordism of $(Y,u_Y)$ such that the obstruction class to it being a null-$\xi$-bordism can be embedded.
					\end{lemma}
					
					\begin{proof}
						Here we make use of \Cref{sec:km} and the Kervaire-Milnor invariant.  The rough idea of the argument is as follows.  First, argue that we can arrange that $\km(c)=0$, and then apply \Cref{lem:km_stable_embedding} to make $c$ embedded.  While doing this we check that $c$ remains a representative for our obstruction class throughout our modifications.
						
						We begin by assuming that $\km(c)=1$, since otherwise we may skip the step of arranging that $\km(c)=0$.  Now we freely arrange that $c$ has an algebraic dual sphere by connected-summing with $\CP^2$ and tubing $c$ into $\CP^1\subset \CP^2$ (note that since $\CP^2$ is simply-connected and oriented, there is no issue in extending the $(B\pi,w_1^{\pi})$-structure).  Since the $\CP^1\subset \CP^2$ is a representative of $\smallO(\CP^2,\CP^2\xrightarrow{*}B\pi)$, it is clear that $c\# \CP^1$ is a representative for $\smallO(X\# \CP^2,Y,u_{X\# \CP^2}^{\pi},u_Y)$.  We subsume the $\CP^2$ into $X$ and assume that our $c$ has an algebraic dual sphere.  
						
						Again form the connected-sum $X\# \CP^2$ and let $b$ be the immersed sphere in $\CP^2$ given by \Cref{lem:km_example} which is characteristic, has $\km(b)=1$ and $\mudash(b)=0$.  Tube $c$ into $b$ to form $c\# b$, noting that $c\# b$ is again a representative for $\smallO(X\# \CP^2,Y,u_{X\# \CP^2}^{\pi},u_Y)$, since $b$ is characteristic.  We split into two cases.  First assume that $c\# b$ is r-characteristic.  Then, by \Cref{thm:kervaire_milnor} and that $c$ has an algebraic dual sphere, we can calculate that $\km(c\# b)=1+1=0$, since a collection of Whitney discs for the double-points of $c\# b$ is given by a union of the Whitney discs for $b$ and $c$ separately, and these Whitney discs are contained inside their respective connected-summands.  In the second case we assume that $c\# b$ is not r-characteristic.  Then by \Cref{thm:kervaire_milnor} and that $c$ has an algebraic dual sphere, $\km(c\#b)=0$.  This completes the argument that $\km(c\# b)=0$.  Again we subsume the $\CP^2$ into $X$ and assume that our $c$ has $\km(c)=0$.  
						
						We now use \Cref{lem:km_stable_embedding} and stably embed $c$. One should note that, since $S^2\times S^2$ is a simply-connected spin manifold, there is no issue in extending the $\xi$-structure when taking connected-sums with $S^2\times S^2$, and hence $c$ is still a representative for $\smallO(X\# (S^2\times S^2), Y, u_Y, u^{\pi}_{X\# (S^2\times S^2)}$).
					\end{proof}
					
					The following three lemmas are a consequence of the proof of \Cref{thm:terInvariants}
					
					\begin{lemma}\label{lem:mu_vanish_homology}
						Assume we are in the situation described by \Cref{stp:vanishing_of_ter}.  If $\mudashj(c)=0\in I/IJ$ then there exists an additional closed $(B\pi,w_1^{\pi})$-manifold $X'$ and an immersion $c'\subset X'$ such that $c\# c'\subset X\# X'$ represents the obstruction to extending the $\xi$-structure across $X\# X'$ and $\mudash(c\# c')=0\in I\subset \Z[\pi]$.
					\end{lemma}
					
					\begin{proof}
						
						By \cref{lem:KernelOfTheMapItoHomology} the kernel of the map $f\colon I\to I/IJ$ is generated by $\Lambda$-linear combinations of $\{2g\mid g\in \pi\}$ and $\{gh-g-h+1\mid g,h\in \pi\}$. By assumption, we have that $\mudash(c)\in \ker(f)$. For each of these generators $\omega\in \ker(f)$ we will exhibit a closed $(B\pi,w_1^{\pi})$-manifold $(M,u_M^{\pi})$ and an element $c_M\in\pi_2(M)$ such that under the Hurewicz map $h_2(c_M)=\PD(\smallO(M,u_M^{\pi}))$. 
                        
						Pick an element $g\in \pi$ other than the identity. Consider $M=S^1\times S^3\#\CP^2$ in the case that $w_1^{\pi}(g)$ is trivial and  $M=S^1\wt{\times} S^3\#\CP^2$ in the case that $w_1^{\pi}(g)$ is nontrivial.  Denote the generator of the fundamental group of $M$ by $t$.  There is a map $M\to B\pi$, which sends $t$ to $g$ and we use this to endow $M$ with a $(B\pi,w_1^{\pi})$-structure.  Let $p=(2-t)[\CP^1]\in\pi_2(M)$. This element has the property that $h_2(p)=\PD(w_2(M))=\PD(\smallO(M,u_M^{\pi})$.  By \Cref{lem:mudash_lambdadash_sum_formula} we have
						\[\begin{cases}
							w_1^{\pi}(t) \text{ trivial}&	\lambdadash(p,p)=(2-t)\ol{(2-t)}-1=4-2t-2t^{-1}=\mudash(p)+\ol{\mudash(p)}\\
							w_1^{\pi}(t) \text{ non-trivial}&	\lambdadash(p,p)=(2-t)\ol{(2-t)}-3=-2t+2t^{-1}=\mudash(p)+\ol{\mudash(p)}
						\end{cases}
                        \]
						and from this we conclude that $\mudash((2-t)[\CP^1])=(2-2t)$. From the connected-sum $X\#M$ form $X'$ by representing the loop $gt^{-1}$ away from the representatives of $c$ and $p$ and surgering that loop in a $\xi$-preserving way. Note that in particular the loop $gt^{-1}$ is framable. It can be seen that the new obstruction class $\smallO(X',Y,u_{X'}^{\pi}, u_Y)$ is represented by the immersion $c\#p$. Since the surgery was performed away from $c$ or $p$ we obtain that $\mudash(w\# p)=\mudash(w)+2-2g$.

						Pick $g,h\in \pi$ two nonzero elements.	Consider 
						\[
						M=\begin{cases}
							S^1\!\times S^3\#S^1\!\times S^3\#\CP^2\quad \text{ if } w_1^{\pi}(g)=0,w_1^{\pi}(h)=0\\
							S^1\wt{\times} S^3\#S^1\! \times S^3\#\CP^2\quad \text{ if } w_1^{\pi}(g)=1,w_1^{\pi}(h)=0\\
							S^1\! \times S^3\#S^1\wt{\times} S^3\#\CP^2\quad \text{ if } w_1^{\pi}(g)=0,w_1^{\pi}(h)=1\\
							S^1\wt{\times} S^3\#S^1\wt{\times} S^3 \#\CP^2\quad  \text{ if } w_1^{\pi}(g)=1,w_1^{\pi}(h)=1\\
						\end{cases}
						\]
						and fix the two generators $t,s$ of $\pi_1(M)$. There is a map $M\to B\pi$ which sends $t$ to $g$ and $s$ to $h$, and using this map endows $M$ with a $\xi(\pi,w_1^{\pi},w_2^{\pi})$ structure. Pick $p=(t+s^{-1}-1)[\CP^{1}]$.  By \Cref{lem:mudash_lambdadash_sum_formula} and for example for $w^{\pi}_1(t)=0=w_1^{\pi}(s)$ we see that
						\[
						\lambdadash(p,p)=(t+s^{-1}-1)\ol{(t+s^{-1}-1)}=ts+st-s-s^{-1}-t-t^{-1}+3,
						\]
						and so $\mudash(p)=ts-s-t$. It can be verified that for all other combinations of $w_1^{\pi}(t)$ and $w_1^{\pi}(t)$ the result is the same. By the same argument as before, form $X\# M$ and perform two $\xi$-preserving surgeries on $gt^{-1}$ and $hs^{-1}$ and observe $\mudash(w\#p)=\mudash(w)+gh-g-h$.
						
						Now performing some finite combination of the above operations completes the proof.
					\end{proof}
					\begin{lemma}\label{lem:mu_vanish_delta_2}
						Assume we are in the situation described by \Cref{stp:vanishing_of_ter}.  If $\mudashj(c)\in\im \delta_2^{(3,1)}$ then there exists a $\Z/2$-null-homologous $\alpha\in\pi_2(X)$ such that $\mudashj(c\# \alpha)=0\in I/IJ\cong H_1(\pi;\Z/2)$.
					\end{lemma}
					\begin{proof}
						Let $z\in H_3(X;\Z/2)$ be such that $\mudash_J(c)=\delta_2^{(3,1)}(z)$ and choose any lift $\alpha\in\pi_2(X)$ of $\mathfrak{f}^{X}(z)$. Then since $\alpha$ is $\Z/2$-nullhomologous we can calculate like in the proof of \Cref{prop:geometricDelta3} that						\[
						\mudashj(c\#\alpha)-\mudashj(c)=\delta^{(3,1)}_2(z),
						\]
						showing that $\mudashj(c\#\alpha)=0$ in $I/IJ$.
					\end{proof}
					
					\begin{lemma}\label{lem:mu_vanish_delta_3}
						Assume we are in the situation described by \Cref{stp:vanishing_of_ter}.  If $\mudashj(c)\in\im \delta_3$ then there exists a different manifold $X'$ with $c'\in\pi_2(X')$ like in the \Cref{stp:vanishing_of_ter} such that $\mudashj(c')=0\in H_1(\pi;\Z/2)/\im \delta_2$.
					\end{lemma}
					\begin{proof}
						Since $\mudashj(c)\in\im \delta^{(4,0)}_3$ there exists an element $a\in H_4(\pi;Z^{w_1^{\pi}})$ which hits it.  Therefore there exists a closed $(B\pi,w_1^{\pi})$-manifold $M^4$ with $\smallO(M,u_M^{\pi})$ represented by $c_M\pi_2(M)$ such that $\mudashj(c_M)=\mudashj(c)$.  Let $(X\# M)'$ denote the result of surgering the map $u_X^{\pi}\# u_M^{\pi}$ to a $\pi_1$-equivalence, where the surgeries are performed away from some fixed representatives of $c$ and $c_M$ in a $\xi$-preserving way.  Then 
						\[
						[\mudashj(c\# c_M)]=[\mudashj(c)+\mudashj(c_M)]=0\in H_1(\pi;\Z/2)/(\im\delta_2^{(3,1)})
						\] 
						where the first equality comes from the fact that $c$ and $c_M$ are disjoint in $(X\# M)'$ along with \Cref{prop:mudashj_lamdadashj_well-defined}, and the second equality is by our assumption.
					\end{proof}
					
					\begin{proof}[Proof of \Cref{thm:ter=0}]
						
						We start with the `if' direction, i.e.\ we assume that there exists a $\xi$-filling $(X,u_X)$ of $(Y,u_Y)$.  We can choose a representative for $\smallO(X,Y,u_Y,u^{\pi}_X)$ to be a local unknot which bounds an embedded $D^3$ and hence $\terG(Y,u_Y)=0$.
						
						We now consider the `only if' direction.  By \Cref{prop:PrimaryObstruction} and \Cref{lem:sec_vanishing_better}  we can assume that we are in situation described in \Cref{stp:vanishing_of_ter}.  We now prove that we can find a new $\xi$-filling for $(Y,u_Y)$ such that the $c$ is embedded.  
						
						By \Cref{lem:mu_vanish_delta_3} and \Cref{lem:mu_vanish_delta_2} we can assume that our filling $X$ has a representative $c$ for $\smallO(X,Y,u_Y,u_X^{\pi})$ with $\mudashj(c)=0\in H_1(\pi;\Z/2)$.  By \Cref{lem:mu_vanish_homology} we can further assume that $c$ has the property that $\mudash(c)=0\in\Z[\pi]$ and by \Cref{lem:mu_vanish_group_ring} this implies that we can assume that $c$ is embedded. All the steps above allow us to change our filling manifold $X$.
						
						Now, generalising an argument due to Kervaire-Milnor \cite{kervaire_milnor}, we can arrange that $c$ has Euler number $+1$ by connected-summing with $\CP^2$ and $\ol{\CP}^2$ and tubing into the $\CP^1\subset \CP^2$ or $\CP^1\subset \ol{\CP}^2$ if necessary.  Blowing down $c$, i.e.\ replacing the tubular neighbourhood of $c$ by a copy of $D^4$, reduces the problem to extending the given $\xi$-structure on $S^3$ to $D^4$, which we now briefly demonstrate can always be achieved.
						
						Consider the following diagram.
						\[
						\begin{tikzcd}
							& & B \arrow[d,"{\xi(\pi,w_1^{\pi},w_2^{\pi})}"] \arrow[r,"p"] & B\pi \arrow[d,"w_1^{\pi}\times w_2^{\pi}"] \\
							S^3 \arrow[urr,"u_{S^3}"] \arrow[r,hook] & D^4 \arrow[r,"\nu_{D^4}"'] & BO \arrow[r,"w_1\times w_2"'] & K(\Z/2,1)\times K(\Z/2,2)
						\end{tikzcd}
						\]
						By the definition of the normal 1-type $B$, we have that $u_{S^3}$ extends to $D^4$ provided that we can extend the map $p\circ u_S^3$ such that the corresponding relative cohomology class in $H^1(D^4,S^3;\Z/2)\oplus H^2(D^4,S^3;\Z/2)$ given by postcomposing with $w_1^{\pi}\times w_2^{\pi}$ is trivial, but this is clear since the cohomology groups are themselves trivial.  This means that we have shown that $(Y,u_Y)$ is $\xi$-nullbordant, which finishes the proof.
					\end{proof}
					
					As a consequence of \Cref{thm:ter=0} we have the following result relating the tertiary and geometric tertiary invariants.
					
					\begin{corollary}\label{cor:tersVanishTogether}
						Let $(Y,u_Y)$ be a $\xi$-3-manifold, with $\pri(Y,u_Y)=0$ and $\secondary(Y,u_Y)=0$. Then the geometric tertiary invariant $\terG(Y,u_Y)$ vanishes if and only if the usual tertiary invariant $\ter(Y,u_Y)$ vanishes.
					\end{corollary}
					
					\begin{proof}
						Assume that $\terG(Y,u_Y)$ vanishes.  Then by \Cref{thm:ter=0} there exists a $\xi$-filling and hence $\ter(Y,u_Y)$ must also vanish.  Similarly if $\ter(Y,u_Y)$ vanishes then the standard AHSS/JSS argument proves that there exists a $\xi$-filling and hence by \Cref{thm:ter=0} we must have that $\terG(Y,u_Y)$ vanishes.
					\end{proof}
					
					We can now prove \Cref{mnthm:1} from the introduction.
					\begin{proof}[Proof of \Cref{mnthm:1}]
						 We have that $u_Y^{\pi}([Y])=0$ if and only if $(Y,u_Y)\in F_{2,1}$ by \Cref{lem:primary}.  Furthermore, $(Y,u_Y)\in F_{2,1}$ if and only if there exists a $(B\pi,w_1^{\pi})$-filling of $(Y,u_Y)$ by \Cref{prop:PrimaryObstruction}.  We now assume this filling $(X,u_X^{\pi})$ exists. We have that $(u_{X}^{\pi})_*\PD(\smallO(X,Y,u_Y,u_X^{\pi}))=0$ if and only if $(Y,u_Y)\in F_{1,2}$ by \Cref{prop:SecXYuY}.  Furthermore $(Y,u_Y)\in F_{1,2}$ if and only if there exits a $(B\pi,w_1^{\pi})$-filling such that $\smallO(X,Y,u_Y,u_X^{\pi})$ has a spherical representative by \Cref{lem:sec_vanishing_better}.  We now assume such a filling $(X',u_{X'}^{\pi})$ exists with $c$ a spherical representative for $\smallO(X,Y,u_Y,u_{X'}^{\pi})$.	We have that $\terG(Y,u_Y)= [\mudashj(c)]=0$ if and only if there exists a $\xi$-filling of $(Y,u_Y)$ by \Cref{thm:ter=0}.  Note that this is equivalent to $(Y,u_Y)\in F_{0,3}$ since $F_{0,3}=\{0\}\in \Omega_3^{\xi}$.
					\end{proof}
					
					\subsection{Comparing \texorpdfstring{$\terG$}{terGeo} and \texorpdfstring{$\ter$}{ter}.}
					
					We now show that $\ter$ and $\terG$ are `equal'.  We will make what this means precise later, but for now we need a realisation result for $\terG$.
					
					\begin{proposition}\label{prop:realisation_for_terG}
						Let $\xi=\xi(\pi,w_1^{\pi},w_2^{\pi})$ be a normal $1$-type and let $g\in \pi$.  Then there is a compact $4$-dimensional manifold $V_g$ and a map $u_{V_g}^{\pi}\colon V_g\to B\pi$ such that
						\begin{enumerate}
							
							\item\label{it:1Properties} $\partial V_g=S^1\times S^1\times S^1$ if $w_1^{\pi}(g)=0$, $\partial V_g=S^1\times K$, for $K$ a Klein bottle if $w_1^{\pi}(g)\neq 0$.
							\item\label{it:2Properties} There exists an identification $\pi_1(V_g)\cong \Z$ such that under this identification $u_{V_g}^{\pi}$ sends $+1$ to $g$.
							\item\label{it:3Properties} $u_{V_g}^{\pi*}(w_1^{\pi})=w_1(V_g)$.
							\item\label{it:4Properties} The (normal) second Stiefel-Whitney class $w_2(V_g)$ vanishes and so the normal $1$-type of $V_g$ is $\xi(\Z,w_1(V_g),0)$.
							\item\label{it:5Properties}  There exists a $\xi(\Z,w_1(V_g),0)$-structure $p$ on $\partial V_g$, inducing the map $\pi_1(\partial V_g)\to \pi_1(V_g)\cong \Z$ such that the obstruction $\smallO(V_g,\partial V_g,p,u_{V_g}^{\pi})$ is the Poincar\'{e} dual of an immersed sphere $S$ in $V_g$ such that $\mudashJ(S)=t-1\in I/IJ$.
						\end{enumerate}
                        
					\end{proposition}
					\begin{proof}
						\begin{figure}[h!]
							\centering
							\includegraphics[width=0.5\linewidth]{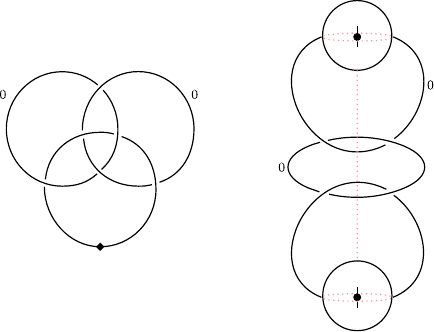}
							\caption{On the left: a Kirby diagram for the manifold $V_g$ if $g$ is orientable.  On the right: a Kirby diagram for the manifold $V_g$ if $g$ is non-orientable.  The notation for the 1-handles is due to Akbulut \cite[Sections 1.1, 1.5]{akbulut}}
							\label{fig:ter}
						\end{figure}
						
						Let
						\[
						Y_g=\begin{cases}
							S^1\times S^1\times S^1&w_1(g)=0\\
							S^1\times K&w_1(g)\neq 0,
						\end{cases}
						\]
						where $K:=S^1\widetilde{\times} S^1$ denotes the Klein bottle. In the non-orientable case the orientation character is given by the composition $S^1\times(S^1\wt{\times} S^1)\xrightarrow{q}S^1\to B\Z/2$ (here $q$ is the projection onto $K$ postcomposed with the fibre bundle map).  In both cases above $V_g$ will be defined as a suitable filling of $Y$.
						
						If $w_1(g)=0$, define $u_{Y_g}^{\pi}\colon Y_g=S^1\times S^1\times S^1\xrightarrow{\pr_3} S^1=B\Z$, where $\pr_3$ denotes projection onto the third coordinate. We describe a $(B\Z,0)$-filling of $Y_g$ by taking the filling $S^1\times D^2\times S^1$ and surgering out the kernel of $\pi_1(u_{Y_g}^{\pi})$ to form the manifold $V_g$, which is described using a Kirby diagram in \Cref{fig:ter}. Fix the induced $\pi_1$-isomorphism $V_g\to B\Z$ which restricts to $u_{Y_g}^{\pi}$ on the boundary.
						
						If $w_1(g)\neq 0$, define $u_{Y_g}^{\pi}\colon Y_g=S^1\times S^1\wt{\times} S^1\xrightarrow{q} B\Z=S^1$. Note that the first Stiefel-Whitney class of $Y_g$ factors through $u_{Y_g}^{\pi}$. We describe a $(B\Z,w_1^{\pi})$-filling of $Y_g$ by taking the filling $S^1\times H$ where $H$ denotes the solid Klein bottle, and then surgering out the kernel of $\pi_1(u_{Y_g}^{\pi})$ to form $V_g$, which is described using a Kirby diagram in \Cref{fig:ter}. Fix the induced $\pi_1$-isomorphism $V_g\to B\Z$ which restricts to $u_{Y_g}^{\pi}$ on the boundary.
						
						Now (\ref{it:1Properties}), (\ref{it:2Properties}) and (\ref{it:3Properties}) are immediate. The integral intersection form on the second homology of manifolds in \Cref{fig:ter} is even in both cases, so the second \emph{tangential} Wu class, equal to the second \emph{normal} Stiefel-Whitney  vanishes, establishing (\ref{it:4Properties}) (this follows since there is no $2$-torsion in the fundamental group of $V_g$). 
						
						Pick $p_B$ to be some $\xi(\Z,w_1^{\pi},0)$-structure on $V_g$. Recall from \cref{rmk:SpinPinxi} that a $\xi(\Z,w_1^{\pi},0)$ structure on a manifold $Y_g$ is a map $u_{Y_g}^{\pi}\colon Y_g\to B\Z$ over $w_1^{\pi}$ and a $\Pin^{+}$-structure on $Y_g$. Since we already have the former, we now define a suitable $\Pin^+$-structure.  The handle decomposition given in \Cref{fig:ter} defines an isomorphism 
						\[
						H^2(V_g,\partial V_g;\Z/2)\cong H_2(V_g;\Z/2)\to (\Z/2)\oplus (\Z/2).
						\] 
						Let $x\in H^{1}(Y_g;\Z/2)$ be any preimage of $(1,1)\in H^2(V_g,Y_g;\Z/2)$ (note that this exists since the map $H^{1}(Y_g;\Z/2)\to H^2(V_g,Y_g;\Z/2)$ is surjective). $\Pin^+$-structures on $Y_g$ are a torsor over $H^1(Y_g;\Z/2)$ and so we can define $p=p_B|_{Y_{g}}\cdot x$.
						
						The $\xi(\Z,w_1^{\pi},0)$-structure on $Y_g$ does not extend to $V_g$ because $x$ does not come from $H^1(V_g;\Z/2)$. By construction, the obstruction $\smallO(V_g,\partial V_g,p,V_g\to B\Z)$ to $p$ extending to $V_g$ is represented in either case by an immersed 2-sphere $S$ consisting of the core of both handles, connected by the obvious annulus in \Cref{fig:ter}.
						Both of these spheres have two self intersections, only one of them with non-trivial group element giving $\mudash_J(S)=t-1$.  This proves (\ref{it:5Properties}).
					\end{proof}
					
					To show that $\ter$ and $\terG$ agree, we must show they land in the same group. In \Cref{prop:delta2isd2} we have already shown that $d_2, \delta_2$ coincide, in particular we have to show that $\im (d_3)=\im( \delta_3)\subseteq H_1(\pi;\Z/2)/\im(d_2)$. We show the two inclusions separately.

                    \begin{lemma}\label{lem:delta3subsetd3}
    For $\xi=\xi(\pi,w_1^{\pi},w_2^{\pi})$ a normal $1$ -type we have for $d_3:E^3_{4,0}\to E^{3}_{1,2}$ the differential of the JSS and $\delta_3$ the geometric map defined in \cref{def:geometricDelta3} that \[\im(\delta_3)\subseteq \im (d_3)\subset H_1(\pi;\Z/2)/\im(d_2).\] 
\end{lemma}
\begin{proof}
 Let $x\in \im \delta_3 \subset H_1(\pi;\Z/2)/\im(d_2)$ be a nonzero element.  Then we can represent the equivalence class $[x] \in H_1(\pi;\Z/2)/\im(d_2,d_3)$ by a $\xi$-3-manifold $(Y,u_Y)$ that lies in the corresponding bordism class.  \Cref{cor:tersVanishTogether} implies that $(Y,u_Y)$ is $\xi$-null-bordant and hence $x\in \im d_3$.
\end{proof}

					 \begin{theorem}\label{thm:tersAgree}
						Fix $\xi=\xi(\pi,w_1^{\pi},w_2^{\pi})$ a normal $1$ -type. Let $(Y,u_Y)$ be a $\xi$-3-manifold, with $\pri(Y,u_Y)=0=\secondary(Y,u_Y)$. 
                        Then 
						\[
						[\terG(Y,u_Y)]=\ter(Y,u_Y)\in H_1(\pi;\Z/2)/\im(d_2,d_3).
						\]
					\end{theorem}

					\begin{proof}
                    If $\terG(Y,u_Y)=0$ we have already seen in \Cref{cor:tersVanishTogether} that $\ter(Y,u_Y)=0$. So we assume $\terG(Y,u_Y)$ is nonzero.
                    
                    Let $(X,u_X^{\pi})$ be  a $(B\pi,w_1^{\pi})$-filling of $(Y,u_Y)$, and further let $c$ be a spherical representative for~$\smallO(X,Y,u_Y,u_X^{\pi})$ (see \Cref{lem:sec_vanishing_better}). Let $g\in\pi$ be any lift of $\mudashj(c)\in H_1(\pi;\Z/2)$.  By assumption $\mudashj(c)$ is nontrivial and so $g$ is non-trivial. We have that $\mudashJ(c)=g-1\in I/IJ$ (see \Cref{lem:tensorIsos}(\ref{it:H1group3}) for the isomorphism $H_1(\pi;\Z/2)\cong I/IJ$). By abuse of notation, we choose a map $g\colon B\Z\to B\pi$ that sends the generator to $g$. Fix a homotopy $H\colon B\Z\times I\to B\Z/2\times B^2\Z/2$   between the following maps
						\[
						\begin{tikzcd}
							& B\pi \arrow[rd, "w_1^{\pi}\times w_2^{\pi}"] &\\
							B\Z \arrow[ru, "g"] \arrow[rr, "w_1^g\times *"] && B\Z/2\times B^2\Z/2
						\end{tikzcd}
						\]
						and use it to obtain a map $\varphi_H\colon B_{\xi(\Z,w_1^g,0)}\to B_{\xi}$ over $BO$
						
						\[
						\begin{tikzcd}
							B_{\xi(\Z,w_1^g,0)}\ar[rr]\ar[ddr]\ar[rd,dashed,"\varphi_H"]&&B\Z\ar[dr,"g"]\ar[ddr,"w_1^g\times *"',pos=0.1]&\\
							&B_{\xi}\ar[rr,crossing over]\ar[d]&&B\pi\ar[d,"w_1^{\pi}\times w_2^{\pi}"]\\
							&BO\ar[rr,"w_1\times w_2"]&&B\Z/2\times B^2\Z/2
						\end{tikzcd}
						\]
						which maps $\xi(\Z,w_1^{\pi},0)$-manifolds to $\xi$-manifolds. Furthermore, this gives a map of Thom spectra $M\xi(\Z,w_1^{\pi},0)\to M\xi$ which induces a map of bordism groups
						\[
                        (\varphi_{H})_*\colon\Omega_3^{\xi(\Z,w_1^g,0)}\rightarrow\Omega_3^{\xi}.
						\]
						  Let $V$ be the manifold given by \Cref{prop:realisation_for_terG} for the element $g$ and the orientation character $w_1^g$. We will use the notation from that proposition throughout this proof. 
						
						By (\ref{it:4Properties}) of \Cref{prop:realisation_for_terG} $V$ admits a 1-smoothing $p_B$ with respect to $\xi(\Z,w_1(V),0)$ which we now fix.  By (\ref{it:5Properties}) of \Cref{prop:realisation_for_terG} there is a $\xi(\Z,w_1(V),0)$-structure $p$ on  $\partial V$ whose obstruction to extending to $V$ is some immersed sphere $S$ with $\mudash(S)=t-1$, where $t$ is the generator of $\pi_1(V)$. There is an action of $H^1(\partial V;\Z/2)$ on the set of $\xi(\Z,w_1(V),0)$-structures on $\partial V$, given by the $H$-space structure on $B^2\Z/2$. Denote by $x\in H^1(\partial V;\Z/2)$ the difference of these structures, i.e.\ $x:=p-p_B\vert_{\partial V}$. It holds that under the connecting homomorphism 
						\[
						H^1(\partial V;\Z/2)\xrightarrow{\mathfrak{d}} H^2(V,\partial V;\Z/2)
						\]
						we have $\mathfrak{d}(x)=\PD^{-1}(S)$. 
						
						Define $[\partial V,u_B]=\varphi_{H*}([\partial V,p_B|_{\partial V}])$, which gives a $\xi$-structure on $\partial V$ extending to $V$. Define $[\partial V,u]:=[\partial V,u_B]\cdot x$. Equivalently we have $\varphi_H([\partial V,p])=[\partial V,u]$. From $\mathfrak{d}(x)=\PD^{-1}(S)$ we see that $S$ is indeed also the obstruction to extending $u$ to $V$. For a proof use an argument analogous to \cite[Page 521]{Kervaire57}.

						We now show that $\terG(Y,u_Y)=\ter(Y,u_Y)$. By assumption we have that $\terG((Y,u_Y)+(\partial V,u))=0$ and by \Cref{cor:tersVanishTogether} we have $\ter((Y,u_Y)+(\partial V,u))=0$. Now $\terG(Y,u_Y)=\ter(Y,u_Y)$ if and only if $\terG(\partial V,u)=\ter(\partial V,u)$.

						Using the James spectral sequence we find that $\ter^{\xi(\Z,w_1^g,0)}\in H_1(B\Z;\Z/2)\cong I/IJ\cong \Z/2$, generated by  $t-1$, is the only obstruction to an element being trivial in the bordism group $\Omega_3^{\xi(\Z,w_1^g,0)}$. But $\terG(\partial V,p)=t-1\in H_1(B\Z;\Z/2)$. By \cref{cor:tersVanishTogether} $\ter(\partial V,p)$ cannot vanish and so $\ter(\partial V,p)=t-1$. Now since $(\varphi_{H})_*$ comes from a map of Thom spectra we have the following naturality of the tertiary obstructions \[g_*(\ter^{\xi(\Z,w_1^g,0)}([\partial V,p]))=\ter^{\xi}(\varphi_{H}([\partial V,p]))=\ter^{\xi}(\partial V,u)\] and so $ \ter^{\xi}(\partial V,u)=g-1$. Since $V$ is a $(B\pi,w_1^g)$-filling of $\partial V$ we use \Cref{lem:surgeryBelowMiddle} to surger $V$ to make the map $V\to B\pi$ a $\pi_1$-isomorphism.  These surgeries can be made away from a representative of some spherical representative of $\smallO(V,\partial V,p,g)$ and so we can calculate $\terG(V,u)$ as follows \[\terG(V,u)=g_*(\terG(V,p))=g_*(t-1)=g-1.\]
                        
                        As explained above the result now follows.
					\end{proof}

                    \begin{proposition}\label{prop:d3delta3Agree}
						For $\xi=\xi(\pi,w_1^{\pi},w_2^{\pi})$ a normal $1$ -type we have for $d_3:E^3_{4,0}\to E^{3}_{1,2}$ the differential of the JSS and $\delta_3$ the geometric map defined in \cref{def:geometricDelta3} that \[\im(d_3)=\im (\delta_3)\subset H_1(\pi;\Z/2)/\im(d_2).\] 
					\end{proposition}
					
					\begin{proof}
						In \Cref{lem:delta3subsetd3} we have already shown that $\im (\delta_3)\subseteq \im(d_3)$.
                        
						Assume that $x\in \im d_3 \subseteq H_1(\pi;\Z/2)/\im (d_2)$. Represent $x$ by an element of $\pi$.  By \Cref{prop:realisation_for_terG} there exists a $(B\pi,w_1^{\pi})$-filling $(V,u_V^{\pi})$ of a $\xi$-3-manifold $(\partial V,u_{\partial V})$ with a spherical representative $c$ for $\smallO(V,\partial V,u_{\partial V}, u_V^{\pi})$ having $\mudashj(c)=x-1\in I/IJ$ and so $\mudashj(c)=x\in H_1(\pi;\Z/2)/\im(d_2)$. Now by \Cref{thm:tersAgree} and our assumption we have that \[0=[x]=[\terG(V,\partial V,u_{\partial V},u_V^{\pi})]=\ter(\partial V, u_{\partial V})\in H_1(\pi;\Z/2)/\im(d_2,d_3).\]  By the definition of $\ter$ this implies that there exists a $\xi$-filling $(X,u_X)$ of $(\partial V,u_{\partial V})$.  Now form a closed $(B\pi,w_1^{\pi})$-manifold $(M,u_M):= (V,u_V^{\pi})\cup_{(\partial V,u_{\partial V})}(X,u_X)$.  By construction $c$ is a representative for $\smallO(M,u_M^{\pi})$ and so $\delta_3((u_M^{\pi})_*[M])=\mudashj(c)=x$ and hence $x\in \im \delta_3$.
					\end{proof}

An immediate consequence of \Cref{thm:tersAgree} and \Cref{prop:d3delta3Agree} is

\begin{corollary}
    For $\xi=\xi(\pi,w_1^{\pi},w_2^{\pi})$ a normal 1-type we have that for any $\xi$-3-manifold $(Y,u_Y)$ such that $\pri(Y,u_Y)=0=\secondary(Y,u_Y)$ we have 

    \[
					\terG(Y,u_Y)=\ter(Y,u_Y)\in H_1(\pi;\Z/2)/\im(d_2,d_3)=H_1(\pi;\Z/2)/\im(\delta_2,\delta_3).
					\]
\end{corollary}

Finally we remark that the above does \emph{not} prove that $d_3$ and $\delta_3$ agree as maps, which we have not been able to determine.

					\appendix
					
					\section{The Serre exact sequence}\label{sec:SerreExSeq}
					
					Let $X$ be any connected space with a finite $4$-skeleton, let $\pi=\pi_1(X)$ and let $u_X^{\pi}\colon X\to B\pi$ be any classifying map. Then there is the following well-known exact sequence:
					\begin{equation}\label{eq:SerreExactSeq}
						H_3(X;R)\xrightarrow{(u_X^{\pi})_*} H_3(\pi;R)\xrightarrow{\mathfrak{f}^X:=d_3}\pi_2(X)\otimes_{\Lambda} R \xrightarrow{h^{\Lambda}_R} H_2(X;R)\xrightarrow{(u_X^{\pi})_*} H_2(\pi;R).
					\end{equation}
					with coefficients in any abelian group, though throughout this section we assume that $R$ is $\Z$ or $\Z/2$. This exact sequences arises as the five-term sequence of the Serre spectral sequence of the fibration $\widetilde{X}\to X\to B\pi$.  The map $\mathfrak{f}^X$ is then given by a certain third differential $d_3$. In the following section we prove that this map can be expressed in terms of the first $k$ invariant of $X$. We do not claim originality of this result, but we have been unable to find it in the literature.

					Let $k$ be the first $k$-invariant of $X$, i.e\ let $k\in H^3(\pi;\pi_2(X))$ be the element which classifies the principal fibration 
					\[B^2\pi_2(X)\to P_2(X)\to B\pi,\]
					where $P_2(X)$ is the second Postnikov stage of $X$. Let $\Lambda\to R$ be the augmentation map, making $R$ into a left $\Lambda$-module. Let  $k^r\in H^3(\pi;\pi_2(X)\otimes_{\Lambda}R)$ be the reduced $k$ invariant and $\ev(k^r)\in \Hom(H_3(\pi),R)$ the evaluation map. Then
					
					\begin{theorem}\label{thm:SerreExactSeq}
						Let $X$ be a space, $\pi=\pi_1(X)$ and $\wt{X}\to X\to B\pi$ be a fibration. Then for the map $d_3$ in the  exact sequence (\ref{eq:SerreExactSeq}) we have $d_3=\ev(k^r)$.    
					\end{theorem}
					\begin{remark}
						\Cref{thm:SerreExactSeq} is a special case of a more general theorem which will be stated in the thesis of the third-named author \cite{VeselaThesis}.
					\end{remark}
					\begin{proof}[Proof of \Cref{thm:SerreExactSeq}]
						For simplicity denote $A=\pi_2(X)\otimes_{\Lambda}R$. There is the following zig-zag of fibrations
                        \[
                        \begin{tikzcd}
                            B^2A\ar[r]\ar[d]&B^2\pi_2(X)\ar[d]&\wt{X}\ar[d]\ar[l]\\
                            P_2^r(X)\ar[r]\ar[d]&P_2(X)\ar[d]&X\ar[l]\ar[d]\\
                            B\pi\ar[r,equals]&B\pi&B\pi\ar[l,equals]
                        \end{tikzcd}
                        \]
                        where the leftmost fibration is the principal fibration given by the reduced $k$-invariant $k^r\in H^3(\pi;A)$. Consider the third differential of each of the following fibrations. Since we have isomorphisms $H_0(\pi;H_2(A))\xrightarrow{\cong} H_0(\pi;H_2(B^2\pi_2(X)))\xleftarrow{\cong} H_0(\pi;H_2(\wt{X}))\cong A$ the zig-zag above shows that the third differentials of the above fibrations are equal.
                        
                        Next we consider only the left-most fibration. It has the trivial action of $\pi$ on the $\pi_2$ of the fibre and so it can be delooped into a fibration
						\[
                        P_2^r(X)\to B\pi\xrightarrow{k^r}B^3A.
                        \]
						In particular the above composition is null-homotopic, and hence the following
						\begin{equation*}
							\begin{tikzcd}
								B^{2}A\ar[r]\ar[d,equals]&P_2^r(X)\ar[d,"*"]\ar[r]&B\pi\ar[d,"k^r"]\\
								B^2 A\ar[r]&EB^3A\ar[r]&B^{3}A,
							\end{tikzcd}
						\end{equation*}
                        is a map of fibrations, where $EB^3\simeq *$ and  $B^2A\to EB^2A\to B^3A$ is the universal fibration.  Hence, by naturality of the differentials of the Serre spectral sequence, we obtain a commutative diagram
						\[\begin{tikzcd}
							H_3(B\pi)\ar[r,"d_3"]\ar[d,"\ev(k^r)"]& A\ar[d,"\id"]\\
							A\ar[r,"\id"]&A,
						\end{tikzcd}\]
						and so $\ev(k^r)=d_3$.
					\end{proof}

					In the present article, we need various equivalent definitions of $\ev(k^r)$.  We now define maps $\mathfrak{f}_1,\mathfrak{f}_2,\mathfrak{f}_2\colon H_3(B\pi)\to \pi_2(X)\otimes_{\Lambda}\Z$ and show that they are equal to $\ev(k^r)$. 
                    
					\begin{definition}[$\mathfrak{f}_1$]\label{def:mathfrakf1}
						
						Consider the following two $\Lambda$-resolutions of $\Z$, the lower of which is a free resolution.
						\[\begin{tikzcd}
							0\ar[r]&\pi_2(X)\ar[r]&C_2(\wt{X})\ar[r,"d_2^{\wt{X}}"]& C_1(\wt{X})\ar[r,"d_1^{X}"]&C_0(\wt{X})\ar[r,"\varepsilon"]& \Z\ar[r]&0\\
							C_4(\wt{B\pi})\ar[r]\ar[u,dashed]&C_3(\wt{B\pi})\ar[u,dashed,"\mathfrak{f}_1^{\text{chain}}"]\ar[r,"d_2^{\pi}"]&C_2(\wt{B\pi})\ar[r,"d_1^{\pi}"]\ar[u,dashed]&C_1(\wt{B\pi})\ar[u,dashed]\ar[r,"d_0^{\pi}"]&C_0(\wt{B\pi})\ar[r,"\varepsilon"]\ar[u,"i_*",dashed]&\Z\ar[r]\ar[u]&0
						\end{tikzcd}\]
						By the fundamental theorem of homological algebra we can lift the identity map $\Z\to\Z$ to get the vertical maps. The map $\mathfrak{f}_1^{\text{chain}}$ gives an element in the cohomology group $ H^3(\pi;\pi_2(X))$. We define $\mathfrak{f}_1$ as the image of $\mathfrak{f}^{\text{chain}}_1$ under the following composition
						\[H^3(\pi;\pi_2(X))\to H^3(\pi;\pi_2(X)\otimes_{\Lambda}R)\to \Hom(H_3(B\pi),\pi_2(X)\otimes_{\Lambda}R).\]
					\end{definition}

					\begin{definition}[$\mathfrak{f}_2$]\label{def:mathfrakf2}
						Fix a representative of $B\pi$ to be $X$ with cells of dimension 3,4, \dots attached. Define $\mathfrak{f}_2\colon H_3(\pi;R)\to \pi_2(X)\otimes_{\Lambda}R$  by sending $\sum_i r_ie_3^i$ to $\sum_i([\partial e_3^i])\otimes r_i)\in \pi_2(X)\otimes_{\Lambda}R$, where $\partial e_3^i$ can be viewed as an element in $\pi_2(X)$ after choosing a whisker to the basepoint. Since we land in $\pi_2(X)\otimes_{\Lambda}R$, the value of $\mathfrak{f}_2(\sum_i r_ie_3^i)$ does not depend on these choices.
					\end{definition}
					
					\begin{definition}[$\mathfrak{f}_3$]\label{def:geomKinv}
						Fix a representative of $B\pi$ to be $X$ with cells of dimension $3,4,\dots$ attached. Define a map $\mathfrak{f}_3\colon H_3(\pi;R)\to \pi_2(X)\otimes_{\Lambda}R$ as follows. Take a class $z\in H_3(\pi;R)$. Thom showed \cite{Thom1954} that every integral class up to dimension $6$ can be represented by a map from an $R$-oriented manifold. Represent $z$ by a $3$-dimensional manifold $z\colon W_z\to B\pi$. Pick some triangulation of $W_z$ and denote $W^{\circ}_z$ to be $W_z$ with the interior of any $3$-simplex removed. Construct a CW-complex $C_z$ with only $0,1,2$-cells with a homotopy equivalence to $h\colon C_z\to W_z^{\circ}$. Homotope the map $z\circ h\colon C_z\to B\pi$ to land in $X$.  Pick a homotopy inverse $h^{-1}:W_z^{\circ}\to C_z$ a and obtain a map $W_z^{\circ}\to B\pi$, homotopic to $z|_{W_z^{\circ}}$, which lands in $X$. Define $\mathfrak{f}_3$ by sending $z$ to $[\partial W_z^{\circ}]\in \pi_2(X)\otimes_{\Lambda}R$.
					\end{definition}
					
					\begin{prop}\label{prop:f_all_Agree}
						We have equalities $\ev(k^r)=\mathfrak{f}_1=\mathfrak{f}_2=\mathfrak{f}_3$ and in particular the  maps $\mathfrak{f}_1,\mathfrak{f}_2$ and $\mathfrak{f}_3$ are well-defined.
					\end{prop}
					
					\begin{proof}
						The equivalence $\ev(k^r)=\mathfrak{f}_1$ is by \cite{Eilenberg_MacLane_1949}, which even equates the cohomology classes in $H^3(\pi;\pi_2(X))$ used in the respective definitions. For a modern treatment see \cite{CK2025RelativekInva}. Thus $\mathfrak{f}_1$ is well-defined.
						
						Next we establish $\mathfrak{f}_2=\mathfrak{f}_1$. The definition of $\mathfrak{f}_1$ does not depend on the model of $B\pi$, so choose a model for the definition of $\mathfrak{f}_2$. Define a map $\mathfrak{f}^{\text{cell}}\colon C_3(\wt{B\pi})\to \pi_2(X)$ by linearly extending the map defined on a $3$-cell $e_3$ by taking the attaching map $\partial e_3$, which lands in $\wt{X}$. The map $\mathfrak{f}^{\text{cell}}$ fits into the following commutative diagram
						\[
                        \begin{tikzcd}
							0\ar[r]&\pi_2(X)\ar[r,"d_2^{\wt{X}}"]& C^{\text{cell}}_2(\wt{X})\\
							C_4^{\text{cell}}(\wt{B\pi})\ar[u]\ar[r,"d_3^{\pi}"]&C^{\text{cell}}_3(\wt{B\pi})\ar[u,"\mathfrak{f}^{\text{cell}}"]\ar[r,"d_2^{\pi}"]&C^{\text{cell}}_2(\wt{B\pi})\ar[u,equals],
						\end{tikzcd}
                        \]
						and therefore defines a cohomology class $[\mathfrak{f}^{\text{cell}}]\in H^3(\pi;\pi_2(X))$. Since $\mathfrak{f}_1$ is well-defined, the map $\mathfrak{f}^{\text{cell}}$ could have been used as $\mathfrak{f}^{\text{chain}}_1$ in \Cref{def:mathfrakf1} and so $\mathfrak{f}^{\text{cell}}$ maps to $\mathfrak{f}_1$ under the composition
						\[
                        H^3(\pi;\pi_2(X))\to H^3(\pi;\pi_2(X)\otimes_{\Lambda}R)\to \Hom(H_3(B\pi),\pi_2(X)\otimes_{\Lambda}R).
                        \]
						
						By construction $\mathfrak{f}$ also maps to $\mathfrak{f}_2$ showing $\mathfrak{f}_1=\mathfrak{f}_2$.

						Finally we establish that $\mathfrak{f}_2=\mathfrak{f}_3$.  Let $B\pi$ be given as above. Let $z\in H_3(\pi;R)$ be any element represented by a map from an $R$-oriented $3$-manifold $z\colon W_z\to B\pi$. For any given triangulation of $W_z$ approximate $z$ by a cellular map such that all but a single $3$-cell map to $X$. The $R$-fundamental class $ [W_z]$ is given by the sum of all $3$-cells in $W_z$ and it follows that $z$ can be represented by an $R$-linear combination of $3$-cells all but one of which lie in $X$. It is now immediate that $\mathfrak{f}_2(z)=\mathfrak{f}_3(z)$ since all boundaries of $3$-cells in $X$ are nullhomotopic.

					\end{proof}

\section{James spectral sequence}\label{sec:jss}
				
				Defined in the second-named author's thesis, the James spectral sequence (henceforth JSS) \cite[Theorem 3.1.1]{teichnerthesis} is used to compute the bordism groups $\Omega_n^{\xi}$ when $w_1(\xi)=0$.  In the present paper we work with potentially non-orientable $\xi$-structures and so it is desirable to have a more general JSS.  The following construction is a slight generalisation of the construction given in loc.\ cit.\ 

                \begin{remark}\label{rmk:topJSS}
                    Although in this paper we work entirely in the smooth category (see \Cref{rmk:smooth_vs_top}) the James spectral sequence as we present it also applies in the topological category and has uses there, for example when working with the 2-type of a topological 4-manifold (like in \cite{OP22}).  To obtain the topological version of the discussion below, replace all stable vector bundles with stable $\Top$ bundles.  There are appropriate versions of Thom spaces/spectra for $\Top$ bundles (see \cite[IV 5.12]{Rudyak}) and these should be substituted in also.  All of the below arguments work verbatim with these changes made.
                \end{remark}
				
				\begin{definition}
					Let $\xi_n:B\rightarrow BO_n$ be an $n$-dimensional vector bundle. Then we define the {\it Thom space} of $\xi_n$ to be the space $\Th(\xi_n):=D\xi_n/S\xi_n$, where $D\xi_n$ resp.\ $S\xi_n$ are the total spaces of the disc resp.\ sphere bundle over $B$ given by $\xi_n$.
				\end{definition}
				
				By spectra, we mean sequential spectra, also known as prespectra. A sequential spectum is a sequence of spaces $E_1,E_2,\dots$ indexed by $n\in \Z_{>0}$  and bonding maps $\Sigma E_n\to E_{n+1}$. For the following definition of Thom spectra see for example
				\cite[Section 2.1.4]{MalkSpectra}.
				
				\begin{definition}
					Let $\xi\colon B\to BO$ be a stable vector bundle. Define the {\it Thom spectrum of $\xi$} as 
					\[
					M\xi:=(\Th(\xi_n),\Sigma \Th(\xi_n)\to \Th(\xi_{n+1}))
					\]
					where $\xi_n$ is given by the homotopy pullback
					\[
					\begin{tikzcd}
						B_n\ar[r]\ar[d,"\xi_n"]&B\ar[d,"\xi"]\\
						BO_n\ar[r]& BO
					\end{tikzcd}
					\]
					and the bonding map $\Sigma \Th(\xi_n)\to \Th(\xi_{n+1}))$ is the usual map given by the homotopy equivalence $\Th(\xi_N\oplus \varepsilon_1)\cong \Sigma \Th(\xi_n)$ and the map $\xi_n\oplus\varepsilon_1\to \xi_{n+1}$ over $B_n\to B_{n+1}$, where $\varepsilon_1$ is the trivial $1$-dimensional vector bundle.\footnote{
						Note that the spectrum $MO$ arises in this way by taking $\xi=\Id\colon BO\to BO$.}

				\end{definition}
				
				Let $F\to B\xrightarrow{f} K$ be a fibration, and let $\xi\colon B\to BO$ be a stable vector bundle. We choose a representative of the map $BO_n\to BO$ such that it is a fibration which we denote by $BO_n\xrightarrow{\text{fib.}} BO$. It now follows that $B_n\xrightarrow{\text{fib.}}B$ is a fibration. Let $\xi_n:B_n\to BO_n$ be the pullback. Define $F_n$ to be the fibre of the composition $B_n\xrightarrow{\text{fib.}} B\xrightarrow{\text{fib.}} K$.

				\begin{lemma}\label{lem:relativeFibration}
					
					Let $\zeta_k$ be a $k$-dimensional vector bundle over $B_n$. Then the composition $D(\zeta_k)\to B_n\to K$ is a fibration with fibre $D(\zeta_k\vert_{F_n})$ and the following is a relative fibration 
					\begin{equation}\label{eq:monodromy}
						\begin{tikzcd}(D(\zeta_k\vert_{F_n}),S(\zeta_k\vert_{F_n}))\ar[r]&(D(\zeta_k),S(\zeta_k))\ar[r,"\text{fib.}"]& K.
						\end{tikzcd}
					\end{equation}
					
				\end{lemma}
				\begin{proof}
					The first claim follows from the following commutative diagram, where the leftmost square is a pullback
					\[\begin{tikzcd}
						F_n\ar[r,hook]&B_n\ar[r,"\text{fib.}"]&B\ar[r,"\text{fib.}"]&K\\
						D(\zeta_k\vert_{F_n})\ar[r,hook]\ar[u,"\text{fib.}"]&D(\zeta_k)\ar[u,"\text{fib.}"]\ar[rru,"\text{fib.}"', bend right]
					\end{tikzcd}\]
					and from the fact that a composition of fibrations is itself a fibration.  The second fact follows since for any vector bundle the inclusion of the sphere bundle into the disc bundle is a cofibration.
				\end{proof}
				
				For the remainder of the section let $h$ be a generalised homology theory with $\pi_i(h)=0$ for all $i<0$.
				
				\begin{definition}\label{def:monodromy}
					
					Define  $\rho_n\colon\pi_1(K)\to \Aut(h_*(Th(\xi_n\vert_{F_n})))$ to be the monodromy (on the level of $h$-homology) of the relative fibration of \Cref{lem:relativeFibration} for $\zeta_k:=\xi_n$. Similarly, define $\rho_n^{+}\colon \pi_1(K)\to \Aut(h_*(\Th(\xi_n\oplus \varepsilon))$ to be the monodromy corresponding to the fibration of \Cref{lem:relativeFibration} for $\zeta_k:=\xi_n\oplus \varepsilon$. 
				\end{definition}

				\begin{lemma}\label{lem:monodromy}
					There is a natural stabilisation map $g_n\colon \Th(\xi_n\vert_{F_n})\to \Th(\xi_{n}\oplus \varepsilon\vert_{F_{n+1}})$  on Thom spaces which fits into the following commutative diagram for any $\gamma\in\pi_1(K)$.
					\begin{equation}\label{eq:wantLimitOfRhos}
						\begin{tikzcd}
							h_{*}\Th(\xi_n\vert_{F_n}) \ar[r,"\cong"']\arrow[r,"s_n"] \arrow[d,"\rho_n(\gamma)"'] & h_{*+1}\Th(\xi_n\oplus \varepsilon \vert_{F_n})  \arrow[r,"(g_{n})_{*}"] \arrow[d,"\rho^{+}_n(\gamma)"] & h_{*+1}\Th(\xi_{n+1}\vert_{F_{n+1}}) \arrow[d,"\rho_{n+1}(\gamma)"] \\
							h_{*}\Th(\xi_n\vert_{F_n}) \arrow[r,"s_n"]\ar[r,"\cong"'] & h_{*+1}\Th(\xi_n\oplus \varepsilon \vert_{F_n})  \arrow[r,"(g_{n})_{*}"'] & h_{*+1}\Th(\xi_{n+1}\vert_{F_{n+1}})
						\end{tikzcd}    
					\end{equation}
					
					Furthermore, there exists a limit of the maps $\rho_n$ from \Cref{def:monodromy} giving a map \[\hat{\rho}\colon \pi_1(K)\to Aut(h_*M(\xi\vert_F)).\]
				\end{lemma}
				\begin{proof}

					First we construct the map $g_n$. In the commutative diagram
					\[
					\begin{tikzcd}
						F_{n} \arrow[d,"{g_n}",dashed] \arrow[r] & B_{n} \arrow[d] \arrow[r,"\text{fib.}"] & K \arrow[d,"="'] \\
						F_{n+1}  \arrow[r] & B_{n+1}  \arrow[r,"\text{fib.}"] & K
					\end{tikzcd}
					\]
					we see that the map $F_n\to B_n\to B_{n+1}\to K$ is null-homotopic and so there exists a map $g_n$ fitting into the diagram as above.
                    
					We get an induced a map between two fibrations of (\ref{eq:monodromy}) 
					\[
					\begin{tikzcd}(D(\xi_n\oplus\varepsilon\vert_{F_n}),S(\xi_n\oplus\varepsilon\vert_{F_n}))\ar[d,"(g_n)_*"]\ar[r]&(D(\xi_n\oplus\varepsilon),S(\xi_n\oplus\varepsilon))\ar[d]\ar[r]& K\ar[d]\\(D(\xi_{n+1}\vert F_{n+1}),S(\xi_{n+1}\vert_{F_{n+1}}))\ar[r]&(D(\xi_{n+1}),S(\xi_{n+1}))\ar[r]& K.
					\end{tikzcd}
					\]
					This in turn induces a map between the corresponding Thom spaces as claimed and shows that the right square of the diagram (\ref{eq:wantLimitOfRhos}) commutes. The left square of (\ref{eq:wantLimitOfRhos}) commutes since, by construction, $\rho^{+}_n$ acts trivially  the $\varepsilon$-summand of the vector bundle.  This concludes the proof of the commutativity of (\ref{eq:wantLimitOfRhos}). The  outer rectangle of the following diagram 
					\[
					\begin{tikzcd}
						F\ar[r]\ar[d]&B\ar[r]\ar[d]&BO\ar[d]\\
						F_n\ar[r]&B_n\ar[r]& BO_n
					\end{tikzcd}
					\]
					is a pullback. This can be seen since the right square is a pullback and the left square is a pullback as can be shown using the universal property, using that $F$ and $F_n$ respectively are fibres of $B\xrightarrow{\text{fib.}} K$ and $B_n\xrightarrow{\text{fib.}} K$.
					
					Let $M\xi\vert_F$ be the Thom spectrum of $\xi\vert_F$, i.e.\ $M\xi\vert_F$ is the sequential spectrum  \[(\Th(\xi_n\vert_F),\Sigma \Th(\xi_n\vert_F)\to \Th(\xi_{n+1}\vert_F))\] with bonding maps given by $\Sigma \Th(\xi_n\vert_F)\xrightarrow{s_n,\simeq} \Th(\xi_n\oplus \varepsilon\vert_F)\xrightarrow{g_n} \Th(\xi_{n+1}\vert_F)$. In particular, by definition of homology of a spectrum, we have the following colimit
					\[\colim_{n\to\infty}h_{n+k} \Th\xi_n\vert_{F_n}=h_k(M\xi\vert_F).\]
					We use this to show that the limit of maps $\rho_n(\gamma)$ exists for $\gamma\in\pi_1(K)$. The above colimit can be expressed as the quotient  $\sqcup h_{n+k}\Th(\xi_n\vert_{F_n})/ \sim$, where $\sim$ is defined by identifying those elements that eventually share an image. This way we can define 
					\[\hat{\rho}(\gamma)\in \Aut(h_*M(\xi\vert_F)),\] for any element $x\in h_*M(\xi\vert_F)$ by choosing a preimage $x_n\in h_{n+k}\Th(\xi_n\vert_{F_n})$ and defining $\hat{\rho}(\gamma)(x)$ to be the image of $\rho_n(\gamma)(x_n)$.
				\end{proof}
				
				\begin{theorem}[James spectral sequence]\label{thm:jss}
					Let $h$ be a generalised homology theory with $\pi_i(h)=0$ for all $i<0$, let \[F\to B\xrightarrow{f} K\] be a fibration, and let $\xi\colon B\to BO$ be a stable vector bundle.  Denote by $\hat{\rho}\colon \pi_1(K)\to\Aut(h_*(M(\xi\vert_F)))$ the induced monodromy from \Cref{lem:monodromy}.   Then there exists a spectral sequence
					\[
					E^{2}_{p,q}\cong H_p(K; h_q(M(\xi\vert F))_{\hat{\rho}})\,  \implies\, h_{p+q}(M\xi).
					\]
					
				\end{theorem}

				\begin{proof}
					Let $F_n \to B_n\xrightarrow{f} K$ for $n\in\mathbb{N}$ be the sequence of fibrations determined by the standard inclusions $BO_n\to BO$, where $F_n$ is defined to be the fibre (as above, we arrange $BO_n\to BO$ to be fibrations).  This yields a sequence of relative fibrations (\ref{eq:monodromy})
					\[\begin{tikzcd}(D(\xi_n\vert_{F_n}),S(\xi_n\vert_{F_n}))\ar[r]&(D(\xi_n),S(\xi_n))\ar[r,"\text{fib.}"]& K.
					\end{tikzcd}\]
					 with monodromies $\rho_n$.  Accordingly there is a corresponding sequence of relative, twisted Serre spectral sequences (see \cite[Theorem 9.34]{davis_kirk_2001})
					\[
					\prescript{n}{}{E}^{2}_{p,q}\cong H_p(K;h^{n}_q(D(\xi_n\vert_{F_n}),S(\xi_n\vert_{F_n}))_{\rho_n}) \implies h^{n}_{p+q}(D(\xi_n),S(\xi_n)) \cong h_{p+q}(M\xi_n\vert_{F_n}),
					\]
					with differentials denoted by $\prescript{n}{}{d}^k_{p,q}$.  The last isomorphism above comes from the Thom-isomorphism theorem, and we have used the homology theory $h^{n}$ defined by $h_{*}^{n}=h_{*+n}$.\footnote{Note that there is an oversight in the argument in \cite{teichnerthesis} at this point: it is not known if the monodromies $\rho_n$ are trivial even if the monodromy $\hat{\rho}$ is. Our argument fixes this by always allowing our Serre spectral sequences to be twisted.}
					
					We now wish to obtain the JSS as a direct limit of these relative, twisted Serre spectral sequences.  Recall the composition of maps
					\[
                    h_{*}\Th(\xi_n\vert_{F_n}) \xrightarrow{s_n}  h_{*+1}\Th(\xi_n\oplus \varepsilon \vert_{F_n})  \xrightarrow{(g_{n})_{*}} h_{*+1}\Th(\xi_{n+1}\vert_{F_{n+1}}) 
                    \]
					where $s_n$ is the suspension isomorphism and $(g_n)_*$ was defined in \Cref{lem:monodromy}. 
                    
					Using the same lemma, this composition induces a homomorphism of spectral sequences, which we denote by $\prescript{n}{}{\varphi}^k_{p,q}$. For example on the second page $\prescript{n}{}{\varphi}^2_{p,q}$ is given by the following composition
					\begin{align*}
						\prescript{n}{}{E}^{2}_{p,q} &\cong H_p(K;h^{n}_q(D(\xi_n\vert_{F_n}),S(\xi_n\vert_{F_n}))_{\rho_n})\\
						&\cong H_p(K;\wt{h}^{n}_q(\Th(\xi_n\vert_{F_n}))_{\rho_n})\xrightarrow{(s_n)_*,\cong} H_p(K;\wt{h}^{n}_{q+1}(\Th(\xi_n\oplus \varepsilon \vert_{F_n}))_{\rho^{+}_n})\\
						&\cong H_p(K;\wt{h}^{n+1}_{q}(\Th(\xi_n\oplus \varepsilon \vert_{F_n}))_{\rho^+_n})\xrightarrow{(g_n)_*}H_p(K;\wt{h}^{n+1}_q(\Th(\xi_{n+1}\vert_{F_{n+1}}))_{\rho_{n+1}})\cong \prescript{n+1}{}{E}^{2}_{p,q}           
					\end{align*}
					 where the first isomorphism of the third line is given by reindexing $h_q^{1}=h_{q+1}$. By the naturality of the Serre spectral sequence we obtain a commutative diagram involving the differentials
					\[
					\begin{tikzcd}[column sep=large, row sep=large]
						\prescript{n}{}{E}^{k}_{p,q} \arrow[r,"\prescript{n}{}{d}^k_{p,q}"] \arrow[d,"\prescript{n}{}{\varphi}^k_{p,q}"'] & \prescript{n}{}{E}^{k}_{p-k,q+k-1} \arrow[d,"\prescript{n}{}{\varphi}^k_{p-k,q+k-1}"] \\
						\prescript{n+1}{}{E}^{k}_{p,q} \arrow[r,"\prescript{n+1}{}{d}^k_{p,q}"'] & \prescript{n+1}{}{E}^{k}_{p-k,q+k-1}.
					\end{tikzcd}
					\]
					
					Since direct limit is an exact functor the colimit 
					\[
					(E^k_{p,q},d^k_{p,q}):= \colim_{n\to\infty} (\prescript{n}{}{E}^k_{p,q}, d^k_{p,q}).
					\]
					preserves differentials and gives that our colimit indeed defines a spectral sequence.
					
					The calculation of the second page is exactly as in \cite{teichnerthesis} and we obtain
					\[
					E^2_{p,q} = \colim_{n\to\infty} H_p(K; h_q(T(\xi_n\vert_{F_n}))_{\rho_n}) \cong H_p(K;h_q(M\xi\vert_F)_{\hat{\rho}}).\qedhere
					\]
				\end{proof}
                				
				In all applications we will take $h=\pi^{s}$ (i.e.\ stable homotopy theory), the stable vector bundle $\xi=\xi(\pi,w_1^{\pi},w_2^{\pi})$ will be a normal 1-type and we will take as input the fibration
				\begin{equation}\label{eq:JSSfibration}
				\BSpin \to B \to B\pi,
				\end{equation}

                that gives our normal $1$-type.

                We now express the monodromy $\hat{\rho}$.

                \begin{lem}\label{lem:twistingJSS}
         
                    Applying the JSS to the fibration (\ref{eq:JSSfibration}) using $h=\pi^s$ yields a monodromy $\hat{\rho}\colon \pi_1(K)\to \Aut(\Omega^{\Spin}_i)$ given in the low dimensions by
                     \[
                     \begin{cases}
                      \hat{\rho}(\gamma)=\Id& \text{for } i=1,2,3\\
                      \hat{\rho}(\gamma)=w_1(\gamma)& \text{for } i=0,4.
                    \end{cases}
                    \]
                \end{lem}
                \begin{proof}

                Any $\Z/2$-action on $\Z/2\cong \Omega_1^{\Spin}\cong\Omega_2^{\Spin}$ is trivial.

                For $i=0,4$ compare the fibration $B\Spin\to B\to K$ to the fibration $BSO\to BO\to K(\Z/2,1)$. This gives an isomorphism of the coefficients $\Omega_i^{\Spin}\to \Omega^{SO}_i$ for $i=0,4$. The monodromy on $BSO\to BO\to K(\Z/2,1)$ is given by the orientation reversal. It follows from considering the signature that orientation reversal gives the nontrivial involution on $\Omega^{SO}_0$ and $\Omega_4^{SO}$ and naturality concludes the proof. 
                \end{proof}
                
				Applying the JSS to this data gives
				\[
				E^{2}_{p,q}\cong H_p(\pi;\pi^s_q(M(\xi\vert \BSpin))_{\hat{\rho}})\cong H_p(\pi;(\Omega^{\Spin}_q)_{w_1^{\pi}}) \ \implies\ \pi^s_{p+q}(M\xi)\cong \Omega^{\xi}_{p+q}.
				\]

				\subsection{The second differential of the JSS}
				
				We now work to prove the following result concerning the second differential $d_2$ in the JSS for the 1-type $\xi(\pi,w_1^{\pi},w_2^{\pi})$.  In what follows the generalised homology theory $h$ used for any JSS is always assumed to be stable homotopy, i.e.\ $h=\pi^s$.
				
				\begin{proposition}[cf.\ {\cite[Theorem 3.1.3]{teichnerthesis}},  {\cite[Proposition 5.9]{OP22}}]\label{prop:JSSdifferential}
					The second differential $d_2$ in the JSS for the 1-type $\xi(\pi,w_1^{\pi},w_2^{\pi})$ is given as
					\[
					d_2=(\Sq^2(-)+\Sq^1(-)\smile w_1+-\smile w_2)^*\colon H_r(\pi;\Omega_1^{\Spin})\to H_{r-2}(\pi;\Omega_2^{\Spin})
					\]
					for  $r\leq 4$ and
					\[
					d_2=(\Sq^2(-)+\Sq^1(-)\smile w_1+-\smile w_2)^*\circ \red_2\colon H_r(\pi;\Omega_0^{\Spin})\to H_{r-2}(\pi;\Omega_1^{\Spin})
					\]
					for $r\leq 5$.
				\end{proposition}

                To prove this we will follow Orson-Powell \cite[Section 5]{OP22}\footnote{A notational difference is that we will write $\xi$ for vector bundles whereas Orson-Powell write $v$.} where they proved a similar result (cited above), expanding on the work of the second-named author also cited above.  The extra details will make it easier to pinpoint the alterations necessary for the non-orientable case.

                \begin{lem}[{Non-orientable version of \cite[Lemma 5.3]{OP22}}]\label{lem:non-orientable_OP_53}
                    Let $\xi\colon B\to BO$ be a (potentially non-orientable) stable vector bundle over a CW complex.  Then
                    \[
				\Sq^{2}\phi^{\infty}(-)=\phi^{\infty}(\Sq^2(-)+w_2(\xi)\smile-+w_1(\xi)\smile \Sq^1(-)),\]
                for $\phi^{\infty}$ the stable Thom isomorphism map.
                \end{lem}

                \begin{proof}
                    The proof proceeds exactly as in \cite[Proof of Lemma 5.3]{OP22} except instead of the term (in their notation) $U_n\smile \rho^*_n(w_1(\xi_n))\smile \Sq^1(\rho_n^*(x))$ vanishing this now gives precisely the additional term $\Sq^1(-)\smile w_1(\xi)$.
                \end{proof}

                To ease notation further we will now denote by 
                \[
                \Sq^2_{w_1(\xi),w_2(\xi)}(-):=\Sq^2(-)+w_2(\xi)\smile-+w_1(\xi)\smile \Sq^1(-).
                \]
                The next step is compute the differentials for a certain degenerate version of the JSS, which we state and prove now.

                \begin{lemma}[{Non-orientable version of \cite[Lemma 5.4]{OP22}}]\label{lem:differentialsOP}
                 Let $B$ be a CW-complex and consider the input for the James spectral sequence 
                 \[
                \begin{aligned}
                    *\to B\xrightarrow{\id} B
                    \qquad &\text{and}\qquad
                    \xi\colon B\to BO
                \end{aligned}
                \]
                Then the James spectral sequence has the following differentials $d_2^{r,s}\colon E_2^{r,s}\to E_2^{r-2,s+1}$ for all $r\geq 2$.
            \begin{enumerate}
                \item Via the universal coefficient theorem, $d^{r,1}_2\colon H_r(B;\Z/2)\to H_{r-2}(B;\Z/2)$ is identified\footnote{In this section only the notation ${(-)}^{\vee}$ denotes the functor $\Hom(-;\Z/2)$.} with $(\Sq^2_{w_1(\xi),w_2(\xi)})^\vee$. More precisely, the following square commutes:
                \[
                \begin{tikzcd}
                H_r(B;\Z/2)  \ar[d,"\cong"]\ar[rr,"d_2^{r,1}"]
                && H_{r-2}(B;\Z/2) \ar[d,"\cong"] \\
                H^r(B;\Z/2)^\vee \ar[rr,"(\Sq^2_{w_1(\xi),w_2(\xi)})^\vee"]
                && H^{r-2}(B;\Z/2)^\vee.
                \end{tikzcd}
                \]
                \item The map $d_2^{r,0}\colon H_r(B;\Z^{w_1(\xi)})\to H_{r-2}(B;\Z/2)$ is identified with reduction mod 2 followed by $(\Sq^2_{w_1(\xi),w_2(\xi)})^\vee$. More precisely, the following diagram commutes:
                \[
                \begin{tikzcd}
                H_r(B;\Z^{w_1(\xi)}) \ar[rrrr,"d_2^{r,0}"] \ar[d,"\operatorname{red}_2"] &&&& H_{r-2}(B;\Z/2) \ar[d,"\cong"] \\
                H_r(B;\Z/2)\ar[rr,"\cong"] &&H^r(B;\Z/2)^\vee \ar[rr, "(\Sq^2_{w_1(\xi),w_2(\xi)})^\vee"] && H^{r-2}(B;\Z/2)^\vee.
                \end{tikzcd}
                \]
            \end{enumerate}
        \end{lemma}

        \begin{proof}
            We follow \cite[Proof of Lemma 5.4]{OP22}.  First note that $\phi_{\infty}$ still induces a natural isomorphism between the JSS for the degenerate fibration listed above and the Atiyah-Hirzebruch spectral sequence for the Thom spectrum $M(\xi)$.  In the latter the differentials are well-known to be given as \begin{align*}
                d_2^{r,1}= (\Sq^2_{w_1(\xi),w_2(\xi)})^{\vee} \qquad &\text{and} \qquad d_2^{r,0}=(\Sq^2_{w_1(\xi),w_2(\xi)})^{\vee}\circ \red_2.
            \end{align*}
            The proof then collapses to showing that a large diagram commutes, which is exactly the same diagram as in \cite[Page 235, Diagram (5)]{OP22} aside from the third square, which is replaced by
            \[
            \begin{tikzcd}
                H^r(M(\xi);\Z/2)^{\vee} \arrow[d, "(\Sq^2)^{\vee}"] \arrow[r,"(\phi^{\infty})^{\vee}"] & H^r(B;\Z/2)^{\vee} \arrow[d, "(\Sq^2_{w_1(\xi),w_2(\xi)})^\vee"] \\
                H^{r-2}(M(\xi);\Z/2)^{\vee} \arrow[r,"(\phi^{\infty})^{\vee}"] & H^{r-2}(B;\Z/2)^{\vee}
            \end{tikzcd}
            \]
            and this square commutes by \Cref{lem:non-orientable_OP_53}; and the first square, which is replaced by
            \[
            \begin{tikzcd}
                H_r(M(\xi);\Z^{w_1(\xi)})\arrow[d, "\red_2"] \arrow[r,"\phi_{\infty}"] & H_r(B;\Z^{w_1(\xi)}) \arrow[d, "\red_2"] \\
                H^{r-2}(M(\xi);\Z/2)\arrow[r,"\phi_{\infty}"] & H^{r-2}(B;\Z/2)
            \end{tikzcd}
            \]
            which also clearly commutes.  All of the other squares commute by the arguments in \cite{OP22} unchanged.
        \end{proof}
				
			With the previous lemmas in hand, we can state and prove a general statement which says that, in situations where the differential can be compared to the degenerate case of the JSS, we can also compute the differential.

            \begin{prop}[{Non-orientable version of \cite[Corollary 5.7]{OP22}}]\label{Prop:OPSurjectiveJSS}
                Let
                \[
                \begin{tikzcd}
                    B\ar[d,"\xi"]\ar[r,"f"]& K\ar[d]\\
                    BO\ar[r]&K(\Z/2,1)\times K(\Z/2,2)
                \end{tikzcd}
                \]
                be a pullback square of fibrations. Suppose there exist classes $w_1\in H^1(K;\Z/2)$ and $w_2\in H^2(K;\Z/2)$ such that $f^*(w_1)=w_1(\xi)$ and $f^*(w_2)=w_2(\xi)$.
                \begin{enumerate}
                \item\label{it:flesh1}Suppose for some $r\in\Z$ that $f_*\colon H_r(B;\Z/2)\to H_r(K;\Z/2)$ is surjective. Then, in the sense of \Cref{lem:differentialsOP}, $d^{r,1}_2\colon H_r(K;\Z/2)\to H_{r-2}(K;\Z/2)$ is $(\Sq^2_{w_1,w_2})^\vee$.
                \item\label{it:flesh2} Suppose for some $r\in\Z$ that $f_*\colon H_r(B;\Z^{w_1(\xi)})\to H_r(K;\Z^{w_1})$ is surjective. Then, in the sense of \Cref{lem:differentialsOP}, $d_2^{r,0}\colon H_r(K;\Z^{w_1})\to H_{r-2}(K;\Z/2)$ is reduction mod 2 followed by $(\Sq^2_{w_1,w_2})^\vee$.
            \end{enumerate}
            \end{prop}
            Note that we have a stronger requirement for the fibration in \Cref{Prop:OPSurjectiveJSS} than in \cite[Corollary 5.7]{OP22}, because we need the twisting in the JSS for $*\to B\to B$ to agree with the twisting in our fibration $F\to B\to K$. With this extra assumption it follows from \Cref{lem:twistingJSS} that this is the case, and hence the domain of the differential $d_0^{r,0}$ is correct as written.
            \begin{proof}
               The additional assumption of \cite[Corollary 5.7]{OP22} is satisfied, since the map $*\to B\Spin\to BO$ induces an isomorphism on the Thom spectra $\pi_i(M(*))\cong \pi_i^{st}\xrightarrow{\cong}\pi_i(M(\xi\circ i))$ for $i=0,1,2$. The proof of (1) is entirely analogous to that in loc.\ cit.; the proof of (2) is left to the reader.  
                
                We give the setup.  By assumption, for any $x\in H_r(K;\Z/2)$ there exists a $y\in H_r(B;Z/2)$ such that $f_*(y)=x$.  The proof then proceeds exactly as in \cite{OP22} except for differences in their final calculation, which we now redo. Let $u\in H^{r-2}(K;\Z/2)$ be arbitrary.  We then have that
            \begin{align*}
                \left((f^*)^\vee\circ \Sq^2_{w_1(\xi),w_2(\xi)}(u)\right)&(y)= \left(\Sq^2_{w_1(\xi),w_2(\xi)}(f^*u)\right)(y)\\
                &= \Sq^2(f^*u)(y)+ \left(w_2(\xi)\cup (f^*u)\right)(y)+\left(w_1(v)\cup\Sq^1(f^*(u))\right)(y)\\
                &= \Sq^2(u)(f_*y)+ \left((f^*w_2)\cup (f^*u)\right)(y)+\left(f^*w_1\cup f^*\Sq^1(u)\right)(y)\\
                &= \Sq^2(u)(f_*y)+ \left(w_2\cup u\right)(f_*y)+\left(w_1\cup \Sq^1(u)\right)(f_*y)\\
                &= \Sq^2(u)(x)+ \left(w_2\cup u+w_1\cup\Sq^1(u)\right)(x)\\
                &=\Sq^2_{w_1,w_2}(x).\qedhere
            \end{align*}
        \end{proof}

        We now finish this section by using \Cref{Prop:OPSurjectiveJSS} to prove \Cref{prop:JSSdifferential}, which was stated earlier.

        \begin{proof}[Proof of \Cref{prop:JSSdifferential}]
            In the proof of \cite[Theorem 3.1.3]{teichnerthesis} the second-named author shows that $H_i(B;\Z/2)\to H_i(\pi;\Z/2)$ is surjective for $i\leq 4$ and $H_i(B;\Z)\to H_i(\pi;\Z)$ is surjective for $i\leq 5$. In the non-orientable case the proof is the same and uses that $H_5(K(\Z/2,1)\times K(\Z/2,2);\Z^{w_1})$ is finite.   We conclude the proof using \Cref{Prop:OPSurjectiveJSS}. 
        \end{proof}

					\bibliography{biblio.bib}{}

@article {norman_69,
    AUTHOR = {Norman, R. A.},
     TITLE = {Dehn's lemma for certain {$4$}-manifolds},
   JOURNAL = {Invent. Math.},
  FJOURNAL = {Inventiones Mathematicae},
    VOLUME = {7},
      YEAR = {1969},
     PAGES = {143--147},
      ISSN = {0020-9910,1432-1297},
   MRCLASS = {57.10 (55.00)},
  MRNUMBER = {246309},
MRREVIEWER = {H.\ Terasaka},
       DOI = {10.1007/BF01389797},
       URL = {https://doi.org/10.1007/BF01389797},
}

@book {akbulut,
    AUTHOR = {Akbulut, Selman},
     TITLE = {4-manifolds},
    SERIES = {Oxford Graduate Texts in Mathematics},
    VOLUME = {25},
 PUBLISHER = {Oxford University Press, Oxford},
      YEAR = {2016},
     PAGES = {xii+262},
      ISBN = {978-0-19-878486-9},
   MRCLASS = {57N13 (57-02)},
  MRNUMBER = {3559604},
MRREVIEWER = {Jonathan\ D.\ Williams},
       DOI = {10.1093/acprof:oso/9780198784869.001.0001},
       URL = {https://doi.org/10.1093/acprof:oso/9780198784869.001.0001},
}

@book {davis_kirk_2001,
    AUTHOR = {Davis, James F. and Kirk, Paul},
     TITLE = {Lecture notes in algebraic topology},
    SERIES = {Graduate Studies in Mathematics},
    VOLUME = {35},
 PUBLISHER = {American Mathematical Society, Providence, RI},
      YEAR = {2001},
     PAGES = {xvi+367},
      ISBN = {0-8218-2160-1},
   MRCLASS = {55-01 (57-01)},
  MRNUMBER = {1841974},
MRREVIEWER = {Don\ Shimamoto},
       DOI = {10.1090/gsm/035},
       URL = {https://doi.org/10.1090/gsm/035},
}

@article {OP22,
    AUTHOR = {Orson, Patrick and Powell, Mark},
     TITLE = {Mapping class groups of simply connected 4-manifolds with
              boundary},
   JOURNAL = {J. Differential Geom.},
  FJOURNAL = {Journal of Differential Geometry},
    VOLUME = {131},
      YEAR = {2025},
    NUMBER = {1},
     PAGES = {199--275},
      ISSN = {0022-040X,1945-743X},
   MRCLASS = {57K40 (57N37 57R52 57S05)},
  MRNUMBER = {4947553},
       DOI = {10.4310/jdg/1755544135},
       URL = {https://doi.org/10.4310/jdg/1755544135},
}

@article {KPT25,
    AUTHOR = {Kasprowski, Daniel and Powell, Mark and Teichner, Peter},
     TITLE = {The {K}ervaire-{M}ilnor invariant in the stable classification
              of spin 4-manifolds},
   JOURNAL = {Tunis. J. Math.},
  FJOURNAL = {Tunisian Journal of Mathematics},
    VOLUME = {7},
      YEAR = {2025},
    NUMBER = {2},
     PAGES = {417--436},
      ISSN = {2576-7658,2576-7666},
   MRCLASS = {57K40},
  MRNUMBER = {4904391},
       DOI = {10.2140/tunis.2025.7.417},
       URL = {https://doi.org/10.2140/tunis.2025.7.417},
}

@article{Kreck99,
	author = {Kreck, M.},
	fjournal = {Annals of Mathematics. Second Series},
	journal = {Ann. of Math. (2)},
	mrreviewer = {Laurence R. Taylor},
	number = {3},
	pages = {707--754},
	title = {Surgery and duality},
	url = {http://dx.doi.org/10.2307/121071},
	volume = {149},
	year = {1999}}

@book{Rudyak,
	author = {Rudyak, Yuli B.},
	isbn = {3-540-62043-5},
	mrclass = {55-02 (55N22 55P42 57-02)},
	mrnumber = {1627486},
	mrreviewer = {Donald M. Davis},
	note = {With a foreword by Haynes Miller},
	pages = {xii+587},
	publisher = {Springer-Verlag, Berlin},
	series = {Springer Monographs in Mathematics},
	title = {On {T}hom spectra, orientability, and cobordism},
	year = {1998}}

@phdthesis{teichnerthesis,
	author = {Teichner, P.},
	note = {Shaker Verlag, ISBN 3-86111-182-9},
	school = {University of Mainz, Germany},
	title = {Topological 4-manifolds with finite fundamental group},
	year = {1992}}

@article{Freedman82,
	author = {Freedman, M. H.},
	fjournal = {Journal of Differential Geometry},
	journal = {J. Differential Geom.},
	mrreviewer = {John J. Walsh},
	number = {3},
	pages = {357--453},
	title = {The topology of four-dimensional manifolds},
	url = {http://projecteuclid.org/euclid.jdg/1214437136},
	volume = {17},
	year = {1982}}

@book{Freedman-Quinn,
	address = {Princeton, NJ},
	author = {Freedman, M. H. and Quinn, F.},
	mrreviewer = {Ian Hambleton},
	pages = {viii+259},
	publisher = {Princeton University Press},
	series = {Princeton Mathematical Series},
	title = {Topology of 4-manifolds},
	volume = {39},
	year = {1990}}

@Article{KPT,
 Author = {Kasprowski, Daniel and Powell, Mark and Teichner, Peter},
 Title = {Four-manifolds up to connected sum with complex projective planes},
 FJournal = {American Journal of Mathematics},
 Journal = {Am. J. Math.},
 ISSN = {0002-9327},
 Volume = {144},
 Number = {1},
 Pages = {75--118},
 Year = {2022},
 DOI = {10.1353/ajm.2022.0001},
 zbMATH = {7465467}
}

@article{Ham-kreck-finite,
	author = {Hambleton, I. and Kreck, M.},
	fjournal = {Mathematische Annalen},
	journal = {Math. Ann.},
	mrreviewer = {Laurence R. Taylor},
	number = {1},
	pages = {85--104},
	title = {On the classification of topological {$4$}-manifolds with finite fundamental group},
	url = {http://dx.doi.org/10.1007/BF01474183},
	volume = {280},
	year = {1988}}

@book{Wall99,
	author = {Wall, C. T. C.},
	edition = {Second},
	note = {Edited and with a foreword by A. A. Ranicki},
	pages = {xvi+302},
	publisher = {American Mathematical Society, Providence, RI},
	series = {Mathematical Surveys and Monographs},
	title = {Surgery on compact manifolds},
	url = {http://dx.doi.org/10.1090/surv/069},
	volume = {69},
	year = {1999}}

@book{Kirby-4-manifold-book,
	author = {Kirby, Robion C.},
	isbn = {3-540-51148-2},
	mrclass = {57N13 (57M25 57R15 57R20 57R75)},
	mrnumber = {1001966 (90j:57012)},
	mrreviewer = {Christopher W. Stark},
	pages = {vi+108},
	publisher = {Springer-Verlag, Berlin},
	series = {Lecture Notes in Mathematics},
	title = {The topology of {$4$}-manifolds},
	volume = {1374},
	year = {1989}}

@article{Milnor-Spin,
	author = {Milnor, John W.},
	journal = {Enseignement Math. (2)},
	pages = {198--203},
	title = {Spin structures on manifolds},
	volume = {9},
	year = {1963}}

@book{milnorstasheff,
  title={Characteristic classes},
  author={Milnor, John and Stasheff, James},
  number={76},
  year={1974},
  publisher={Princeton university press}
}

@article{Kervaire57,
title={Relative Characteristic Classes},
author={Michel Kervaire},
journal={American Journal of Mathematics},
pages={517-558},
volume= {79},
number={3},
year={1957},
month={July},
publisher={ The Johns Hopkins University Press}}

@article{Thom1954,
author = {Thom, Réné},
journal = {Commentarii mathematici Helvetici},
pages = {17-86},
title = {Quelques propriétés globales des variétés différentiables.},
url = {http://eudml.org/doc/139072},
volume = {28},
year = {1954},
}

@article {kasprowski_nicholson_vesela_24,
    AUTHOR = {Kasprowski, Daniel and Nicholson, John and Vesel\'a, Simona},
     TITLE = {Stable equivalence relations on 4-manifolds},
   JOURNAL = {Proc. Lond. Math. Soc. (3)},
  FJOURNAL = {Proceedings of the London Mathematical Society. Third Series},
    VOLUME = {131},
      YEAR = {2025},
    NUMBER = {5},
     PAGES = {Paper No. e70101},
      ISSN = {0024-6115,1460-244X},
   MRCLASS = {57K40 (19J25 57N65 57Q10 57R67)},
  MRNUMBER = {4989093},
       DOI = {10.1112/plms.70101},
       URL = {https://doi.org/10.1112/plms.70101},
}

@article {kosanovic2023newapproachlightbulb,
    AUTHOR = {Kosanovi\'c, Danica and Teichner, Peter},
     TITLE = {A new approach to light bulb tricks: disks in 4-manifolds},
   JOURNAL = {Duke Math. J.},
  FJOURNAL = {Duke Mathematical Journal},
    VOLUME = {173},
      YEAR = {2024},
    NUMBER = {4},
     PAGES = {673--721},
      ISSN = {0012-7094,1547-7398},
   MRCLASS = {57K45 (57R40 58D10)},
  MRNUMBER = {4734552},
MRREVIEWER = {Richard\ Stong},
       DOI = {10.1215/00127094-2023-0036},
       URL = {https://doi.org/10.1215/00127094-2023-0036},
}

@Article{cairns_1940,
 Author = {Cairns, Stewart S.},
 Title = {Homeomorphisms between topological manifolds and analytic manifolds},
 FJournal = {Annals of Mathematics. Second Series},
 Journal = {Ann. Math. (2)},
 ISSN = {0003-486X},
 Volume = {41},
 Pages = {796--808},
 Year = {1940},
 Language = {English},
 DOI = {10.2307/1968860},
 zbMATH = {3040528},
 Zbl = {0025.23502}
}

@Article{moise_1952,
 Author = {Moise, Edwin E.},
 Title = {Affine structures in 3-manifolds. {V}: {The} triangulation theorem and {Hauptvermutung}},
 FJournal = {Annals of Mathematics. Second Series},
 Journal = {Ann. Math. (2)},
 ISSN = {0003-486X},
 Volume = {56},
 Pages = {96--114},
 Year = {1952},
 Language = {English},
 DOI = {10.2307/1969769},
 zbMATH = {3075448},
 Zbl = {0048.17102}
}

@Article{bing_1954,
 Author = {Bing, R. H.},
 Title = {Locally tame sets are tame},
 FJournal = {Annals of Mathematics. Second Series},
 Journal = {Ann. Math. (2)},
 ISSN = {0003-486X},
 Volume = {59},
 Pages = {145--158},
 Year = {1954},
 Language = {English},
 DOI = {10.2307/1969836},
 zbMATH = {3087489},
 Zbl = {0055.16802}
}

@Article{kervaire_milnor,
 Author = {Kervaire, M. A. and Milnor, John W.},
 Title = {On 2-spheres in 4-manifolds},
 FJournal = {Proceedings of the National Academy of Sciences of the United States of America},
 Journal = {Proc. Natl. Acad. Sci. USA},
 ISSN = {0027-8424},
 Volume = {47},
 Pages = {1651--1657},
 Year = {1961},
 Language = {English},
 DOI = {10.1073/pnas.47.10.1651},
 URL = {www.ncbi.nlm.nih.gov/pmc/articles/PMC223187},
 zbMATH = {3175821},
 Zbl = {0107.40303}
}

@Article{tristram,
 Author = {Tristram, A. G.},
 Title = {Some cobordism invariants for links},
 FJournal = {Proceedings of the Cambridge Philosophical Society},
 Journal = {Proc. Camb. Philos. Soc.},
 ISSN = {0008-1981},
 Volume = {66},
 Pages = {251--264},
 Year = {1969},
 Language = {English},
 zbMATH = {3304750},
 Zbl = {0191.54703}
}

@Article{kasprowski_powell_ray_teichner,
 Author = {Kasprowski, Daniel and Powell, Mark and Ray, Arunima and Teichner, Peter},
 Title = {Embedding surfaces in 4-manifolds},
 FJournal = {Geometry \& Topology},
 Journal = {Geom. Topol.},
 ISSN = {1465-3060},
 Volume = {28},
 Number = {5},
 Pages = {2399--2482},
 Year = {2024},
 Language = {English},
 DOI = {10.2140/gt.2024.28.2399},
 zbMATH = {7927933}
}

@article {stong_1994,
    AUTHOR = {Stong, Richard},
     TITLE = {Existence of {$\pi_1$}-negligible embeddings in
              {$4$}-manifolds. {A} correction to {T}heorem 10.5 of
              {F}reedmann and {Q}uinn},
   JOURNAL = {Proc. Amer. Math. Soc.},
  FJOURNAL = {Proceedings of the American Mathematical Society},
    VOLUME = {120},
      YEAR = {1994},
    NUMBER = {4},
     PAGES = {1309--1314},
      ISSN = {0002-9939,1088-6826},
   MRCLASS = {57N35 (57N13 57Q25)},
  MRNUMBER = {1215031},
       DOI = {10.2307/2160253},
       URL = {https://doi.org/10.2307/2160253},
}

@unpublished{MalkSpectra,
    author ={Cary Malkiewich} ,
    title ={Spectra and stable homotopy theory} ,
    note={Lecture Notes},
    year = {retrieved: 08/2025},
    url={https://people.math.binghamton.edu/malkiewich/}
}

@article{Eilenberg_MacLane_1949, 
title={Homology of spaces with operators. II}, 
volume={65}, 
url={http://dx.doi.org/10.1090/s0002-9947-1949-0033001-0}, 
DOI={10.1090/s0002-9947-1949-0033001-0}, 
number={1}, 
journal={Transactions of the American Mathematical Society}, 
publisher={American Mathematical Society (AMS)}, author={Eilenberg, Samuel and MacLane, Saunders}, 
year={1949}, 
pages={49–99}, 
language={en} }

@misc{CK2025RelativekInva,
      title={The first relative k-invariant}, 
      author={Anthony Conway and Daniel Kasprowski},
      note={Preprint, available at \url{https://arxiv.org/abs/2510.18796}},
      year={2025},
      eprint={2510.18796},
      archivePrefix={arXiv},
      primaryClass={math.GT},
      url={https://arxiv.org/abs/2510.18796}, 
}

@misc{VeselaThesis,
    title={The $1.5$-type and surfaces in $4$-manifolds},
    author={Simona Vesel\'a},
    year={2026},
    note={PhD. thesis, (to appear)}
}
					\bibliographystyle{amsalpha}
					
				\end{document}